\documentclass[12pt]{article}
\usepackage{amsmath,amsfonts,amssymb,amsthm}
\input amssym.def
\begin{document}
\def \Z{\Bbb Z}
\def \C{\Bbb C}
\def \R{\Bbb R}
\def \Q{\Bbb Q}
\def \N{\Bbb N}
\def \bR{\bf R}
\def \A{{\mathcal{A}}}
\def \D{{\cal{D}}}
\def \E{{\cal{E}}}
\def \S{{\cal{S}}}
\def \R{{\cal{R}}}
\def \wt{{\rm wt}}
\def \tr{{\rm tr}}
\def \span{{\rm span}}
\def \Res{{\rm Res}}
\def \End{{\rm End}}
\def \Ind {{\rm Ind}}
\def \Irr {{\rm Irr}}
\def \Aut{{\rm Aut}}
\def \Hom{{\rm Hom}}
\def \mod{{\rm mod}}
\def \ann{{\rm Ann}}
\def \ad{{\rm ad}}
\def \rank{{\rm rank}\;}

\def \<{\langle}
\def \>{\rangle}
\def \t{\tau }
\def \a{\alpha }
\def \e{\epsilon }
\def \l{\lambda }
\def \L{\Lambda }
\def \g{{\frak{g}}}
\def \h{{\hbar}}
\def \k{{\frak{k}}}
\def \sl{{\frak{sl}}}
\def \b{\beta }
\def \om{\omega }
\def \o{\omega }
\def \c{\chi}
\def \ch{\chi}
\def \cg{\chi_g}
\def \ag{\alpha_g}
\def \ah{\alpha_h}
\def \ph{\psi_h}
\def \be{\begin{equation}\label}
\def \ee{\end{equation}}
\def \bex{\begin{example}\label}
\def \eex{\end{example}}
\def \bl{\begin{lem}\label}
\def \el{\end{lem}}
\def \bt{\begin{thm}\label}
\def \et{\end{thm}}
\def \bp{\begin{prop}\label}
\def \ep{\end{prop}}
\def \br{\begin{rem}\label}
\def \er{\end{rem}}
\def \bc{\begin{coro}\label}
\def \ec{\end{coro}}
\def \bd{\begin{de}\label}
\def \ed{\end{de}}
\def \pf{{\bf Proof. }}

\newcommand{\n}{\:^{\times}_{\times}\:}
\newcommand{\nno}{\nonumber}
\newcommand{\nord}{\mbox{\scriptsize ${\circ\atop\circ}$}}
\newtheorem{thm}{Theorem}[section]
\newtheorem{prop}[thm]{Proposition}
\newtheorem{coro}[thm]{Corollary}
\newtheorem{conj}[thm]{Conjecture}
\newtheorem{example}[thm]{Example}
\newtheorem{lem}[thm]{Lemma}
\newtheorem{rem}[thm]{Remark}
\newtheorem{de}[thm]{Definition}
\newtheorem{hy}[thm]{Hypothesis}
\makeatletter
\@addtoreset{equation}{section}
\def\theequation{\thesection.\arabic{equation}}
\makeatother
\makeatletter

\begin{center}
{\Large \bf $\hbar$-adic quantum vertex algebras and their modules}
\end{center}

\begin{center}
{Haisheng Li\footnote{Partially supported by NSF grant DMS-0600189}\\
Department of Mathematical Sciences\\
Rutgers University, Camden, NJ 08102}
\end{center}

\begin{abstract}
This is a paper in a series to study vertex algebra-like structures
arising from  various algebras including quantum affine algebras and
Yangians. In this paper, we study notions of $\hbar$-adic nonlocal
vertex algebra and $\hbar$-adic (weak) quantum vertex algebra,
slightly generalizing Etingof-Kazhdan's notion of quantum vertex
operator algebra.  For any topologically free $\C[[\h]]$-module $W$,
we study $\hbar$-adically compatible subsets and $\hbar$-adically
$\S$-local subsets of $(\End W)[[x,x^{-1}]]$. We prove that any
$\hbar$-adically compatible subset generates an $\hbar$-adic
nonlocal vertex algebra with $W$ as a module and that any
$\hbar$-adically $\S$-local subset generates an $\hbar$-adic weak
quantum vertex algebra with $W$ as a module. A general construction
theorem of $\hbar$-adic nonlocal vertex algebras and $\hbar$-adic
quantum vertex algebras is obtained. As an application we associate
the centrally extended double Yangian of $\sl_{2}$ to $\hbar$-adic
quantum vertex algebras.
\end{abstract}

\section{Introduction}
In \cite{ek}, one of an important series of papers, Etingof and
Kazhdan introduced a fundamental notion of quantum vertex operator
algebra and they constructed a family of quantum vertex operator
algebras which are formal deformations of vertex operator algebras
associated with affine Lie algebras $\widehat{\sl_{n+1}}$. For a
quantum vertex operator algebra in this sense, the underlying space
is a topologically free $\C[[\hbar]]$-module $V=V^{0}[[\hbar]]$ with
$V^{0}$ a vector space over $\C$, and the vertex operator map $Y$ is
a $\C[[\hbar]]$-module map from $V$ to $\Hom (V,
V^{0}((x))[[\hbar]])$, where the key axioms are a quasi
commutativity property, called $\S$-locality, and an associativity
property. Furthermore,  the $\S$-locality by assumption is governed
by a unitary rational quantum Yang-Baxter operator $\S(x)$ on $V$.
It follows from the definition that $V/\hbar V$ is an ordinary
vertex algebra (over $\C$), so quantum vertex operator algebras in
this sense are formal deformations of vertex algebras. As it was
mentioned therein, a generalization of this theory to the super case
is straightforward.

Inspired by \cite{ek},  in a series of papers (\cite{li-qva1},
\cite{li-qva2}, \cite{li-infinity}, \cite{kl}) we have extensively
studied a notion of (weak) quantum vertex algebra, as a
generalization of the notions of vertex algebra and vertex
superalgebra. For a (weak) quantum vertex algebra $V$ in this sense,
the underlying space is a vector space over $\C$ and the vertex
operator map $Y$ is a linear map {}from $V$ to $\Hom (V,V((x)))$,
which satisfies a certain braided Jacobi identity, or equivalently
an $\S$-locality and associativity. This theory of (weak) quantum
vertex algebras has many of the features of the theory of ordinary
vertex (super-)algebras. For example, as it was proved in
\cite{li-qva1}, weak quantum vertex algebras and their modules can
be constructed {}from what were called $\S$-local sets of vertex
operators on an arbitrarily given vector space, just as vertex
(super)algebras and modules can be constructed from ``mutually
local'' vertex operators (see \cite{li-local}). Examples of quantum
vertex algebras and modules were constructed in \cite{li-qva2}
{}from Zamolodchikov-Faddeev algebras of a certain type and in
\cite{li-infinity} from $q$-versions of double Yanigan
$DY_{\hbar}(\sl_{2})$ with $q$ a nonzero complex number.

In this paper, we come back to Etingof-Kazhdan's notion of quantum
vertex operator algebra with a slight generalization such that the
classical limits are more general quantum vertex algebras. More
specifically, we systematically study notions of $\hbar$-adic
nonlocal vertex algebra and $\hbar$-adic (weak) quantum vertex
algebra, and we establish general construction theorems, with the
ultimate goal to associate such $\hbar$-adic quantum vertex algebras
to centrally extended double Yangians essentially in the same way
that affine Lie algebras were associated with vertex operator
algebras. An $\hbar$-adic nonlocal vertex algebra will be a
topologically free $\C[[\hbar]]$-module $V$ equipped with a
$\C[[\hbar]]$-module map $Y$ from $V$ to $(\End V)[[x,x^{-1}]]$ and
a vector ${\bf 1}\in V$ such that for every positive integer $n$,
$V/\hbar^{n}V$ is a nonlocal vertex algebra over $\C$, while an
$\hbar$-adic weak quantum vertex algebra is an $\hbar$-adic nonlocal
vertex algebra $V$ that satisfies $\S$-locality in the sense of
\cite{ek} with $\S(x)$ only a $\C[[\hbar]]$-module map without any
other assumption. Furthermore, an $\hbar$-adic quantum vertex
algebra is an $\hbar$-adic weak quantum vertex algebra $V$ such that
the $\S$-locality operator $\S(x)$ is a unitary rational quantum
Yang-Baxter operator and satisfies the shift condition and hexagon
identity as in \cite{ek}.

For each finite-dimensional simple Lie algebra $\g$, Drinfeld (see
\cite{dr1}) introduced a Hopf algebra $Y(\g)$, called Yangian, as a
deformation of the universal enveloping algebra $U(\g[t])$ of Lie
algebra $\g[t]$. Then the double $DY_{\hbar}(\g)$ of $Y(\g)$ in the
sense of Drinfeld was studied in \cite{kt}. Furthermore, centrally
extended double Yangian $\widehat{DY_{\hbar}(\g)}$ was studied in
\cite{kh} and \cite{ioko}, as a deformation of the universal
enveloping algebra $U(\hat{\g})$ of the affine Lie algebra
$\hat{\g}$, where a vertex operator representation was also given.
Our objective is to establish a canonical association of the
centrally extended double Yangians (in which the parameter $\hbar$
is a formal variable, instead of a complex number) with vertex
algebra-like structures. This is our main motivation to study
$\hbar$-adic (weak) quantum vertex algebras.

In this paper we build the foundation for this theory of
$\hbar$-adic (weak) quantum vertex algebras. As one of the main
results, we establish a general construction of $\hbar$-adic weak
quantum vertex algebras and their modules. This is an $\hbar$-adic
version of the general construction in \cite{li-qva2} of weak
quantum vertex algebras and their modules, and we here extensively
use the results therein. Let $W$ be a general $\C[[\hbar]]$-module.
Consider formal series
$$a(x)=\sum_{m\in \Z}a_{m}x^{-m-1}\in (\End W)[[x,x^{-1}]]$$
satisfying the condition that for every $w\in W$ and for every
positive integer $n$, there exists an integer $k$ such that
$a_{m}w\in \hbar^{n}W$ for $m\ge k$, and let $\E_{\hbar}(W)$ consist
of all such $a(x)$. In the case that $W=W^{0}[[\hbar]]$ is
topologically free (with $W^{0}$ a vector space over $\C$), we have
$\E_{\hbar}(W)=\E(W^{0})[[\hbar]]$ which is also topologically free.
(Recall that for a vector space $U$ over $\C$, $\E(U)=\Hom
(U,U((x)))$.) We then study what we call $\hbar$-adically compatible
subsets and $\hbar$-adically $\S$-local subsets of $\E_{\hbar}(W)$.
We prove that any $\hbar$-adically compatible subset of
$\E_{\hbar}(W)$ generates an $\hbar$-adic nonlocal vertex algebra
with $W$ as a canonical module, while an $\hbar$-adically $\S$-local
subset of $\E_{\hbar}(W)$ generates an $\hbar$-adic weak quantum
vertex algebra with $W$ as a canonical module.

The fact is that the generating functions of a centrally extended
double Yangian  on a highest weight module $W$ together with the
identity operator $1_{W}$ form an $\hbar$-adically $\S$-local subset
of $\E_{\hbar}(W)$, and hence one has an $\hbar$-adic weak quantum
vertex algebra generated by those generating functions. It was known
(see \cite{li-local}, \cite{ll}) that if $W$ is a highest weight
module for affine Lie algebra $\hat{\g}$ of level $\ell\in \C$, then
the canonical generating functions of $\hat{\g}$ generate a vertex
algebra, which can be identified as a so-called vacuum
$\hat{\g}$-module of the same level $\ell$. For centrally extended
double Yangians, the situation is different; the generated
$\hbar$-adic weak quantum vertex algebra is not a module for
$\widehat{DY_{\hbar}(\g)}$, but it is a module for a certain cover
of $\widehat{DY_{\hbar}(\g)}$. This is mainly due to that the field
associated to a Cartan element is broken into two fields in the
quantum case. In this paper, we pick up the simplest case with
$\g=\sl_{2}$ and work out the details. More specifically, we
introduce a cover $\widetilde{DY_{\hbar}(\sl_{2})}$ of
$\widehat{DY_{\hbar}(\sl_{2})}$ and by using our general
construction we show that on a universal vacuum
$\widetilde{DY_{\hbar}(\sl_{2})}$-module of a generic level (which
is defined suitably), there exists a canonical $\hbar$-adic quantum
vertex algebra structure with every highest weight
$\widehat{DY_{\hbar}(\sl_{2})}$-module of the same level as a
module. In principle, a generalization to the centrally extended
double Yangian of a general finite-dimensional simple Lie algebras
can be done in a similar way, but one has to deal with the
complicated Serre type relations. We plan to study this in a future
publication.

In this paper, we also construct a family of $\hbar$-adic quantum
vertex algebras as deformations of certain quantum vertex algebras
which were studied in \cite{kl}. Those quantum vertex algebras were
constructed by using certain generalized Weyl-Clifford algebras, or
Zamolodchikov-Faddeev algebras. In one special case, we obtain a
quantum $\beta\gamma$-system, and in another we obtain a formal
deformation of the vertex operator superalgebra $V_{L}$ associated
with the lattice $L=\Z\alpha$ with $\<\alpha,\alpha\>=1$.

There is a very interesting paper \cite{ab}, in which Anguelova and
Bergvelt studied a broad class of vertex algebra-like structures
called $H_{D}$-quantum vertex algebras, using some ideas of
Borcherds from \cite{b-qva}, and they  constructed certain
interesting examples  by employing Borcherds' bicharacter
construction.  This notion of $H_{D}$-quantum vertex algebra
generalizes the notion of braided vertex operator algebra in
\cite{ek} in several aspects. For example, the braiding operator
$\S$ (describing quasi locality) is allowed to have two
(independent) spectral parameters, instead of one. What we here call
$\hbar$-adic quantum vertex algebras can be considered as a
subfamily of $H_{D}$-quantum vertex algebras. A drawback of this
generality is that general $H_{D}$-quantum vertex algebras, just as
Etingof-Kazhdan's braided vertex operator algebras, fail to satisfy
the usual associativity for vertex algebras, though they do satisfy
a braided associativity. On the other hand, weak quantum vertex
algebras in the sense of \cite{li-qva1} and $\hbar$-adic weak
quantum vertex algebras all satisfy the usual associativity, which
promises a transparent representation theory. Especially, examples
of $\hbar$-adic quantum vertex algebras (and their modules) are
constructed by using vertex operators on potential modules from a
representation point of view.

This paper is organized as follows: In Section 2,  we study
$\hbar$-adic nonlocal vertex algebras and $\hbar$-adic (weak)
quantum vertex algebras and present some basic results. In Section
3, we present some technical results. In particular we discuss
$\hbar$-adic nonlocal vertex subalgebras. In Section 4, we give a
general construction of $\hbar$-adic (weak) quantum vertex algebras
and their modules using $\hbar$-adic $\S$-local sets of (formal)
vertex operators. In Section 5, we construct some $\hbar$-adic
quantum vertex algebras which are deformations of certain quantum
vertex algebras. In Section 6, as an application of the general
construction we associate the centrally extended double Yangian of
$\sl_{2}$ with $\hbar$-adic quantum vertex algebras.

Acknowledgement: Part of this paper was finished during my visit at
Shanghai Jiaotong University, China, in May 2008. I am very grateful
to Professor Cuipo Jiang for her hospitality. I would like to thank
the referees for valuable suggestions to put this paper in better
shape.

\section{$\hbar$-adic nonlocal vertex algebras and $\hbar$-adic
weak quantum vertex algebras}

In this section we study the notions of $\hbar$-adic  nonlocal
vertex algebra and $\hbar$-adic (weak) quantum vertex algebra, and
we present basic properties of $\hbar$-adic nonlocal vertex
algebras. The notion of $\hbar$-adic (weak) quantum vertex algebra
slightly generalizes Etingof-Kazhdan's notion of quantum vertex
operator algebra. In this paper, we use the standard formal variable
notations and conventions as established in \cite{flm} (cf.
\cite{ll}). The scalar field will be the field $\C$ of complex
numbers, $\N$ and $\Z_{+}$ denote the set of nonnegative integers
and the set of positive integers, respectively.

We start by recalling the notion of nonlocal vertex algebra (cf.
\cite{bk}, \cite{li-g1}). A {\em nonlocal vertex algebra} is a
vector space $V$, equipped with a linear map
\begin{eqnarray*}
Y:V&\rightarrow& \Hom (V,V((x)))
\subset (\End V)[[x,x^{-1}]]\nonumber\\
v&\mapsto& Y(v,x)=\sum_{n\in \Z}v_{n}x^{-n-1}\;\;\; (v_{n}\in \End V)
\end{eqnarray*}
and equipped with a distinguished vector ${\bf 1}\in V$, satisfying
the conditions that
\begin{eqnarray}
& & Y({\bf 1},x)=1,\\
& &Y(v,x){\bf 1}\in V[[x]]\;\;\;\mbox{ and }\;\;\;
\lim_{x\rightarrow 0}Y(v,x){\bf 1}\;(=v_{-1}{\bf 1})=v\ \ \
\mbox{for }v\in V,\label{edefcreation}
\end{eqnarray}
and that for $u,v,w\in V$, there exists $l\in \N$ such that
\begin{eqnarray}\label{eweakassoc}
(x_{0}+x_{2})^{l}Y(u,x_{0}+x_{2})Y(v,x_{2})w
=(x_{0}+x_{2})^{l}Y(Y(u,x_{0})v,x_{2})w
\end{eqnarray}
in $V[[x_{0}^{\pm 1},x_{2}^{\pm 1}]]$ (the {\em weak associativity}).

For a nonlocal vertex algebra $V$, define a linear operator
${\cal{D}}$ on $V$ by
\begin{eqnarray}
{\cal{D}}v=v_{-2}{\bf 1} \left(=\lim_{x\rightarrow
0}\frac{d}{dx}Y(v,x){\bf 1}\right) \;\;\;\mbox{ for }v\in V.
\end{eqnarray}
We have (\cite{li-g1}, Proposition 2.6)
\begin{eqnarray}\label{edproperty}
[{\cal{D}},Y(v,x)]=Y({\cal{D}}v,x)={d\over dx}Y(v,x) \;\;\;\mbox{
for }v\in V,
\end{eqnarray}
and furthermore,
\begin{eqnarray}
& &e^{x{\cal{D}}}Y(v,x_{1})e^{-x{\cal{D}}}=Y(e^{xD}v,x_{1})=Y(v,x_{1}+x),
\label{econjugationformula1}\\
& &Y(v,x){\bf 1}=e^{x{\cal{D}}}v.\label{ecreationwithd}
\end{eqnarray}

In \cite{li-qva1}, the following class of nonlocal vertex algebras
was singled out:

\bd{dweak-qva} {\em A {\em weak quantum vertex algebra} is a
nonlocal vertex algebra $V$, satisfying $\S$-locality: For $u,v\in
V$, there exist
$$\sum_{i=1}^{r}u^{(i)}\otimes v^{(i)}\otimes f_{i}(x)
\in V\otimes V\otimes \C((x)) $$ and a nonnegative integer $k$ such
that
\begin{eqnarray}\label{eSlocality-wqva-1}
(x_{1}-x_{2})^{k}Y(u,x_{1})Y(v,x_{2})
=(x_{1}-x_{2})^{k}\sum_{i=1}^{r}f_{i}(x_{2}-x_{1})
Y(v^{(i)},x_{2})Y(u^{(i)},x_{1}),
\end{eqnarray}
where $f_{i}(x_{2}-x_{1})$ is to be expanded in the nonnegative
powers of $x_{1}$, i.e., in view of the formal Taylor theorem,
$$f_{i}(x_{2}-x_{1})=e^{-x_{1}\frac{\partial}{\partial
    x_{2}}}f_{i}(x_{2})\in \C((x_{2}))[[x_{1}]].$$}
\ed

\br{rwqva-g1va} {\em  Let $V$ be a weak quantum vertex algebra. Let
$u,v,w\in V$ and assume that (\ref{eSlocality-wqva-1}) holds. Then
weak associativity relation (\ref{eweakassoc}) and
(\ref{eSlocality-wqva-1}) imply
\begin{eqnarray}\label{ejacobi-wva-def}
& &x_{0}^{-1}\delta\left(\frac{x_{1}-x_{2}}{x_{0}}\right)
Y(u,x_{1})Y(v,x_{2})w\nonumber\\
& &\ \ \ \ -x_{0}^{-1}\delta\left(\frac{x_{2}-x_{1}}{-x_{0}}\right)
\sum_{i=1}^{r} f_{i}(-x_{0})Y(v^{(i)},x_{2})Y(u^{(i)},x_{1})w\nonumber\\
&=&x_{2}^{-1}\delta\left(\frac{x_{1}-x_{0}}{x_{2}}\right)
Y(Y(u,x_{0})v,x_{2})w
\end{eqnarray}
(the {\em $\S$-Jacobi identity}). In fact, the notion of weak
quantum vertex algebra can be alternatively defined by using all the
axioms that define the notion of nonlocal vertex algebra except that
the weak associativity axiom is replaced by $\S$-Jacobi identity. }
\er

\bd{dyang-baxter} {\em Let $U$ be a vector space. A {\em unitary
rational quantum Yang-Baxter operator} (with one parameter) on $U$
is a linear map
$$\S(x): U\otimes U\rightarrow U\otimes U\otimes \C((x))$$
satisfying the condition that
\begin{eqnarray*}
& &\S^{21}(x)\S(-x)=1,\\
& &\S^{12}(x_{1})\S^{13}(x_{1}+x_{2}) \S^{23}(x_{2})
=\S^{23}(x_{2})\S^{13}(x_{1}+x_{2}) \S^{12}(x_{1}),
\end{eqnarray*}
where $\S^{21}(x)=P \S(x)P$ with $P$ the permutation operator on
$U\otimes U$.} \ed

The following notion is essentially due to Etingof and Kazhdan
\cite{ek}:

\bd{dqvoa-R} {\em A {\em quantum vertex algebra} is a nonlocal
vertex algebra $V$ equipped with a unitary rational quantum
Yang-Baxter operator $\S(x)$ on $V$, satisfying the conditions that
\begin{eqnarray}\label{eshift-constant}
[\D \otimes 1, \S(x)]=-\frac{d\S(x)}{dx}\ \ \ (\mbox{the {\em shift
condition}}),
\end{eqnarray}
and that for any $u,v\in V$, there exists a nonnegative integer $k$
such that
\begin{eqnarray}\label{es-locality-def}
& &(x_{1}-x_{2})^{k}Y(x_{1})(1\otimes Y(x_{2}))
(\S(x_{1}-x_{2})(u\otimes v)\otimes w)\nonumber\\
&=&(x_{1}-x_{2})^{k}Y(x_{2})(1\otimes Y(x_{1}))(v\otimes u\otimes w)
\end{eqnarray}
for all $w\in V$, and that
\begin{eqnarray}
\S(x)(Y(z)\otimes 1)=(Y(z)\otimes 1)\S^{23}(x)\S^{13}(x+z)
\end{eqnarray}
(the {\em hexagon identity}), where $Y(x): V\otimes V\rightarrow
V((x))$ is the map associated to the vertex operator map
$Y(\cdot,x): V\rightarrow \Hom (V,V((x)))$.} \ed

The following notion  is due to Etingof and Kazhdan \cite{ek}:

\bd{dnondeg} {\em Let $V$ be a nonlocal vertex algebra. For each
positive integer $n$, define a linear map
\begin{eqnarray}
Z_{n}: \C((x_{1}))\cdots ((x_{n}))\otimes V^{\otimes n} \rightarrow
V((x_{1}))\cdots ((x_{n}))
\end{eqnarray}
by
\begin{eqnarray}
Z_{n}(f\otimes v^{(1)}\otimes \cdots \otimes v^{(n)})=
f(x_{1},\dots,x_{n}) Y(v^{(1)},x_{1})\cdots Y(v^{(n)},x_{n}){\bf 1}
\end{eqnarray}
for $v^{(1)},\dots,v^{(n)}\in V,\; f\in \C((x_{1}))\cdots
((x_{n}))$. If all the linear maps $Z_{n}$ for $n\ge 1$ are
injective, $V$ is said to be {\em nondegenerate}. } \ed

The following proposition (\cite{li-qva1}, Theorems 4.8 and 5.11)
was lifted {}from \cite{ek}:

\bp{pfirst} Let $V$ be a weak quantum vertex algebra. Assume that
$V$ is nondegenerate. Then there exists a unique linear map $\S(x):
V\otimes V\rightarrow V\otimes V\otimes \C((x))$ satisfying the
condition that for any $u,v\in V$, there exists a nonnegative
integer $k$ such that (\ref{es-locality-def}) holds with
$$\S(x)(v\otimes u)=\sum_{i=1}^{r}v^{(i)}\otimes u^{(i)}\otimes
f_{i}(x).$$ Furthermore, $V$ equipped with $\S(x)$ is a quantum
vertex algebra. \ep

\bd{dquasi-module} {\em Let $V$ be a nonlocal vertex algebra. A {\em
$V$-module} is a vector space $W$ equipped with a linear map
\begin{eqnarray*}
Y_{W}: V &\rightarrow & \Hom (W,W((x)))
\subset (\End W)[[x,x^{-1}]]\nonumber\\
v&\mapsto& Y_{W}(v,x)=\sum_{n\in \Z}v_{n}x^{-n-1}\;\;\; (v_{n}\in
\End W)
\end{eqnarray*}
satisfying the conditions that
\begin{eqnarray}
Y_{W}({\bf 1},x)=1_{W}\;\; \mbox{ (where $1_{W}$ denotes the
identity operator on $W$)},
\end{eqnarray}
and that for any $u, v\in V,\; w\in W$, there exists $l\in \N$ such
that
\begin{eqnarray}\label{emoduleweakassoc}
(x_{0}+x_{2})^{l}Y_{W}(u,x_{0}+x_{2})Y_{W}(v,x_{2})w
=(x_{0}+x_{2})^{l}Y_{W}(Y(u,x_{0})v,x_{2})w.
\end{eqnarray}
A {\em quasi $V$-module} is defined by using all the above axioms
except that the last weak associativity axiom is replaced by a
weaker axiom: For any $u, v\in V,\; w\in W$, there exists $0\ne
p(x_{1},x_{2})\in \C[x_{1},x_{2}]$ such that
\begin{eqnarray}\label{equasimoduleweakassoc}
p(x_{0}+x_{2},x_{2})Y_{W}(u,x_{0}+x_{2})Y_{W}(v,x_{2})w
=p(x_{0}+x_{2},x_{2})Y_{W}(Y(u,x_{0})v,x_{2})w.
\end{eqnarray} }
\ed

Next, we study $\hbar$-adic analogues. Let $\hbar$ be a formal
variable throughout this paper. A $\C[[\hbar]]$-module $W$ is said
to be {\em torsion-free} if $\hbar w\ne 0$ for every $0\ne w\in W$,
and is said to be {\em separated} if $\cap_{n\ge 1}\hbar^{n}W=0$.
For a $\C[[\hbar]]$-module $W$, using subsets $w+\hbar^{n}W$ for
$w\in W,\; n\ge 1$ as the basis of open sets one obtains a topology
on $W$, which is called the {\em $\hbar$-adic topology}. A
$\C[[\hbar]]$-module $W$ is said to be {\em $\hbar$-adically
complete} if every Cauchy sequence in $W$ with respect to this
$\hbar$-adic topology has a limit in $W$.
 A $\C[[\hbar]]$-module $W$ is
{\em topologically free} if $W=W^{0}[[\h]]$ for some vector space
$W^{0}$ over $\C$. A fact is that a $\C[[\hbar]]$-module is
topologically free  if and only if it is torsion-free, separated,
and $\hbar$-adically complete (cf. \cite{kassel}).

\bd{dEhW} {\em Let $W$ be a $\C[[\h]]$-module. Define $\E_{\h}(W)$
to consist of each formal series
$$a(x)=\sum_{m\in \Z}a_{m}x^{-m-1}\in (\End W)[[x,x^{-1}]]$$
such that for every $w\in W$, $a_{m}w\rightarrow 0$, that is, for
every positive integer $n$,
$$a_{m}w\in \hbar^{n}W\;\;\;\mbox{ for $m$ sufficiently large}.$$}
\ed

For every $\C[[\h]]$-endomorphism $F$ of $W$, $F$ preserves the
submodules $\h^{n}W$ for $n\in \N$, so that $F$ gives rise to an
endomorphism of $W/\h^{n}W$ for each $n\in \N$. In this way, we have
natural $\C[[\h]]$-module homomorphisms
$$\tilde{\pi}_{n}: \End W\rightarrow \End (W/\h^{n}W)$$
for $n\in \N$. We also use $\tilde{\pi}_{n}$ for its canonical
extensions---the $\C[[\h]]$-module homomorphisms from $(\End
W)[[x_{1}^{\pm 1},\dots,x_{r}^{\pm 1}]]$ to $(\End
(W/\h^{n}W))[[x_{1}^{\pm 1},\dots, x_{r}^{\pm 1}]]$ for $r\ge 1$. In
terms of the maps $\tilde{\pi}_{n}$ we have
\begin{eqnarray}
\E_{\h}(W)=\left\{ a(x)\in (\End W)[[x,x^{-1}]]\;|\;
\tilde{\pi}_{n}(a(x))\in \E(W/\h^{n}W)\;\;\;\mbox{ for all }n\in
\N\right\}.
\end{eqnarray}

\bd{dqnoncomva} {\em An {\em $\h$-adic nonlocal vertex algebra} is a
 topologically free $\C[[\hbar]]$-module $V$, equipped with a
$\C[[\hbar]]$-module map
\begin{eqnarray*}
Y:V&\rightarrow& \E_{\hbar}(V)\subset (\End V)[[x,x^{-1}]]
\nonumber\\
v&\mapsto& Y(v,x)=\sum_{n\in \Z}v_{n}x^{-n-1}
\end{eqnarray*}
and equipped with a distinguished vector ${\bf 1}\in V$, satisfying
the conditions that
\begin{eqnarray*}
& &Y({\bf 1},x)=1,\\
& & Y(v,x){\bf 1}\in V[[x]]\;\;\;\mbox{ and }\;\;\;
\lim_{x\rightarrow 0}Y(v,x){\bf 1}\;(=v_{-1}{\bf 1})=v\ \ \ \mbox{
for }v\in V,\label{edefcreation-q}
\end{eqnarray*}
and that for $u,v,w\in V$ and for $n\in \N$, there exists $l\in \N$
such that
\begin{eqnarray}\label{eweakassoc-q}
(x_{0}+x_{2})^{l}Y(u,x_{0}+x_{2})Y(v,x_{2})w
\equiv (x_{0}+x_{2})^{l}Y(Y(u,x_{0})v,x_{2})w
\end{eqnarray}
modulo $\h^{n}V[[x_{0}^{\pm 1},x_{2}^{\pm 1}]]$ (the {\em $\h$-adic
weak associativity}).} \ed

\br{rinfinite-sum} {\em Notice that for $r,s\in \Z$, the coefficient
of the monomial $x_{0}^{r}x_{2}^{s}$ in the expression
$Y(u,x_{0}+x_{2})Y(v,x_{2})w$ is
$$\sum_{i\ge 0}\binom{r+i}{i}u_{-r-1-i}v_{-s-1+i}w,   $$
which is an infinite sum in general, even though it converges to an
element of $V$. This is one of the few places where $V$ needs to be
$\hbar$-adically complete. } \er

The following is a characterization of an $\h$-adic nonlocal vertex
algebra in terms of nonlocal vertex algebras over $\C$:

\bp{pinverse-limit} Let $V$ be a topologically free
$\C[[\h]]$-module, equipped with a vector ${\bf 1}\in V$ and a
$\C[[\h]]$-module map $Y$ {}from $V$ to $(\End V)[[x,x^{-1}]]$. Then
$(V,Y, {\bf 1})$ carries the structure of an $\h$-adic nonlocal
vertex algebra if and only if for every $n\in \N$, $(V/\h^{n}V,
Y^{(n)}, {\bf 1}+\hbar^{n}V)$ is a nonlocal vertex algebra over
$\C$, where $Y^{(n)}: V/\hbar^{n}V\rightarrow (\End
(V/\hbar^{n}V))[[x,x^{-1}]]$ is the canonical map reduced from $Y$.
\ep

\begin{proof}  {}From definition it is clear that
if $(V,Y,{\bf 1})$ is an $\h$-adic nonlocal vertex algebra,
$V/\h^{n}V$ is a nonlocal vertex algebra over $\C$ for every $n\in
\N$. Now, assume that for every $n\in \N$, $V/\h^{n}V$ is a nonlocal
vertex algebra over $\C$. For $v\in V$,  we have
$\tilde{\pi}_{n}(Y(v,x))\in \E(V/\h^{n}V)$ for $n\in \N$. Thus
$Y(v,x)\in \E_{\h}(V)$. On the other hand, for each $n\in \N$, with
${\bf 1}+\h^{n}V$ being the vacuum vector of $V/\h^{n}V$, we have
$Y({\bf 1},x)v-v\in \h^{n}V[[x,x^{-1}]]$ for $v\in V$. Because $V$
is separated, we must have $Y({\bf 1},x)v-v=0$. Similarly, we can
show $Y(v,x){\bf 1}\in V[[x]]$ and $\lim_{x\rightarrow 0}Y(v,x){\bf
1}=v$. The weak associativity of $V/\h^{n}V$ for $n\in \N$ exactly
amounts to the $\h$-adic weak associativity. Then $(V,Y,{\bf 1})$ is
an $\h$-adic nonlocal vertex algebra.
\end{proof}

\br{rhvoa-limit} {\em Let $V$ be an $\h$-adic nonlocal vertex
algebra. We have a projective inverse system of nonlocal vertex
algebras over $\C$ (or over $\C[[\h]]$):
\begin{eqnarray}
0\leftarrow V/\h V \leftarrow V/\h^{2}V \leftarrow V/\h^{3}V
\leftarrow \cdots
\end{eqnarray}
and the $\h$-adic nonlocal vertex algebra $V$ can be considered as
an inverse limit. } \er

Using Proposition \ref{pinverse-limit} (and the arguments in the
proof) we immediately have:

\bl{lDproperty}
Let $V$ be an $\h$-adic nonlocal vertex algebra. Define
$\D\in \End V$ by
\begin{eqnarray}
\D (v)=v_{-2}{\bf 1}\;\;\;\mbox{ for }v\in V.
\end{eqnarray}
Then
\begin{eqnarray*}
[\D, Y(v,x)]=Y(\D v,x)=\frac{d}{dx}Y(v,x)\;\;\;\mbox{ for }v\in V.
\end{eqnarray*}
\el

Let $V$ be an $\hbar$-adic nonlocal vertex algebra. Following
\cite{ek}, let
$$Y(x): V\hat{\otimes} V\rightarrow V[[x,x^{-1}]]$$
be the $\C[[\hbar]]$-module map associated to the vertex operator
map $Y: V\rightarrow (\End V)[[x,x^{-1}]]$, where here and forth
$V\hat{\otimes} V$ and $V\hat{\otimes} V\hat{\otimes}
\C((x))[[\hbar]]$ stand for the $\hbar$-adically completed tensor
products. If $V=V^{0}[[\h]]$ with $V^{0}$ a $\C$-vector space, we
have
$$V\hat{\otimes} V=(V^{0}\otimes V^{0})[[\hbar]]\ \
\mbox{ and }\ \ V\hat{\otimes} V\hat{\otimes}
\C((x))[[\hbar]]=(V^{0}\otimes V^{0}\otimes \C((x)))[[\hbar]].$$

\bd{dYx12} {\em Let $V$ be an $\hbar$-adic nonlocal vertex algebra.
Define a $\C[[\hbar]]$-module map
$$Y(x_{1},x_{2}): V\hat{\otimes} V\rightarrow (\End V)[[x_{1}^{\pm
1},x_{2}^{\pm 1}]]$$ by
\begin{eqnarray}
Y(x_{1},x_{2})(u\otimes v)(w)=Y(x_{1})(1\otimes Y(x_{2}))(u\otimes
v\otimes w)=Y(u,x_{1})Y(v,x_{2})w
\end{eqnarray}
for $u,v,w\in V$.} \ed

 \bd{dwqvoa} {\em An {\em $\h$-adic weak
quantum vertex algebra} is an $\h$-adic nonlocal vertex algebra $V$
which satisfies {\em $\hbar$-adic $\S$-locality}: For $u,v\in V$,
there exists
$$F(u,v,x)\in V\hat{\otimes} V\hat{\otimes} \C((x))[[\hbar]]$$
satisfying the condition that for every positive integer $n$, there
exists $k\in \N$ such that
\begin{eqnarray}
(x_{1}-x_{2})^{k}Y(u,x_{1})Y(v,x_{2})w \equiv
(x_{1}-x_{2})^{k}Y(x_{2})(1\otimes
Y(x_{1}))(F(u,v,x_{2}-x_{1})\otimes w) \ \ \
\end{eqnarray}
modulo $\h^{n}V[[x_{1}^{\pm 1},x_{2}^{\pm 1}]]$  for all $w\in V$.}
 \ed

\br{rwqva-jacobi} {\em Let $V$ be an $\hbar$-adic weak quantum
vertex algebra. We see that for every positive integer $n$,
$V/\h^{n}V$ is a weak quantum vertex algebra over $\C$. For
$u,v,w\in V$, by Remark \ref{rwqva-g1va} we have
\begin{eqnarray}
& &x_{0}^{-1}\delta\left(\frac{x_{1}-x_{2}}{x_{0}}\right)
Y(u,x_{1})Y(v,x_{2})w\nonumber\\
& &\;\;\;\; -x_{0}^{-1}\delta\left(\frac{x_{2}-x_{1}}{-x_{0}}\right)
Y(x_{2})(1\otimes Y(x_{1}))(F(u,v,-x_{0})\otimes w)\nonumber\\
&\equiv& x_{2}^{-1}\delta\left(\frac{x_{1}-x_{0}}{x_{2}}\right)
Y(Y(u,x_{0})v,x_{2})w
\end{eqnarray}
modulo $\h^{n}V[[x_{0}^{\pm 1},x_{1}^{\pm 1},x_{2}^{\pm 1}]]$. Since
$V$ is $\hbar$-adically complete, the coefficient of each monomial
$x_{0}^{r}x_{1}^{p}x_{2}^{q}$ for $r,p,q\in \Z$ in each of the three
main terms is an element of $V$. As $\cap _{n\ge 1}\h^{n}V=0$, we
obtain
\begin{eqnarray}\label{eqjacobi}
& &x_{0}^{-1}\delta\left(\frac{x_{1}-x_{2}}{x_{0}}\right)
Y(u,x_{1})Y(v,x_{2})w\nonumber\\
& &\;\;\;\; -x_{0}^{-1}\delta\left(\frac{x_{2}-x_{1}}{-x_{0}}\right)
Y(x_{2})(1\otimes Y(x_{1}))(F(u,v,-x_{0})\otimes w)\nonumber\\
&=&x_{2}^{-1}\delta\left(\frac{x_{1}-x_{0}}{x_{2}}\right)
Y(Y(u,x_{0})v,x_{2})w
\end{eqnarray}
(the {\em $\S$-Jacobi identity}). Clearly, $\S$-Jacobi identity is
equivalent to $\hbar$-adic weak associativity and $\hbar$-adic
$\S$-locality. In view of this, the notion of $\hbar$-adic weak
quantum vertex algebra can be defined alternatively  by using the
$\S$-Jacobi identity. } \er

\bd{dsim-relation} {\em Let $U$ be a $\C[[\h]]$-module and let $r$
be a positive integer. For
$$A(x_{1},\dots,x_{r}),\ \ B(x_{1},\dots,x_{r})
\in U[[x_{1}^{\pm 1},\dots,x_{r}^{\pm 1}]],$$ we write $A\sim_{\pm}
B$ if for every positive integer $n$, there exists a (nonzero)
polynomial
$$p(x_{1},\dots,x_{r})\in \< x_{i}\pm x_{j}\;|\; 1\le i<j\le r\>
\subset \C[x_{1},\dots,x_{r}]$$ such that
\begin{eqnarray}
p(x_{1},\dots,x_{r})(A(x_{1},\dots,x_{r})- B(x_{1},\dots,x_{r}))\in
\h^{n}U[[x_{1}^{\pm 1},\dots,x_{r}^{\pm 1}]].
\end{eqnarray}} \ed

Clearly, relations $\sim_{\pm}$ on $U[[x_{1}^{\pm
1},\dots,x_{r}^{\pm 1}]]$ are equivalence relations. It is also
clear that the left multiplication by a Laurent polynomial and the
formal partial differential operators $\partial/\partial
x_{1},\dots,
\partial/\partial x_{r}$ preserve the equivalence relations.
For convenience, we shall also use the notation $\sim$ for
$\sim_{-}$. If $V$ is an $\hbar$-adic nonlocal vertex algebra, for
$u,v,w\in V$ we have
$$Y(u,x_{0}+x_{2})Y(v,x_{2})w\sim_{+} Y(Y(u,x_{0})v,x_{2})w$$
in $V[[x_{0}^{\pm 1},x_{2}^{\pm 1}]]$. Furthermore, if $V$ is an
$\hbar$-adic weak quantum vertex algebra, for any $u,v\in V$, there
exists
$$F(u,v,x)\in V\hat{\otimes} V\hat{\otimes} \C((x))[[\hbar]]$$
such that
\begin{eqnarray*}
Y(u,x_{1})Y(v,x_{2}) \sim Y(x_{2}, x_{1})F(u,v,x_{2}-x_{1})
\end{eqnarray*}
in $(\End V)[[x_{1}^{\pm 1},x_{2}^{\pm 1}]]$.

\br{rfact-sim} {\rm Note that the equivalence relations $\sim_{\pm}$
when restricted to certain subspaces of $U[[x_{1}^{\pm
1},\dots,x_{r}^{\pm 1}]]$ amount to the equality relation. For
example, let $U=U^{0}[[\hbar]]$ be a topologically free
$\C[[\hbar]]$-module. For
$$A(x_{1},\dots,x_{r}),\ \ B(x_{1},\dots,x_{r})
\in (U^{0}((x_{1}))\cdots ((x_{r})))[[\hbar]],$$ if $A\sim_{\pm} B$,
we must have $A=B$. This is simply because $(U^{0}((x_{1}))\cdots
((x_{r})))[[\hbar]]$ is a vector space over the field
$\C((x_{1}))\cdots ((x_{r}))$ which contains $x_{i}\pm x_{j}$ for
$1\le i<j\le r$.} \er

\bp{p-skew-symmetry} Let $V$ be an $\hbar$-adic nonlocal vertex
algebra and let
$$u,v\in V,\;A(x)\in
V\hat{\otimes} V\hat{\otimes} \C((x))[[\hbar]].$$ Then
$$Y(u,x_{1})Y(v,x_{2})\sim Y(x_{2},x_{1})(A(x_{2}-x_{1}))$$
if and only if
\begin{eqnarray}\label{eskew-symmetry-1}
Y(u,x)v=e^{x\D}Y(-x)(A(-x)).
\end{eqnarray}
Furthermore, $V$ is an $\hbar$-adic weak quantum vertex algebra if
and only if there exists a $\C[[\hbar]]$-module map
$$\S(x): V\hat{\otimes} V\rightarrow V\hat{\otimes} V\hat{\otimes} \C((x))[[\hbar]]$$
such that
\begin{eqnarray}
Y(u,x)v=e^{x\D}Y(-x)\S(-x)(v\otimes u)\ \ \ \mbox{ for }u,v\in V.
\end{eqnarray}
 \ep

\begin{proof} We only need to prove the first part.
Note that for every positive integer $n$, $V/\hbar^{n}V$ is a
nonlocal vertex algebra over $\C$. By Corollary 5.3 of
\cite{li-qva1}, there exists a nonnegative integer $k$ such that
$$(x_{1}-x_{2})^{k}Y(u,x_{1})Y(v,x_{2})w\equiv (x_{1}-x_{2})^{k}Y(x_{2},x_{1})(A(x_{2}-x_{1}))w
\ \ \ \mod\;\hbar^{n}V$$ for all $w\in V$ if and only if
$$Y(u,x)v\equiv e^{x\D}Y(-x)(A(-x))\ \ \ \mod\;\hbar^{n}V.$$
Since $V$ is $\hbar$-adically complete, the coefficient of each
power of $x$ in $e^{x\D}Y(-x)(A(-x))$ is an element of $V$. As
$\cap_{n\ge 1}\hbar^{n}V=0$, the relations $Y(u,x)v\equiv
e^{x\D}Y(-x)(A(-x))\ \ \ \mod\;\hbar^{n}V$ for all $n\ge 1$ amount
to (\ref{eskew-symmetry-1}).
\end{proof}

The following is a reformulation and a slight generalization of
Etingof and Kazhdan's notion of quantum vertex operator algebra (see
\cite{ek}):

\bd{dqvoa} {\em  An {\em $\hbar$-adic quantum vertex algebra} is an
$\h$-adic nonlocal vertex algebra $V$
 equipped with a $\C[[\hbar]]$-module map
 \begin{eqnarray}
\S(x): V\hat{\otimes} V\rightarrow V\hat{\otimes} V\hat{\otimes}
\C((x))[[\hbar]],
\end{eqnarray}
which satisfies the {\em shift condition}
\begin{eqnarray}\label{eshift}
[\D \otimes 1, \S(x)]=-\frac{d\S(x)}{dx},
\end{eqnarray}
the {\em quantum Yang-Baxter equation}:
\begin{eqnarray}\label{eqYBe}
\S^{12}(x_{1})\S^{13}(x_{1}+x_{2})\S^{23}(x_{2})
=\S^{23}(x_{2})\S^{13}(x_{1}+x_{2})\S^{12}(x_{1})
\end{eqnarray}
and the {\em unitarity condition}:
\begin{eqnarray}\label{eunitarity}
\S^{21}(x)\S(-x)=1,
\end{eqnarray}
subject to the following axioms:

(QA1) The {\em $\hbar$-adic $\S$-locality}: For any $u,v\in V$ and
for any positive integer $n$, there exists $k\ge 0$ such that for
any $w\in V$ the series
$$(x_{1}-x_{2})^{k}Y(x_{1})(1\otimes Y(x_{2}))
(\S(x_{1}-x_{2})(u\otimes v)\otimes w)$$ and
$$(x_{1}-x_{2})^{k}Y(x_{2})(1\otimes Y(x_{1}))(v\otimes u\otimes w)$$
coincide modulo $\h^{n}V[[x_{1}^{\pm 1},x_{2}^{\pm 1}]]$.

(QA4) The {\em hexagon identity}:
\begin{eqnarray}\label{ehexagon}
\S(x_{1})(Y(x_{2})\otimes 1)=(Y(x_{2})\otimes
1)\S^{23}(x_{1})\S^{13}(x_{1}+x_{2}).
\end{eqnarray}}
\ed

Let $V$ be an $\hbar$-adic nonlocal vertex algebra. For a positive
integer $n$, we define a $\C[[\hbar]]$-module map
$$Z_{n}^{V}:V^{\hat{\otimes} n}\hat{\otimes} \C((x_{1})) \cdots
((x_{n}))[[\hbar]]\rightarrow V[[x_{1}^{\pm 1},\dots,x_{n}^{\pm
1}]]$$ as in Definition \ref{dnondeg}. Recall that $V/\hbar V$ is a
nonlocal vertex algebra over $\C$.

\bl{lstepone} Let $V=V^{0}[[\h]]$ be an $\h$-adic nonlocal vertex
algebra such that the nonlocal vertex algebra $V/\hbar V$ over $\C$
is nondegenerate. Then for every positive integer $n$, the
$\C[[\hbar]]$-module map $Z_{n}^{V}$ is injective. \el

\begin{proof} Let $n$ be any positive integer and
let $0\ne A\in V^{\hat{\otimes} n}\hat{\otimes} \C((x_{1})) \cdots
((x_{n}))[[\hbar]]$. Then
$$A=\h^{k}(A_{0}+\h A_{1}+\h^{2}A_{2}+\cdots )$$
for some $k\in \N$, $A_{i}\in (V^{0})^{\otimes n}\otimes \C((x_{1}))
\cdots ((x_{n}))$ with $A_{0}\ne 0$. Writing
$$Y(x)=Y_{0}(x)+\h Y_{1}(x)+\h^{2}Y_{2}(x)+\cdots$$
with $Y_{i}\in (\End V^{0})[[x,x^{-1}]]$ for $i\ge 0$, we have
\begin{eqnarray*}
& &Y(x_{1})(1\otimes Y(x_{2}))\cdots
(1^{\otimes (n-1)}\otimes Y(x_{n}))A\nonumber\\
&=&\h^{k}Y_{0}(x_{1})(1\otimes Y_{0}(x_{2}))\cdots
(1^{\otimes (n-1)}\otimes Y_{0}(x_{n}))(A_{0})
+O(\h^{k+1}).
\end{eqnarray*}
As $V/\hbar V$ is nondegenerate we have
\begin{eqnarray*}
Y_{0}(x_{1})(1\otimes Y_{0}(x_{2}))\cdots
(1^{\otimes (n-1)}\otimes Y_{0}(x_{n}))(A_{0})\ne 0,
\end{eqnarray*}
so that $Z_{n}^{V}(A)\ne 0$. This proves that $Z_{n}^{V}$ is
injective.
\end{proof}

The following is a reformulation of Proposition 1.11 of \cite{ek}:

\bp{pmain1} Let $V$ be an $\h$-adic weak quantum vertex algebra such
that the nonlocal vertex algebra $V/\hbar V$ over $\C$ is
nondegenerate. Then $\S$-locality defines a unique
$\C[[\hbar]]$-module map
\begin{eqnarray}
\S (x): V\hat{\otimes} V\rightarrow V\hat{\otimes} V\hat{\otimes}
\C((x)))[[\hbar]]
\end{eqnarray}
with $\S(x)(u\otimes v)=F(v,u,x)$ for $u,v\in V$ as in Definition
\ref{dwqvoa} and $(V,Y,{\bf 1},\S)$ carries the structure of an
$\h$-adic quantum vertex algebra. \ep

Next, we study modules and quasi modules for $\hbar$-adic nonlocal
vertex algebras.

\bd{dhvoa-module} {\em  Let $V$ be an $\h$-adic nonlocal vertex
algebra. A {\em $V$-module} is a topologically free
$\C[[\hbar]]$-module $W$, equipped with a $\C[[\hbar]]$-module map
$$Y_{W}: V\rightarrow \E_{\hbar}(W)\subset (\End W)[[x,x^{-1}]]$$
satisfying the conditions that $Y_{W}({\bf 1},x)=1_{W}$ and that for
$u,v\in V,\; w\in W$ and for every positive integer $n$, there
exists $l\in \N$ such that
\begin{eqnarray}\label{eweakassoc-qmodule}
(x_{0}+x_{2})^{l}Y_{W}(u,x_{0}+x_{2})Y_{W}(v,x_{2})w \equiv
(x_{0}+x_{2})^{l}Y_{W}(Y(u,x_{0})v,x_{2})w
\end{eqnarray}
modulo $\h^{n}W[[x_{0}^{\pm 1},x_{2}^{\pm 1}]]$. We define a notion
of {\em quasi $V$-module} by replacing the $\h$-adic weak
associativity with the following axiom: For $u,v\in V,\; w\in W$ and
for every positive integer $n$, there exists $0\ne p(x_{1},x_{2})\in
\C[x_{1},x_{2}]$ such that
\begin{eqnarray}\label{eweakassoc-quasimodule}
p(x_{0}+x_{2},x_{2})Y_{W}(u,x_{0}+x_{2})Y_{W}(v,x_{2})w \equiv
p(x_{0}+x_{2},x_{2})Y_{W}(Y(u,x_{0})v,x_{2})w
\end{eqnarray}
modulo $\h^{n}W[[x_{0}^{\pm 1},x_{2}^{\pm 1}]]$.} \ed

We have the following straightforward analogue of Proposition
\ref{pinverse-limit}:

\bp{phadic-module} Let $V$ be an $\h$-adic nonlocal vertex algebra,
let $W$ be a topologically free $\C[[\h]]$-module, and let $Y_{W}$
be a $\C[[\h]]$-module map {}from $V$ to $(\End W)[[x,x^{-1}]]$.
Then $(W,Y_{W})$ is a (quasi) $V$-module if and only if for every
positive integer $n$, $W/\h^{n}W$ is a (quasi) $V/\h^{n}V$-module.
\ep

The following is an $\hbar$-adic version of a result of
\cite{li-qva1}:

\bp{phmodule-algebra-relations} Let $V$ be an $\h$-adic nonlocal
vertex algebra, let
$$m\in \Z,\; u,v,c^{(0)},c^{(1)},\dots \in V,\;
A(x)\in V\hat{\otimes} V\hat{\otimes} \C((x))[[\hbar]]$$ such that
$\lim_{j\rightarrow \infty}c^{(j)}=0$, and let $(W,Y_{W})$ be a
$V$-module. If
\begin{eqnarray}\label{e-hcross-bracket}
& &(x_{1}-x_{2})^{m}Y(u,x_{1})Y(v,x_{2}) -(-x_{2}+x_{1})^{m}
Y(x_{2},x_{1})(A(x))
\nonumber\\
& &\ \ =\sum_{j\ge 0}Y(c^{(j)},x_{2})\frac{1}{j!}
\left(\frac{\partial}{\partial x_{2}}\right)^{j}
x_{2}^{-1}\delta\left(\frac{x_{1}}{x_{2}}\right),
\end{eqnarray}
then
\begin{eqnarray}\label{e-hcross-bracket-module}
& &(x_{1}-x_{2})^{m}Y_{W}(u,x_{1})Y_{W}(v,x_{2}) -(-x_{2}+x_{1})^{m}
Y_{W}(x_{2},x_{1})(A(x))
\nonumber\\
& &\ \ =\sum_{j\ge 0}Y_{W}(c^{(j)},x_{2})\frac{1}{j!}
\left(\frac{\partial}{\partial x_{2}}\right)^{j}
x_{2}^{-1}\delta\left(\frac{x_{1}}{x_{2}}\right).
\end{eqnarray}
If $(W,Y_{W})$ is a faithful $V$-module, the converse also holds.
\ep

\begin{proof} For any positive integer $n$, $V/\hbar^{n}V$ is a
nonlocal vertex algebra over $\C$ and  $W/\hbar^{n}W$ is a
$(V/\hbar^{n}V)$-module. It follows from \cite{li-qva1} (Proposition
6.7) that after applied to a vector in $W$,
(\ref{e-hcross-bracket-module}) holds modulo $\hbar^{n}W$. With $W$
topologically free, $W$ is $\hbar$-adically complete and separated.
Then (\ref{e-hcross-bracket-module}) must hold. For the converse,
for any positive integer $n$, denote by $\rho_{n}$ the
$\C[[\hbar]]$-module map from $V$ to $\E(W/\hbar^{n}W)$. Clearly,
$\hbar^{n}V\subset \ker \rho_{n}$. Then $V/\ker\rho_{n}$ is a
nonlocal vertex algebra over $\C$ and $W/\hbar^{n}W$ is a faithful
$(V/\ker \rho_{n})$-module. Again, from \cite{li-qva1} (Proposition
6.7), (\ref{e-hcross-bracket}) modulo $\ker\rho_{n}$ holds. For
$v\in \cap_{n\ge 1}\ker \rho_{n}$, with $W$ separated we have
$Y_{W}(v,x)=0$. As $W$ is faithful, we must have $\cap_{n\ge 1}\ker
\rho_{n}=0$. Since $V$ is $\hbar$-adically complete,
(\ref{e-hcross-bracket}) holds on $V$.
\end{proof}

We shall also need the following variation:

\bp{phmodule-algebra-relations-v2} Let $V$ be an $\h$-adic nonlocal
vertex algebra, let
$$m\in \Z,\; u,v,c^{(0)},c^{(1)},\dots \in V,\;
A(x)\in V\hat{\otimes} V\hat{\otimes} \C((x))[[\hbar]]$$ such that
$\lim_{j\rightarrow \infty}c^{(j)}=0$, and let $(W,Y_{W})$ be a
$V$-module.
 If
\begin{eqnarray}\label{e-hcross-bracket-v2}
& &(x_{1}-x_{2})^{m}Y(u,x_{1})Y(v,x_{2}) -(-x_{2}+x_{1})^{m}
Y(x_{2},x_{1})(A(x))
\nonumber\\
& &\ \ =\sum_{j\ge 0}Y(c^{(j)},x_{1})\frac{1}{j!}
\left(\frac{\partial}{\partial x_{2}}\right)^{j}
x_{2}^{-1}\delta\left(\frac{x_{1}}{x_{2}}\right),
\end{eqnarray}
(note that the only change is on the variable $x$ for vertex
operators $Y(c^{(j)},x)$), then
\begin{eqnarray}\label{e-hcross-bracket-module-v2}
& &(x_{1}-x_{2})^{m}Y_{W}(u,x_{1})Y_{W}(v,x_{2}) -(-x_{2}+x_{1})^{m}
Y_{W}(x_{2},x_{1})(A(x))
\nonumber\\
& &\ \ =\sum_{j\ge 0}Y_{W}(c^{(j)},x_{1})\frac{1}{j!}
\left(\frac{\partial}{\partial x_{2}}\right)^{j}
x_{2}^{-1}\delta\left(\frac{x_{1}}{x_{2}}\right).
\end{eqnarray}
If $(W,Y_{W})$ is a faithful $V$-module, the converse also holds.
\ep

\begin{proof} Let $(U,Y_{U})$ be a $V$-module, e.g., $U=V$ or $U=W$.
For any $v\in V,\; w\in U$ and for any $n\ge 1$, with $U/\hbar^{n}U$
a module for $V/\hbar^{n}V$ viewed as a nonlocal vertex algebra over
$\C$, we have
\begin{eqnarray*}
Y_{U}(\D v,x)w\equiv \frac{d}{dx}Y_{U}(v,x)w\ \ \ \mod\; \hbar^{n}U.
\end{eqnarray*}
Since $U$ is separated, we have
\begin{eqnarray*}
Y_{U}(\D v,x)w=\frac{d}{dx}Y_{U}(v,x)w.
\end{eqnarray*}
Using this we get
\begin{eqnarray*}
& &\sum_{j\ge 0}Y_{U}(c^{(j)},x_{1})\frac{1}{j!}
\left(\frac{\partial}{\partial x_{2}}\right)^{j}
x_{2}^{-1}\delta\left(\frac{x_{1}}{x_{2}}\right)\\
&=&\sum_{j\ge 0}\frac{1}{j!} \left(\frac{\partial}{\partial
x_{2}}\right)^{j}\left(Y_{U}(c^{(j)},x_{2})
x_{2}^{-1}\delta\left(\frac{x_{1}}{x_{2}}\right)\right)\\
&=&\sum_{j\ge 0}\sum_{i=0}^{j}\frac{1}{(j-i)!}
\left(\left(\frac{\partial}{\partial
x_{2}}\right)^{j-i}Y_{U}(c^{(j)},x_{2})\right)\frac{1}{i!}
\left(\frac{\partial}{\partial x_{2}}\right)^{i}
x_{2}^{-1}\delta\left(\frac{x_{1}}{x_{2}}\right)\\
&=&\sum_{j\ge 0}\sum_{i=0}^{j}\frac{1}{(j-i)!}
Y_{U}(\D^{j-i}c^{(j)},x_{2})\frac{1}{i!}
\left(\frac{\partial}{\partial x_{2}}\right)^{i}
x_{2}^{-1}\delta\left(\frac{x_{1}}{x_{2}}\right)\\
&=&\sum_{i\ge 0}\sum_{r\ge 0}\frac{1}{r!}
Y_{U}(\D^{r}c^{(r+i)},x_{2})\frac{1}{i!}
\left(\frac{\partial}{\partial x_{2}}\right)^{i}
x_{2}^{-1}\delta\left(\frac{x_{1}}{x_{2}}\right).
\end{eqnarray*}
Notice that for $i\ge 0$, $\sum_{r\ge
0}\frac{1}{r!}\D^{r}c^{(r+i)}\in V$ (as $V$ is $\hbar$-adically
complete). Then it follows from Proposition
\ref{phmodule-algebra-relations}.
\end{proof}

\section{Some technical results}
In this section we present certain technical results which we need
in Section 4. In particular, we study $\hbar$-adic nonlocal vertex
subalgebras and subalgebras generated by subsets of an $\hbar$-adic
nonlocal vertex algebra.

\bd{dsubalgebra} {\em Let $V$ be an $\h$-adic nonlocal vertex
algebra. An {\em $\hbar$-adic nonlocal vertex subalgebra} is a
$\C[[\hbar]]$-submodule containing ${\bf 1}$ such that
 $(U,Y, {\bf 1})$ carries the structure of
 an $\hbar$-adic nonlocal vertex algebra. In particular,
 $U$ is a topologically free submodule.}
 \ed

We say that a $\C[[\hbar]]$-submodule $U$ of an $\hbar$-adic
nonlocal vertex algebra $V$ is {\em $Y$-closed} if  $u_{m}v\in U$
for all $u,v\in U,\; m\in \Z$. (This is to distinguish the algebraic
closedness from the topological closedness.)

\br{rinduced-top} {\em Let $U$ be a $\C[[\hbar]]$-submodule of a
topologically free $\C[[\hbar]]$-module $V$. With $\hbar^{n}U\subset
U\cap \hbar^{n}V$ for $n\in \N$, we see that the induced topology on
$U$ from $V$ (with the $\hbar$-adic topology) coincides with the
$\hbar$-adic topology of $U$ if and only if for any $n\in \N$, there
exists $k\in \N$ such that $U\cap \hbar^{k}V\subset \hbar^{n}U$.  }
\er

\bp{psubalgebra1} Let $V$ be an $\h$-adic nonlocal vertex algebra
and let $U$ be a $\C[[\hbar]]$-submodule satisfying the conditions
that ${\bf 1}\in U$, $U$ is $Y$-closed, and that the induced
topology on $U$ from $V$ coincides with its own $\hbar$-adic
topology. In addition we assume that $U$ is $\hbar$-adically
complete. Then $U$ is an $\hbar$-adic nonlocal vertex subalgebra of
$V$. \ep

\begin{proof} Notice that as a submodule of $V$, $U$ is torsion-free and separated.
Since $U$ is also $\hbar$-adically complete, $U$ is topologically
free. Let $u,v\in U$ and $n\in \N$. {}From assumption, there exists
$k\in \N$ such that $U\cap \hbar^{k}V\subset \hbar^{n}U$. With
$u,v\in U\subset V$ and $k\in \N$ being fixed, there exists $r\in
\N$ such that $u_{m}v\in \hbar^{k}V$ for $m\ge r$.
 Then
 $$u_{m}v\in U\cap \hbar^{k}V\subset \hbar^{n}U\ \ \ \mbox{ for
 }m\ge r.$$
That is,
$$Y(u+\hbar^{n}U,x)(v+\hbar^{n}U)\in (U/\hbar^{n}U)((x)).$$
Furthermore, let $w\in V$. By $\hbar$-adic weak associativity, there
exists $l\in \N$ such that
$$(x_{0}+x_{2})^{l}Y(u,x_{0}+x_{2})Y(v,x_{2})w\equiv
(x_{0}+x_{2})^{l}Y(Y(u,x_{0})v,x_{2})w\ \ (\mod \;\hbar^{k}V).$$
Because $U$ is $Y$-closed and $\hbar$-adically complete, and because
$U\cap \hbar^{k}V\subset \hbar^{n}U$, we have
$$(x_{0}+x_{2})^{l}Y(u,x_{0}+x_{2})Y(v,x_{2})w\equiv
(x_{0}+x_{2})^{l}Y(Y(u,x_{0})v,x_{2})w\ \ (\mod \;\hbar^{n}U).$$
Now, $(U,Y,{\bf 1})$ satisfies all the axioms for an $\hbar$-adic
nonlocal vertex algebra.
\end{proof}

\bd{dbracketK} {\em Let $M$ be a $\C[[\h]]$-module. For any
$\C[[\h]]$-submodule $K$, we define
\begin{eqnarray}
[K]=\{ w\in M\;|\; \h^{n}w\in K\;\;\mbox{ for some }n\in \N\}.
\end{eqnarray}}
\ed

The following two lemmas are straightforward:

\bl{lsomeproperty} Let $M$ be a $\C[[\h]]$-module and let $K$ be a
$\C[[\h]]$-submodule such that $[K]=K$. Then $K\cap \h^{n}M=\h^{n}K$
for all $n\in \N$. In particular, the induced topology on $K$ from
$M$ (with the $\hbar$-adic topology) coincides with the $\h$-adic
topology of $K$. \el

\bl{lbracket-closed} Let $V$ be an $\hbar$-adic nonlocal vertex
algebra and let $K$ be a $Y$-closed $\C[[\hbar]]$-submodule. Then
both $[K]$ and the $\hbar$-adic completion $\widehat{K}$ of $K$ in
$V$ are $Y$-closed. \el

Furthermore, we have:

\bp{pbracket-closed-sub} Let $V$ be an $\hbar$-adic nonlocal vertex
algebra and let $K$ be a $Y$-closed $\C[[\hbar]]$-submodule
containing ${\bf 1}$. Then
$\left[\widehat{[K]}\right]=\widehat{[K]}$ and $\widehat{[K]}$ is an
$\hbar$-adic nonlocal vertex subalgebra of $V$. \ep

\begin{proof} Since $[[K]]=[K]$, by Lemma \ref{lsomeproperty} we have $[K]\cap
\hbar^{n}V=\hbar^{n}[K]$ for $n\in \N$. Then $\widehat{[K]}$ is
$\hbar$-adically complete with respect to its only $\hbar$-adic
topology. Then $\widehat{[K]}$ is topologically free. {}From Lemma
\ref{lbracket-closed}, $\widehat{[K]}$ is $Y$-closed. If we can
prove $\left[\widehat{[K]}\right]=\widehat{[K]}$, then by
Proposition \ref{psubalgebra1}, $\widehat{[K]}$ is an $\hbar$-adic
nonlocal vertex subalgebra of $V$. Let $u\in
\left[\widehat{[K]}\right]$. By definition, there exists $k\in \N$
such that $\hbar^{k}u\in \widehat{[K]}$. Furthermore, there exists a
Cauchy sequence $\{a_{n}\}$ in $[K]$ with  $\hbar^{k}u$ as the
limit. Then there exists $r\ge 1$ such that $a_{n}-\hbar^{k}u\in
\hbar^{k}V$ for $n\ge r$. As $V$ is torsion-free, for each $n\ge r$,
there exists uniquely $b_{n}\in V$ such that $a_{n}=\hbar^{k}b_{n}$.
As $[[K]]=[K]$ and $\hbar^{n}b_{n}=a_{n}\in [K]$, we have $b_{n}\in
[K]$ for $n\ge r$. Using the fact that $V$ is torsion-free, we see
that $\{b_{n}\}_{n\ge r}$ is a Cauchy sequence in $[K]$, converging
to $u$. Thus $u\in \widehat{[K]}$. This proves
$\left[\widehat{[K]}\right]=\widehat{[K]}$, concluding the proof.
\end{proof}

Now, let $U$ be a subset of an $\hbar$-adic nonlocal vertex algebra
$V$. Let $U^{(1)}$ be the $\C[[\hbar]]$-span of $U\cup \{{\bf 1}\}$
and then inductively define $U^{(n+1)}$ for $n\ge 1$ to be the
$\C[[\hbar]]$-span of the vectors $a_{m}b$ for $a,b\in U^{(n)},\;
m\in \Z$. {}From definition we have ${\bf 1}\in U^{(n)}\subset
U^{(n+1)}$ for $n\ge 1$.  Set
\begin{eqnarray}
\<U\>^{o}=\left[\cup_{n\ge 1}U^{(n)}\right]\subset V
\end{eqnarray}
and furthermore, set
\begin{eqnarray}
\<U\>=\widehat{\<U\>^{o}} \ \ \ \mbox{(the $\hbar$-adic
completion)}.
\end{eqnarray}
 We have:

\bp{psubalgebra} Let $V$ be an $\h$-adic nonlocal vertex algebra and
let $U$ be a subset of $V$. Then $\<U\>$ is an $\h$-adic nonlocal
vertex subalgebra satisfying the condition that $U\subset \<U\>$ and
$[\<U\>]=\<U\>$. Furthermore, any $\hbar$-adic nonlocal vertex
subalgebra $H$, which satisfies the condition that $U\subset H$ and
$[H]=H$, contains $\<U\>$.  \ep

\begin{proof} It follows from the definition that
 $\cup_{n\ge 1}U^{(n)}$ is $Y$-closed and contains $\{{\bf 1}\}\cup U$.
Then the first assertion follows from Proposition
\ref{pbracket-closed-sub}.
 Let $H$ be an $\hbar$-adic nonlocal vertex subalgebra such that $[H]=H$ and $U\subset H$.
 Then $\cup_{n\ge 1}U^{(n)}\subset H$. Furthermore, we have
$$\<U\>^{o}=\left[\cup_{n\ge 1}U^{(n)}\right]\subset [H]=H.$$
Since $[\<U\>^{o}]=\<U\>^{o}$, the induced topology on $\<U\>^{o}$
(from $H$ or $V$) coincides with its own $\hbar$-adic topology. Then
$\<U\>\subset H$ as $H$ is $\hbar$-adically complete.
\end{proof}

We shall need the following result later:

\bl{lvacuum-like} Let $V$ be an $\h$-adic nonlocal vertex algebra
with a generating subset $U$ in the sense that $V=\<U\>$ and let
$(W,Y_{W})$ be a quasi $V$-module equipped with a
$\C[[\hbar]]$-linear operator $\D$ such that
\begin{eqnarray}
[\D,Y_{W}(u,x)]=\frac{d}{dx}Y_{W}(u,x)\;\;\;\mbox{ for }u\in U.
\end{eqnarray}
Assume that $w$ is a vector of $W$ such that $\D w=0$. Then
$Y(v,x)w\in V[[x]]$ for all $v\in V$ and the linear map $\phi$
defined by $\phi(v)= v_{-1}w$ for $v\in V$ is a $V$-module
homomorphism. \el

\begin{proof} For every positive integer $n$, $V/\h^{n}V$ is a nonlocal vertex
algebra over $\C$ and $W/\h^{n}W$ is a quasi $(V/\h^{n}V)$-module.
Set
$$K_{n}=(\cup_{k\ge 1}U^{(k)}+\hbar^{n}V)/\hbar^{n}V\subset V/\hbar^{n}V.$$
It is clear that $K_{n}$ is a nonlocal vertex subalgebra of
$V/\hbar^{n}V$.
 Let
$\phi_{n}: K_{n}\rightarrow W/\h^{n}W$ be the map induced {}from
$\phi$. {}From \cite{li-qva1} (Proposition 6.2), we have
$$v_{r}w\in \hbar^{n}W \ \ \ \mbox{ for }
v\in \cup_{k\ge 1}U^{(k)},\ \ r\ge 0$$ and $\phi_{n}$ is a
$K_{n}$-module homomorphism, which amounts to
\begin{eqnarray*}
\phi (u_{m}v)\equiv  u_{m}\phi(v)\;\;\;\mod\; \h^{n}W \ \ \mbox{ for
}u,v\in \cup_{k\ge 1}U^{(k)},\ m\in \Z.
\end{eqnarray*}
As $\cap_{n\ge 1}\h^{n}W=0$, we have $v_{r}w=0$ for $v\in \cup_{k\ge
1}U^{(k)},\; r\ge 0$ and
$$\phi(u_{m}v)= u_{m}\phi(v)\;\;\;\mbox{ for }u,v\in \cup_{k\ge 1}U^{(k)},\ m\in \Z.$$

Now, let $u,v\in [\cup_{k\ge 1}U^{(k)}]$. There exists $t\in \N$
such that $\hbar^{t}u,\hbar^{t}v\in \cup_{k\ge 1}U^{(k)}$. Then
$$\hbar^{t}v_{r}w=0\ \ \mbox{ for }r\ge 0\ \mbox{ and }\
\hbar^{2t}\phi(u_{m}v)=\hbar^{2t} u_{m}\phi(v) \;\;\;\mbox{ for
}m\in \Z.$$ Since $W$ is torsion-free, we have $v_{r}w=0$ for $r\ge
0$ and $\phi(u_{m}v)=u_{m}\phi(v)$ for $m\in \Z$. Recall that
$\<U\>$ is the completion of $[\cup_{k\ge 1}U^{(k)}]$. Let
$u^{(i)},v^{(i)}$ $(i\ge 1)$ be sequences in $[\cup_{k\ge
1}U^{(k)}]$, converging to $u, v\in V$, respectively. We have
$$v^{(i)}_{r}w=0\ \ \mbox{ for }r\ge 0\ \mbox{ and }\
\phi(u^{(i)}_{m}v^{(j)})=u^{(i)}_{m}\phi(v^{(j)}) \;\;\;\mbox{ for
}i,j\ge 1,\; m\in \Z.$$ {}From this we see that for every positive
integer $n$,
$$v_{r}w\in \hbar^{n}W \ \ \mbox{ for }r\ge 0\ \ \mbox{ and } \
\phi(u_{m}v)-u_{m}\phi(v)\in \hbar^{n}W \ \ \ \mbox{ for }m\in \Z.$$
Again as $\cap_{n\ge 1}\hbar^{n}W=0$, we get
$$v_{r}w=0 \ \ \mbox{ for }r\ge 0\ \ \mbox{ and } \ \phi(u_{m}v)= u_{m}\phi(v)
\ \ \ \mbox{ for }m\in \Z.$$  Now we have proved
$$Y(v,x)w\in V[[x]],\ \ \ \ \phi(u_{m}v)= u_{m}\phi(v)
\ \ \ \mbox{ for all }u, v\in \<U\>,\; m\in \Z,$$ as needed.
\end{proof}

Let $W$ be a topologically free $\C[[\hbar]]$-module and let $A$ be
a $\C$-subspace of $W$. Notice that for any sequence
$\{a_{n}\}_{n\ge 0}$ in $A$, $\sum_{n\ge 0}a_{n}\hbar^{n}\in W$. Set
\begin{eqnarray}
A[[\hbar]]'=\left\{ \sum_{n\ge 0}a_{n}\hbar^{n}\;|\; a_{n}\in
A\right\},
\end{eqnarray}
which is a $\C[[\hbar]]$-submodule of $W$.

\bd{dlocalsubset-va} {\em Let $A$ and $B$ be subsets of an
$\hbar$-adic nonlocal vertex algebra $V$. We say that the ordered
pair $(A,B)$ is {\em $\hbar$-adically $\S$-local} if for any $a\in
A,\; b\in B$, there exists
$$P(a,b; x)\in \left((\C B)\otimes (\C A)\otimes \C((x))\right)[[\hbar]]'
\subset  V\hat{\otimes} V \hat{\otimes} \C((x)))[[\hbar]] $$ such
that
\begin{eqnarray}\label{esim-proof}
Y(a,x_{1})Y(b,x_{2})\sim Y(x_{2},x_{1})P(a,b;x_{2}-x_{1}).
\end{eqnarray}
We say that a subset $A$ of $V$ is {\em $\hbar$-adically $\S$-local}
if $(A,A)$ is $\hbar$-adically $\S$-local. } \ed

We have the following technical result:

\bl{linduction} Let $A$ and $B$ be $\C$-subspaces of $V$ such that
$(A,B)$ is $\hbar$-adically $\S$-local. Then $(A,B^{(2)})$ and
$(A^{(2)},B)$ are $\hbar$-adically $\S$-local. \el

\begin{proof} From definition, there exists a $\C[[\hbar]]$-module map
$$\S(x): B[[\hbar]]'\hat{\otimes} A[[\hbar]]'
\rightarrow B[[\hbar]]'\hat{\otimes} A[[\hbar]]'\hat{\otimes}
\C((x))[[\hbar]]$$  such that for $a\in A,\; b\in B$,
\begin{eqnarray}\label{esimilar-relation}
Y(a,x_{1})Y(b,x_{2})\sim Y(x_{2},x_{1})\S(x_{2}-x_{1})(b\otimes a).
\end{eqnarray}
We have the maps
\begin{eqnarray*}
&&\S^{32}(x): B[[\hbar]]'\hat{\otimes} A[[\hbar]]'\hat{\otimes}
B[[\hbar]]'\rightarrow B[[\hbar]]'\hat{\otimes}
A[[\hbar]]'\hat{\otimes} B[[\hbar]]' \hat{\otimes}
\C((x))[[\hbar]],\\
 &&\S^{13}(x): B[[\hbar]]'\hat{\otimes}
B[[\hbar]]'\hat{\otimes} A[[\hbar]]'\rightarrow
B[[\hbar]]'\hat{\otimes} B[[\hbar]]'\hat{\otimes} A[[\hbar]]'
\hat{\otimes} \C((x))[[\hbar]].
\end{eqnarray*}
Note that by Proposition
\ref{p-skew-symmetry}, (\ref{esimilar-relation}) is equivalent to
$$Y(a,x)b=e^{x\D}Y(-x)\S(-x)(b\otimes a).$$
Let $a\in A,\; u,v\in B$. Using Lemma \ref{lDproperty} we get
\begin{eqnarray*}
Y(a,x)Y(u,z)v
&\sim & Y(z)(1\otimes Y(x))(\S(z-x)(u\otimes a)\otimes v)\\
&=&Y(z)(1\otimes e^{x\D}Y(-x))\S^{32}(-x)(\S(z-x)(u\otimes a)\otimes v)\\
&=&e^{x\D}Y(z-x)(1\otimes Y(-x))\S^{32}(-x)(\S(z-x)(u\otimes a)\otimes v)\\
&\sim&e^{x\D}Y(-x)(Y(z)\otimes 1)\S^{32}(-x)(\S(z-x)(u\otimes
a)\otimes v)\\
&\sim&e^{x\D}Y(-x)(Y(z)\otimes 1)\S^{32}(-x)(\S(-x+z)(u\otimes
a)\otimes v)
\end{eqnarray*}
in $V[[x^{\pm 1},z^{\pm 1}]]$. In view of Remark \ref{rfact-sim} we
have
\begin{eqnarray*}
Y(a,x)Y(u,z)v=e^{x\D}Y(-x)(Y(z)\otimes
1)\S^{32}(-x)(\S(-x+z)(u\otimes a)\otimes v).
\end{eqnarray*}
It follows that $(A,B^{(2)})$ is $\hbar$-adically $\S$-local.

Similarly,  for $a,b\in A,\; u\in U$, we have
\begin{eqnarray*}
Y(Y(a,x_{0})b,x)u &\sim_{+} & Y(a,x_{0}+x)Y(b,x)u\\
&=& Y(a,x_{0}+x)e^{x\D}Y(-x)\S(-x)(u\otimes b)\\
&=& e^{x\D}Y(a,x_{0})Y(-x)\S(-x)(u\otimes b)\\
&\sim_{+}& e^{x\D}Y(-x)(1\otimes
Y(x_{0}))\S^{13}(-x-x_{0})(\S(-x)(u\otimes b)\otimes a)
\end{eqnarray*}
in $V[[x_{0}^{\pm 1},x^{\pm 1}]]$, which by Remark \ref{rfact-sim}
implies
\begin{eqnarray*}
Y(Y(a,x_{0})b,x)u=e^{x\D}Y(-x)(1\otimes
Y(x_{0}))\S^{13}(-x-x_{0})(\S(-x)(u\otimes b)\otimes a).
\end{eqnarray*}
It follows that $(A^{(2)},B)$ is $\hbar$-adically $\S$-local.
\end{proof}

Now, we have (cf. \cite{li-qva2}, \cite{ltw}):

\bp{pgenerator-Slocal} Let $V$ be an $\hbar$-adic nonlocal vertex
algebra and let $U$ be an $\hbar$-adically $\S$-local subset such
that $V=(\cup_{n\ge 1}U^{(n)})[[\hbar]]'$. Then $V$ is an
$\hbar$-adic weak quantum vertex algebra. \ep

\begin{proof} We must prove that $V$ as a subset of $V$ is
$\hbar$-adically $\S$-local. Because $U$ is $\hbar$-adically
$\S$-local, it follows from Lemma \ref{linduction} (and induction)
that $\cup_{n\ge 1}U^{(n)}$ is $\hbar$-adically $\S$-local. Let
$u,v\in V$. {}From assumption, we have
$$u=\sum_{i\ge 0}a(i)\hbar^{i},\ \ v=\sum_{j\ge 0}b(j)\hbar^{j}
\ \ \mbox{ with }a(i),b(j)\in \cup_{n\ge 1}U^{(n)}.$$ For any
$i,j\in \N$, there exists
$$A_{i,j}(x)\in V\hat{\otimes} V\hat{\otimes} \C((x))[[\hbar]]$$
such that
$$Y(a(i),x_{1})Y(b(j),x_{2})\sim Y(x_{2},x_{1})A_{i,j}(x_{2}-x_{1}). $$
Notice that
$$\sum_{i,j\in \N}A_{i,j}(x)\hbar^{i+j}
\in V\hat{\otimes} V\hat{\otimes} \C((x))[[\hbar]]$$
and
$$Y(u,x_{1})Y(v,x_{2})
\sim Y(x_{2},x_{1})\left(\sum_{i,j\in \N}A_{i,j}(x_{2}-x_{1})\right). $$
This proves that $V$ is $\hbar$-adically $\S$-local.
\end{proof}

Using the proof of Proposition 2.8 in \cite{li-qva2} and Proposition
\ref{pgenerator-Slocal}, we immediately have:

\bp{pvacuum-generating} Let $V$ be an $\hbar$-adic nonlocal vertex
algebra, let $U$ be an $\hbar$-adically $\S$-local subset, and let
$W$ be a $V$-module with $e\in W$ such that $Y_{W}(u,x)e\in V[[x]]$
for $u\in U$. Set $K=(\cup_{n\ge 1}U^{(n)})[[\hbar]]'\subset V$.
Then $Y_{W}(v,x)e\in V[[x]]$ for all $v\in K$. Furthermore, the
$\C[[\hbar]]$-module map $\theta: K\rightarrow W$, defined by
$\theta (v)=v_{-1}e$ for $v\in K$, satisfies that
$$\phi(Y(u,x)v)=Y_{W}(u,x)\phi(v)\ \ \ \mbox{ for }u,v\in K.$$
\ep

\section{A general construction of $\hbar$-adic quantum vertex
algebras}

In this section we give a general construction of $\hbar$-adic
nonlocal vertex algebras and $\hbar$-adic weak quantum vertex
algebras by using what we call $\hbar$-adic quasi-compatible sets of
vertex operators on topologically free $\C[[\hbar]]$-modules.

We start by recalling from \cite{li-qva1} (cf. \cite{li-g1}) the
general construction of nonlocal vertex algebras from
quasi-compatible sets of formal vertex operators. Let $W^{0}$ be a
vector space over $\C$. Set
\begin{eqnarray}
\E(W^{0})=\Hom (W^{0},W^{0}((x))).
\end{eqnarray}
The identity operator on $W^{0}$, denoted by $1_{W^{0}}$, is a
typical element of $\E(W^{0})$, and the formal differential operator
$\frac{d}{dx}$ is an endomorphism of $\E(W^{0})$.

\bd{dcompatibility} {\em An (ordered) sequence
$(\psi^{(1)}(x),\dots,\psi^{(r)}(x))$ in $\E(W^{0})$ is said to be
{\em quasi-compatible} if there exists $0\ne p(x_{1},x_{2})\in
\C[x_{1},x_{2}]$ such that
\begin{eqnarray}
\left(\prod_{1\le i<j\le r}p(x_{i},x_{j})\right)
\psi^{(1)}(x_{1})\cdots \psi^{(r)}(x_{r}) \in \Hom
(W^{0},W^{0}((x_{1},\dots,x_{r}))).
\end{eqnarray}
A subset $U$ of $\E(W^{0})$ is said to be {\em quasi-compatible} if
every (ordered) finite sequence in $U$ is quasi-compatible. We also
define a notion of {\em compatibility} by assuming that
$p(x_{1},x_{2})$ is of the form $(x_{1}-x_{2})^{k}$ with $k\in \N$.}
\ed

Assume that $(a(x),b(x))$ is a quasi compatible pair in $\E(W^{0})$.
By definition, there exists $0\ne p(x_{1},x_{2})\in \C[x_{1},x_{2}]$
such that
\begin{eqnarray}\label{epab-condition}
p(x_{1},x_{2})a(x_{1})b(x_{2})\in \Hom (W^{0},W^{0}((x_{1},x_{2}))).
\end{eqnarray}
Recall from \cite{li-qva1} that
$$\iota_{x_{1},x_{2}}:\ \
 \C_{*}(x_{1},x_{2})\rightarrow \C((x_{1}))((x_{2}))$$
is the algebra-embedding that preserves each element of
$\C[[x_{1},x_{2}]]$, where $\C_{*}(x_{1},x_{2})$ denotes the algebra
extension of $\C[[x_{1},x_{2}]]$ by inverting every nonzero
polynomial. We have
$$\iota_{x,x_{0}}(1/p(x_{0}+x,x))
\left(p(x_{1},x)a(x_{1})b(x)\right)|_{x_{1}=x+x_{0}} \in \left(\Hom
(W^{0},W^{0}((x)))\right)((x_{0})).
$$

\bd{danbdefinition} {\em Let $(a(x), b(x))$ be a quasi-compatible
pair in $\E(W^{0})$. Define $a(x)_{n}b(x)$ for $n\in \Z$, elements
of $\E(W^{0})$,  in terms of generating function
\begin{eqnarray*}
Y_{\cal{E}}(a(x),x_{0})b(x) =\sum_{n\in \Z}a(x)_{n}b(x) x_{0}^{-n-1}
\end{eqnarray*}
by
\begin{eqnarray*}
Y_{\cal{E}}(a(x),x_{0})b(x)=\iota_{x,x_{0}}(1/p(x_{0}+x,x))
\left(p(x_{1},x)a(x_{1})b(x)\right)|_{x_{1}=x+x_{0}},
\end{eqnarray*}
where $p(x_{1},x_{2})$ is any nonzero polynomial such that
(\ref{epab-condition}) holds. } \ed

A quasi-compatible subspace $U$ of $\E(W^{0})$ is said to be {\em
$Y_{\E}$-closed} if
\begin{eqnarray}
a(x)_{n}b(x)\in U\;\;\;\mbox{ for }a(x),b(x)\in U,\; n\in \Z.
\end{eqnarray}

We have (\cite{li-qva1}, Theorem 2.19; cf. \cite{li-g1}):

\bt{tgeneratingthem} Let $W^{0}$ be a vector space over $\C$ and let
$U$ be any quasi-compatible subset of $\E(W^{0})$. There exists a
(unique) smallest $Y_{\E}$-closed quasi-compatible subspace $\<U\>$
of $\E(W^{0})$, containing $U$ and $1_{W^{0}}$, and
$(\<U\>,Y_{\cal{E}},1_{W^{0}})$ carries the structure of a nonlocal
vertex algebra with $U$ as a generating subset and $W^{0}$ is a
quasi-module for $\<U\>$ with
$Y_{W^{0}}(\alpha(x),x_{0})=\alpha(x_{0})$ for $\alpha(x)\in \<U\>$.
Furthermore, if $U$ is compatible, $W^{0}$ is a module for $\<U\>$.
\et

\bd{dS-local} {\em Let $W^{0}$ be a vector space over $\C$. A subset
$U$ of $\E(W^{0})$ is said to be {\em $\S$-local} if for any
$a(x),b(x)\in U$, there exist $a^{(i)}(x),b^{(i)}(x)\in U,\;
f_{i}(x)\in \C((x))$ $(i=1,\dots,r)$ such that
\begin{eqnarray}\label{eRslocal-def}
(x_{1}-x_{2})^{k}a(x_{1})b(x_{2})
=(x_{1}-x_{2})^{k}\sum_{i=1}^{r}f_{i}(x_{2}-x_{1})
b^{(i)}(x_{2})a^{(i)}(x_{1})
\end{eqnarray}
for some $k\in \N$.} \ed

{}From \cite{li-qva1} (Lemma 3.2), every $\S$-local subset of
$\E(W^{0})$ is quasi compatible. In fact, the same proof shows that
every $\S$-local subset is compatible. Furthermore, we have (see
\cite{li-qva1}, Theorem 5.8):

\bt{tweak-qva} Let $W^{0}$ be a vector space over $\C$ and let $U$
be an $\S$-local subset of $\E(W^{0})$. Then the nonlocal vertex
algebra $\<U\>$ generated by $U$ is a weak quantum vertex algebra
with $W^{0}$ as a faithful module. \et

Now, let $W$ be a $\C[[\hbar]]$-module. Recall from Section 2 that
$\E_{\hbar}(W)$ is the $\C[[\hbar]]$-submodule of $(\End
W)[[x,x^{-1}]]$, consisting of each $a(x)=\sum_{m\in
\Z}a_{m}x^{-m-1}$ satisfying the condition that for any $w\in W,\
n\in \N$, there exists $k\in \Z$ such that
$$a_{m}w\in \hbar^{n}W\ \ \ \mbox{ for }m\ge k.$$

For the rest of this section we assume that {\em $W=W^{0}[[\hbar]]$
is a fixed topologically free $\C[[\h]]$-module.} Then $W$ is
torsion-free, separated in the sense that $\cap_{n\ge 1}\h^{n}W=0$,
and $\hbar$-adically complete.  We have the following projective
inverse system
\begin{eqnarray}\label{einverse-system}
0\leftarrow W/\h W \leftarrow W/\h^{2}W\leftarrow
W/\h^{3}W\leftarrow \cdots
\end{eqnarray}
(equipped with the canonical maps from $W/\h^{n+1}W$ to
$W/\hbar^{n}W$ for $n\ge 0$) with $W$ as an inverse limit. Let $F$
be an endomorphism of $W$. For every nonnegative integer $n$, $F$
gives rise to an endomorphism $F_{n}$ of $W/\h^{n}W$. Then we have
an endomorphism $\{F_{n}\}$ of the projective inverse system
(\ref{einverse-system}). Conversely, given any endomorphism, a
sequence $\{ f_{n}\}$,
 of the inverse system (\ref{einverse-system}),
we have an endomorphism $f$ of $W$.

The $\C[[\hbar]]$-module $\End W$ can be naturally identified with
$(\End W^{0})[[\hbar]]$ and we have
$$(\End W)[[x,x^{-1}]]=(\End W^{0})[[x,x^{-1}]][[\hbar]],$$
which is topologically free. Furthermore, we have
\begin{eqnarray}
\E_{\h}(W)=\E (W^{0})[[\h]],
\end{eqnarray}
which is also topologically free.

\bl{ltorsion-free} For $a(x)\in (\End W)[[x,x^{-1}]]$, if
$\h^{k}a(x)\in \E_{\h}(W)$ for some $k\in \N$, then $a(x)\in
\E_{\h}(W)$. \el

\begin{proof}  For any $n\in \N,\; w\in W$, with $\h^{k}a(x)\in \E_{\h}(W)$,
there exists $q\in \Z$ such that $\h^{k}a_{m}w\in \h^{k+n}W$ for
$m\ge q$, where $a(x)=\sum_{m\in \Z}a_{m}x^{-m-1}$. Since $W$ is
torsion-free, we have $a_{m}w\in \h^{n}W$ for $m\ge q$. This proves
$a(x)\in \E_{\h}(W)$.
\end{proof}

\br{rexp1} {\em  For each $n\in \N$, we have a canonical
$\C[[\hbar]]$-module map
\begin{eqnarray}\label{epin}
\tilde{\pi}_{n}: (\End W)[[x,x^{-1}]]\rightarrow (\End
(W/\hbar^{n}W))[[x,x^{-1}]].
\end{eqnarray}
 As $W$ is torsion-free, we have
$$\ker \tilde{\pi}_{n}=\hbar^{n}(\End W)[[x,x^{-1}]].$$
Recall from Section 2 that an element $a(x)$ of $(\End
W)[[x,x^{-1}]]$ lies in $\E_{\hbar}(W)$ if and only if
$\tilde{\pi}_{n}(a(x))\in \E(W/\hbar^{n}W)$ for $n\in \N$. Then we
have canonical $\C[[\hbar]]$-module maps
\begin{eqnarray}
\pi_{n}: \E_{\hbar}(W)\rightarrow \E(W/\hbar^{n}W)
\end{eqnarray}
 for $n\in \N$, where by Lemma \ref{ltorsion-free}
\begin{eqnarray}
\ker \pi_{n}=\E_{\hbar}(W)\cap \hbar^{n}(\End W)[[x,x^{-1}]]
=\hbar^{n}\E_{\hbar}(W).
\end{eqnarray}
For every $n\in \N$, we have a canonical $\C[[\hbar]]$-module map
$$\theta_{n}:\E(W/\h^{n+1}W)\rightarrow
\E(W/\h^{n}W).$$
We have the following projective inverse system
\begin{eqnarray}\label{einverse-system-E}
0\leftarrow \E(W/\h W) \leftarrow \E(W/\h^{2}W)\leftarrow
\E(W/\h^{3}W)\leftarrow \cdots
\end{eqnarray}
with $\E_{\hbar}(W)$ equipped with $\C[[\hbar]]$-module maps
$\pi_{n}$ as an inverse limit. Then for any sequence
$\{\psi_{n}(x)\}$ with $\psi_{n}(x)\in \E(W/\h^{n}W)$ for $n\in \N$,
satisfying the condition that
$\theta_{n}(\psi_{n+1}(x))=\psi_{n}(x)$, there exists a unique
$\psi(x)\in \E_{\h}(W)$ such that $\pi_{n}(\psi(x))=\psi_{n}(x)$ for
$n\in \N$.} \er

\bd{dcompatibility-topo} {\em A finite sequence $a^{1}(x),\dots,
a^{r}(x)$ in $\E_{\h}(W)$ is said to be {\em $\h$-adically
quasi-compatible} if for every positive integer $n$, the sequence
$\pi_{n}(a^{1}(x)),\dots, \pi_{n}(a^{r}(x))$ in $\E(W/\h^{n}W)$ is
quasi-compatible. A subset $U$ of $\E_{\h}(W)$ is said to be {\em
$\h$-adically quasi-compatible} if every finite sequence in $U$ is
$\h$-adically quasi-compatible. Correspondingly, we define notions
of {\em $\h$-adically compatible} sequence and {\em $\h$-adically
compatible} subset.} \ed

Let $r$ be a positive integer. For each $n\in \N$, we have a
canonical $\C[[\hbar]]$-module map
$$\tilde{\pi}_{n}^{(r)}:
(\End W)[[x_{1}^{\pm 1},\dots,x_{r}^{\pm 1}]]\rightarrow (\End
(W/\hbar^{n}W))[[x_{1}^{\pm 1},\dots,x_{r}^{\pm 1}]],$$ where
$\tilde{\pi}_{n}^{(1)}= \tilde{\pi}_{n}$ defined in (\ref{epin}). We
also have a canonical $\C[[\hbar]]$-module map:
$$\tilde{\theta}_{n}^{(r)}: (\End (W/\hbar^{n+1}W))[[x_{1}^{\pm 1},\dots,x_{r}^{\pm 1}]]
\rightarrow (\End (W/\hbar^{n}W))[[x_{1}^{\pm 1},\dots,x_{r}^{\pm
1}]].$$ It is clear that for $n\in \N$,
\begin{eqnarray}\label{ecomposition}
\tilde{\pi}_{n}^{(r)}=\tilde{\theta}_{n}^{(r)}\circ
\tilde{\pi}_{n+1}^{(r)}.
\end{eqnarray}

For any vector space $U$ over $\C$, we set
\begin{eqnarray}
\E^{(r)}(U) =\Hom(U, U((x_{1},\dots,x_{r}))),
\end{eqnarray}
which is naturally a $\C((x_{1},\dots,x_{r}))$-module.

\bd{dernw} {\em Let $r$ be a positive integer. For every $n\in \N$,
define $\E_{n}^{(r)}(W)$ to be the $\C[[\hbar]]$-submodule of $(\End
W)[[x_{1}^{\pm 1},\dots,x_{r}^{\pm 1}]]$, consisting of each formal
series
$$\psi(x_{1},\dots,x_{r})=\sum_{m_{1},\dots,m_{r}\in
  \Z}\psi(m_{1},\dots,m_{r})x_{1}^{-m_{1}-1}\cdots x_{r}^{-m_{r}-1},$$
 satisfying the condition that
$$\tilde{\pi}_{n}^{(r)}(\psi(x_{1},\dots,x_{r}))\in \E^{(r)}(W/\hbar^{n}W),$$
or equivalently, for every $w\in W$, there exists $k\in \Z$ such
that
$$\psi(m_{1},\dots,m_{r})w\in \hbar^{n}W\ \ \ \mbox{whenever
}m_{i}\ge k\ \mbox{ for some }1\le i\le r.$$}
\ed

We see that $\E^{(r)}_{n}(W)$ are also
$\C((x_{1},\dots,x_{r}))$-modules and we have
$$\E_{0}^{(r)}(W)\supset \E_{1}^{(r)}(W) \cdots \supset
\E_{n}^{(r)}(W)\supset \E_{n+1}^{(r)}(W)\supset \cdots. $$
 In terms of  $\E_{n}^{(r)}(W)$, a
sequence $\psi^{1}(x),\dots, \psi^{r}(x)$ in $\E_{\h}(W)$ is
$\hbar$-adically quasi-compatible if and only if for every $n\in
\N$, there exists $0\ne p(x,y)\in \C[x,y]$ such that
$$\left(\prod_{1\le i<j\le r}p(x_{i},x_{j})\right)
\psi^{1}(x_{1})\cdots \psi^{r}(x_{r}) \in \E_{n}^{(r)}(W).$$
Generalizing the maps $\pi_{n}$ and $\theta_{n}$, we have canonical
$\C[[\hbar]]$-module maps for $n\in \N$:
\begin{eqnarray*}
\pi_{n}^{(r)}:& &\E_{n}^{(r)}(W)\rightarrow \E^{(r)}(W/\hbar^{n}W),\\
\theta_{n}^{(r)}: & &\E^{(r)}(W/\hbar^{n+1}W) \rightarrow
\E^{(r)}(W/\hbar^{n}W),
\end{eqnarray*}
which satisfy
\begin{eqnarray*}
\pi_{n}^{(r)}=\theta_{n}^{(r)}\circ \pi_{n+1}^{(r)}.
\end{eqnarray*}
Set
\begin{eqnarray}
\E_{\hbar}^{(r)}(W)=\cap_{n\ge 1}\E_{n}^{(r)}(W)\subset (\End
W)[[x_{1}^{\pm 1},\dots,x_{r}^{\pm 1}]].
\end{eqnarray}

Note that if $(a(x), b(x))$ is an $\hbar$-adically quasi-compatible
 pair in $\E_{\h}(W)$, then for every $n\in \N$,
$(\pi_{n}(a(x)),\pi_{n}(b(x)))$ is a quasi-compatible pair in
$\E(W/\h^{n}W)$ and hence $\pi_{n}(a(x))_{m}\pi_{n}(b(x))$ are
defined for all $m\in \Z$.

\bl{lmorphism} Let $(a(x), b(x))$ be an $\hbar$-adically
quasi-compatible pair in $\E_{\h}(W)$. We have
\begin{eqnarray}
\theta_{n+1}\left(\pi_{n+1}(a(x))_{m}\pi_{n+1}(b(x))\right)
=\pi_{n}(a(x))_{m}\pi_{n}(b(x))
\end{eqnarray}
for $n\in \N,\; m\in \Z$.
\el

\begin{proof} For any fixed $n\in \N$,
let $p(x,y)\in \C[x,y]$ be a nonzero
  polynomial such that
$$p(x_{1},x_{2})a(x_{1})b(x_{2})
\in  \E^{(2)}_{n+1}(W)\subset \E^{(2)}_{n}(W).$$ {}From Definition
\ref{danbdefinition}, we have
\begin{eqnarray*}
p(x_{0}+x,x)Y_{\E}(\pi_{n}(a(x)),x_{0})\pi_{n}(b(x))
&=&\left(p(x_{1},x)\pi_{n}(a(x_{1}))\pi_{n}(b(x))\right)|_{x_{1}=x+x_{0}},\\
&=&\pi_{n}^{(2)}\left(p(x_{1},x)a(x_{1})b(x)\right)|_{x_{1}=x+x_{0}},\\
p(x_{0}+x,x)Y_{\E}(\pi_{n+1}(a(x)),x_{0})\pi_{n+1}(b(x))
&=&\left(p(x_{1},x)\pi_{n+1}(a(x_{1}))\pi_{n+1}(b(x))\right)|_{x_{1}=x+x_{0}}\\
&=&\pi_{n+1}^{(2)}\left(p(x_{1},x)a(x_{1})b(x))\right)|_{x_{1}=x+x_{0}}.
\end{eqnarray*}
With (\ref{ecomposition}) it follows that
\begin{eqnarray*}
p(x_{0}+x,x)Y_{\E}(\pi_{n}(a(x)),x_{0})\pi_{n}(b(x))=
p(x_{0}+x,x)\theta_{n+1}\left(Y_{\E}(\pi_{n+1}(a(x)),x_{0})\pi_{n+1}(b(x))\right),
\end{eqnarray*}
which implies
\begin{eqnarray*}
Y_{\E}(\pi_{n}(a(x)),x_{0})\pi_{n}(b(x))=
\theta_{n+1}\left(Y_{\E}(\pi_{n+1}(a(x)),x_{0})\pi_{n+1}(b(x))\right),
\end{eqnarray*}
as desired.
\end{proof}

Using Lemma \ref{lmorphism} we define the following partial
operations on $\E_{\h}(W)$:

 \bd{dinverse-limit-def} {\em Let $(a(x), b(x))$ be an
$\hbar$-adically quasi-compatible (order) pair in $\E_{\h}(W)$. For
$m\in \Z$, we define
\begin{eqnarray}
a(x)_{m}b(x)=\lim_{\leftarrow}\pi_{n}(a(x))_{m}\pi_{n}(b(x))
\in \E_{\hbar}(W).
\end{eqnarray}
Form the generating function
\begin{eqnarray}
Y_{\E}(a(x),x_{0})b(x)=\sum_{m\in \Z}(a(x)_{m}b(x)) x_{0}^{-m-1}.
\end{eqnarray}}
\ed

{}From definition, for every positive integer $n$ we have
\begin{eqnarray}
\pi_{n}(a(x)_{m}b(x))=\pi_{n}(a(x))_{m}\pi_{n}(b(x))
\end{eqnarray}
for $m\in \Z$. Namely,
\begin{eqnarray} \pi_{n}\left(Y_{\E}(a(x),x_{0})b(x)\right)
=Y_{\E}(\pi_{n}(a(x)),x_{0})\pi_{n}(b(x)).
\end{eqnarray}

Recall that $\E_{\hbar}^{(2)}(W)$ consists of each
$\psi(x_{1},x_{2})\in (\End W)[[x_{1}^{\pm 1},x_{2}^{\pm 1}]]$, such
that for every $n\in \N$, $\tilde{\pi}_{n}^{(2)}(\psi)\in
\E^{(2)}(W/\hbar^{n}W).$

 \bp{panother-version} Let $a(x),b(x)\in \E_{\hbar}(W)$.
Assume that there exists $p(x_{1},x_{2},\hbar)\in
\C[x_{1},x_{2},\hbar]$ with $p(x_{1},x_{2},0)\ne 0$ such that
$$p(x_{1},x_{2},\hbar)a(x_{1})b(x_{2})\in \E_{\hbar}^{(2)}(W).$$
Then $(a(x),b(x))$ is $\hbar$-adically quasi-compatible and
\begin{eqnarray*}
Y_{\E}(a(x),x_{0})b(x)=\iota_{x,x_{0},\hbar}(1/p(x_{0}+x,x,\hbar))
\left(p(x_{1},x,\hbar)a(x_{1})b(x)\right)|_{x_{1}=x+x_{0}}.
\end{eqnarray*}
\ep

\begin{proof} Set $$f(x_{1},x_{2})=p(x_{1},x_{2},0)\in
\C[x_{1},x_{2}],\ \ A=p(x_{1},x_{2},\hbar)a(x_{1})b(x_{2})\in
\E_{\hbar}^{(2)}(W).$$
 Then
$$p(x_{1},x_{2},\hbar)=f(x_{1},x_{2})-\hbar g(x_{1},x_{2},\hbar)$$
for some $g(x_{1},x_{2},\hbar)\in \C[x_{1},x_{2},\hbar]$. We have
$$a(x_{1})b(x_{2})= \iota_{x_{1},x_{2},\hbar}(1/p(x_{1},x_{2},\hbar))A
=\sum_{k\ge
0}\iota_{x_{1},x_{2}}(f(x_{1},x_{2})^{-k-1})g(x_{1},x_{2},\hbar)^{k}\hbar^{k}A.
$$
For any positive integer $n$, we have
$$f(x_{1},x_{2})^{n}a(x_{1})b(x_{2})\equiv \sum_{k=0}^{n-1}
f(x_{1},x_{2})^{n-k-1}g(x_{1},x_{2},\hbar)^{k}\hbar^{k}A \ \ \ \mod
\hbar^{n}(\End W)[[x_{1}^{\pm 1},x_{2}^{\pm 1}]],$$ so that
\begin{eqnarray}
\tilde{\pi}^{(2)}_{n}\left(f(x_{1},x_{2})^{n}a(x_{1})b(x_{2})\right)\in
\E^{(2)}(W/\hbar^{n}W),
\end{eqnarray}
as $\tilde{\pi}^{(2)}_{n}(A)\in \E^{(2)}(W/\hbar^{n}W)$. This proves
that $(a(x),b(x))$ is $\hbar$-adically quasi-compatible.
Furthermore, for $n\in \N$ we have
\begin{eqnarray*}
f(x_{0}+x,x)^{n}\pi_{n}(Y_{\E}(a(x),x_{0})b(x))
&=&\left(f(x_{1},x)^{n}\pi_{n}(a(x_{1}))\pi_{n}(b(x))\right)|_{x_{1}=x_{2}+x_{0}}\\
&=&\pi_{n}^{(2)}\left(f(x_{1},x)^{n}a(x_{1})b(x)\right)|_{x_{1}=x_{2}+x_{0}}.
\end{eqnarray*}
 Then
\begin{eqnarray*}
& &\pi_{n}^{(2)} \iota_{x,x_{0},\hbar}(p(x_{0}+x,x,\hbar)^{-1})
(p(x_{1},x,\hbar)a(x_{1})b(x))|_{x_{1}=x+x_{0}}\\
&=&\pi_{n}^{(2)}
\iota_{x,x_{0},\hbar}(p(x_{0}+x,x,\hbar)^{-1}f(x_{0}+x,x)^{-n})
(p(x_{1},x,\hbar)f(x_{1},x)^{n}a(x_{1})b(x))|_{x_{1}=x+x_{0}}\\
&=&\pi_{n}^{(2)}
\iota_{x,x_{0},\hbar}(p(x_{0}+x,x,\hbar)^{-1}f(x_{0}+x,x)^{-n})p(x_{0}+x,x,\hbar)
(f(x_{1},x)^{n}a(x_{1})b(x))|_{x_{1}=x+x_{0}}\\
&=&\iota_{x,x_{0},\hbar}(f(x_{0}+x,x)^{-n})\pi_{n}^{(2)}(f(x_{1},x)^{n}a(x_{1})b(x))|_{x_{1}=x+x_{0}}\\
 &=&\pi_{n}\left(Y_{\E}(a(x),x_{0})b(x)\right),
\end{eqnarray*}
{}from which the second part follows. \end{proof}

An $\hbar$-adically quasi-compatible $\C[[\hbar]]$-submodule $K$ of
$\E_{\h}(W)$ is said to be {\em $Y_{\E}$-closed} if
$$a(x)_{m}b(x)\in K\ \ \ \mbox{ for }a(x),b(x)\in K,\; m\in \Z.$$

\bp{phvoa-1} Let $V$ be a $Y_{\E}$-closed $\hbar$-adically
quasi-compatible $\C[[\h]]$-submodule of $\E_{\h}(W)$, containing
$1_{W}$. Suppose that $[V]=V$ and $V$ is $\hbar$-adically complete.
Then $(V,Y_{\E},1_{W})$ carries the structure of an $\h$-adic
nonlocal vertex algebra and $W$ is a faithful quasi $V$-module with
$Y_{W}(a(x),x_{0})=a(x_{0})$ for $a(x)\in V$. Furthermore, if $V$ is
$\hbar$-adically compatible, $W$ is a module (instead of a quasi
module).\ep

\begin{proof} Recall that $\E_{\h}(W)$ is topologically free. As a submodule of
$\E_{\h}(W)$, $V$ is torsion-free and separated. Being assumed to be
$\h$-adically complete, $V$ is topologically free. For any $n\in
\N$, as $\pi_{n}(a(x))_{m}\pi_{n}(b(x))=\pi_{n}(a(x)_{m}b(x))$ for
$a(x),b(x)\in V,\; m\in \Z$, we see that $\pi_{n}(V)$ is a
$Y_{\E}$-closed quasi-compatible $\C[[\h]]$-submodule of
$\E(W/\h^{n}W)$, containing $1_{W}$. It follows from Theorem
\ref{tgeneratingthem} that $\pi_{n}(V)$ is a nonlocal vertex algebra
over $\C$ with $W/\h^{n}W$ as a quasi module.

Now we prove that the map $\pi_{n}$ from $V$ to $\pi_{n}(V)$ reduces
to a $\C[[\h]]$-isomorphism from $V/\h^{n}V$ onto $\pi_{n}(V)$, so
that $V/\h^{n}V$ is a nonlocal vertex algebra over $\C$. Let
$a(x)\in V$ be such that $\pi_{n}(a(x))=0$ in $\E(W/\h^{n}W)$. Then
$a(x)W\subset \h^{n}W[[x,x^{-1}]]$. So $a(x)=\h^{n}b(x)$ for some
$b(x)\in (\End W)[[x,x^{-1}]]$. By Lemma \ref{ltorsion-free},
$b(x)\in \E_{\h}(W)$. Then we have $b(x)\in [V]=V$. Thus
$a(x)=\h^{n}b(x)\in \h^{n}V$. This proves that $V\cap \ker
\pi_{n}=\h^{n}V$, which implies $V/\h^{n}V\simeq \pi_{n}(V)\subset
\E(W/\h^{n}W)$. Consequently, $V/\h^{n}V$ is a nonlocal vertex
algebra over $\C$.  By Propositions \ref{pinverse-limit} and
\ref{phadic-module}, $V$ is an $\h$-adic nonlocal vertex algebra
with $W$ as a quasi $V$-module. The last part follows from Theorem
\ref{tgeneratingthem} and Proposition \ref{phadic-module}.
\end{proof}

For convenience, we call any $\hbar$-adic nonlocal vertex algebra
$V$ in Proposition \ref{phvoa-1} an {\em $\hbar$-adic nonlocal
vertex subalgebra of $\E_{\hbar}(W)$}.

 \bl{lpre} Let $U$ be an $\hbar$-adically quasi-compatible $\C[[\hbar]]$-submodule
of $\E_{\h}(W)$. Then $[U]$ is $\hbar$-adically quasi-compatible. If
$U$ is $Y_{\E}$-closed, so is $[U]$. \el

\begin{proof} Notice that $[U]\subset \E_{\h}(W)$ by Lemma \ref{ltorsion-free}.
As $W$ is torsion-free, for any $w\in W,\; s,n\in \N$,
the relation $\h^{s}w\in \h^{s+n}W$ implies $w\in \h^{n}W$.
Furthermore, for $\psi\in (\End W)[[x_{1}^{\pm 1},\dots,x_{r}^{\pm
1}]]$, $s,n\in \N$, the relation $\h^{s}\psi \in \E^{(r)}_{n+s}(W)$
implies $\psi\in \E^{(r)}_{n}(W)$. Let $a^{1}(x),\dots,a^{r}(x)\in
[U]$. There exists $k\in \N$ such that $\h^{k}a^{i}(x)\in U$ for
$i=1,\dots,r$.  As the sequence
$\h^{k}a^{1}(x),\dots,\h^{k}a^{r}(x)$ in $U$ is $\h$-adically
quasi-compatible, for every $n\in \N$, there exists $0\ne p(x,y)\in
\C[x,y]$ such that
$$\h^{rk}\left(\prod_{1\le i<j\le r}p(x_{i},x_{j})\right)
a^{1}(x_{1})\cdots a^{r}(x_{r}) \in \E^{(r)}_{n+rk}(W),$$ which
gives
$$\left(\prod_{1\le i<j\le r}p(x_{i},x_{j})\right)a^{1}(x_{1})\cdots a^{r}(x_{r})
\in \E^{(r)}_{n}(W).$$ This proves that the sequence
$a^{1}(x),\dots,a^{r}(x)$ is $\hbar$-adically quasi-compatible.
Therefore, $[U]$ is $\hbar$-adically quasi-compatible.

Assume that $U$ is $Y_{\E}$-closed. Let $a(x),b(x)\in [U],\; m\in
\Z$. By definition, there exists $k\in \N$ such that $\h^{k}a(x),\;
\h^{k}b(x)\in U$. Then
$$\h^{2k}(a(x)_{m}b(x))=(\h^{k}a(x))_{m}(\h^{k}b(x))\in U.$$
Thus $a(x)_{m}b(x)\in [U]$. This proves that $[U]$ is
$Y_{\E}$-closed.
\end{proof}

\bt{tmain-canonical} Let $K$ be a maximal $\hbar$-adically
quasi-compatible $\C[[\h]]$-submodule of $\E_{\h}(W)$. Then $[K]=K$,
$K$ is $\hbar$-adically topologically free and $Y_{\E}$-closed.
Furthermore, $(K,Y_{\E},1_{W})$ carries the structure of an
$\h$-adic nonlocal vertex algebra with $W$ as a quasi module with
$Y_{W}(\alpha(x),x_{0})=\alpha(x_{0})$ for $\alpha(x)\in K$. If $K$
is $\hbar$-adically compatible, $W$ is a module (instead of a quasi
module). \et

\begin{proof} By Lemma \ref{lpre}, $[K]$ is $\hbar$-adically quasi-compatible.
As $K$ is maximal, we have $[K]\subset K$. Thus $[K]=K$. Let
$a(x),b(x)\in K,\; m\in \Z$. For every $n\in \N$, $\pi_{n}(K)$ is
quasi-compatible, then by Theorem \ref{tgeneratingthem},
$\pi_{n}(K)$ generates a nonlocal vertex algebra $\<\pi_{n}(K)\>$
over $\C$ and we have
$$\pi_{n}(K+\C a(x)_{m}b(x))
=\pi_{n}(K)+\C\pi_{n}(a(x))_{m}\pi_{n}(b(x))\subset \<
\pi_{n}(K)\>,$$ a quasi-compatible $\C$-subspace of $\E(W/\h^{n}W)$.
This proves that $K+\C a(x)_{m}b(x)$ is $\hbar$-adically
quasi-compatible in $\E_{\hbar}(W)$. Again, with $K$ maximal, we
have $a(x)_{m}b(x)\in K$. Thus $K$ is $Y_{\E}$-closed. Now, let
$\{\psi_{m}(x)\}$ be a sequence in $K$, satisfying the condition
that for any $r\ge 0$, there exists $k\ge 0$ such that
$\psi_{m}(x)-\psi_{n}(x)\in \h^{r}K$ whenever $m,n\ge k$. Since
$\E_{\h}(W)$ is $\hbar$-adically complete, the sequence
$\{\psi_{m}(x)\}$ has a limit, say $\psi(x)$, in $\E_{\h}(W)$. For
any $n\in \N$, there exists $m\in \N$ such that
$\psi_{m}(x)-\psi(x)\in \h^{n}\E_{\h}(W)$, which implies
$\pi_{n}(\psi_{m}(x))=\pi_{n}(\psi(x))$. Thus
$\pi_{n}(\psi(x))=\pi_{n}(\psi_{m}(x))\in \pi_{n}(K)$. Consequently,
$\pi_{n}(K+\C \psi(x))\subset \pi_{n}(K)$, which is
quasi-compatible. This proves that $K+\C \psi(x)$ is
$\hbar$-adically quasi-compatible in $\E_{\hbar}(W)$ and then it
follows that $\psi(x)\in K$. Thus $K$ is $\h$-adically complete, so
that it is topologically free. Now, in view of Proposition
\ref{phvoa-1}, $K$ is an $\h$-adic nonlocal vertex algebra with $W$
as a quasi module. Furthermore, $W$ is a module if $K$ is
$\h$-adically compatible.
\end{proof}

Now, let $U$ be an $\hbar$-adically quasi-compatible subset of
$\E_{\h}(W)$. In view of Zorn's lemma, there exists a maximal
$\hbar$-adically quasi-compatible $\C[[\h]]$-submodule $K$ of
$\E_{\h}(W)$, containing $U$ and $1_{W}$. Set
$U^{(1)}=\C[[\hbar]]U+\C[[\hbar]]1_{W}$. Define $U^{(2)}$ to be the
$\C[[\hbar]]$-span of the vectors $a(x)_{m}b(x)$ for $a(x),b(x)\in
U^{(1)},\; m\in \Z$.   Since $K$ is $Y_{\E}$-closed by Theorem
\ref{tmain-canonical}, $U^{(2)}\subset K$. Then $U^{(2)}$ is
$\hbar$-adically quasi-compatible. For $n\ge 1$, we inductively
define $U^{(n+1)}=(U^{(n)})^{(2)}$. In this way, we obtain an
increasing sequence of $\hbar$-adically quasi-compatible
$\C[[\hbar]]$-submodules:
$$U^{(1)}\subset U^{(2)}\subset U^{(3)}\subset\cdots.$$
 Set
\begin{eqnarray}
\<U\>^{o}=\{ a(x)\in \E_{\hbar}(W)\;|\; \hbar^{k}a(x)\in \cup_{n\ge
2}U^{(n)} \ \ \mbox{ for some }k\ge 1 \}.
\end{eqnarray}
That is, $\<U\>^{o}=[\cup_{n\ge 2}U^{(n)}]$. In view of Lemma
\ref{lsomeproperty} we have
\begin{eqnarray}
\<U\>^{o}\cap \hbar^{n}\E_{\hbar}(W)=\hbar^{n}\<U\>^{o}\ \ \ \mbox{
for }n\ge 1.
\end{eqnarray}
In particular, the induced topology of $\<U\>^{o}$ from
$\E_{\hbar}(W)$ coincides with the $\hbar$-adic topology of
$\<U\>^{o}$. Then we define $\<U\>$ to be the $\hbar$-adic
completion of $\<U\>^{o}$.

\bt{tmain-canonical-U} Let $U$ be an $\hbar$-adically
quasi-compatible subset of $\E_{\h}(W)$. Then $[\<U\>]=\<U\>$,
$\<U\>$ is topologically free, $\hbar$-adically quasi-compatible,
and $Y_{\E}$-closed, and $(\<U\>,Y_{\E},1_{W})$ carries the
structure of an $\h$-adic nonlocal vertex algebra and $W$ is a
faithful quasi $\<U\>$-module with
$Y_{W}(\alpha(x),x_{0})=\alpha(x_{0})$ for $\alpha(x)\in \<U\>$.
Furthermore, for any $\hbar$-adic nonlocal vertex subalgebra $V$ of
$\E_{\hbar}(W)$, containing $U$, such that $[V]=V$, we have
$\<U\>\subset V$. \et

\begin{proof} To prove $[\<U\>]=\<U\>$, let $a(x)\in
[\<U\>]$. By definition, $a(x)\in \E_{\hbar}(W)$ and there exists
$k\ge 0$ such that $\hbar^{k}a(x)\in \<U\>$. As $\<U\>$ is the
$\hbar$-adic completion of $\<U\>^{o}$, there exists a Cauchy
sequence $\{\psi_{m}(x)\}$ in $\<U\>^{o}$ with $\hbar^{k}a(x)$ as a
limit. There exists $r\ge 1$ such that
$$\psi_{m}(x)-\hbar^{k}a(x)\in
\hbar^{k}\E_{\hbar}(W)\ \ \ \mbox{ for }m\ge r.$$ Then
$\psi_{m}(x)\in \hbar^{k}\E_{\hbar}(W)$ for $m\ge r$. Set
$\psi_{m}(x)=\hbar^{k}\phi_{m}(x)$ for $m\ge r$ with $\phi_{m}(x)\in
\E_{\hbar}(W)$. Noticing that $[\<U\>^{o}]=\<U\>^{o}$, we have
$\phi_{m}(x)\in \<U\>^{o}$ for $m\ge r$. We see that
$\{\phi_{m}(x)\}_{m\ge r}$ is a Cauchy sequence in $\<U\>^{o}$ with
$a(x)$ as a limit. Thus $a(x)\in \<U\>$. This proves
$[\<U\>]=\<U\>$. As $\<U\>$ is torsion-free, separated, and
$\hbar$-adically complete by definition, $\<U\>$ is topologically
free. It follows from definition that $\cup_{n\ge 2}U^{(n)}$ is
$\hbar$-adically quasi-compatible and $Y_{\E}$-closed. By Lemma
\ref{lpre}, $\<U\>^{o}\;(=[\cup_{n\ge 2}U^{(n)}])$ is
$\hbar$-adically quasi-compatible and $Y_{\E}$-closed. Let
$\psi_{1}(x),\dots,\phi_{r}(x)$ be a sequence in $\<U\>$ and let $n$
be any positive integer. For $1\le i\le r$, there exists a sequence
$\{ \psi_{im}(x)\}$ in $\<U\>^{o}$, which converges to
$\psi_{i}(x)$. Let $k$ be a positive integer such that
$$\psi_{im}(x)-\psi_{i}(x)\in \hbar^{n}\E_{\hbar}(W)
\ \ \ \mbox{ for }1\le i\le r,\; m\ge k.$$ As $\<U\>^{o}$ is
$\hbar$-adically quasi-compatible, there exists $0\ne p(x,y)\in
\C[x,y]$ such that
$$\pi_{n}\left(\coprod_{1\le i<j\le
r}p(x_{i},x_{j})\right)\phi_{1k}(x_{1})\cdots \phi_{rk}(x_{r}) \in
\Hom ((W/\hbar^{n}W), (W/\hbar^{n}W)((x_{1},\dots,x_{r}))).$$ Then
$$\pi_{n}\left(\coprod_{1\le i<j\le
r}p(x_{i},x_{j})\right)\phi_{1}(x_{1})\cdots \phi_{r}(x_{r}) \in
\Hom ((W/\hbar^{n}W), (W/\hbar^{n}W)((x_{1},\dots,x_{r}))).$$ This
proves that $\psi_{1}(x),\dots,\phi_{r}(x)$ is $\hbar$-adically
quasi-compatible. It follows from Lemma \ref{lbracket-closed} that
$\<U\>$ is $Y_{\E}$-closed.   Now, by Proposition \ref{phvoa-1},
$(\<U\>,Y_{\E},1_{W})$ carries the structure of an $\h$-adic
nonlocal vertex algebra and $W$ is a faithful quasi $\<U\>$-module.
Let $V$ be an $\hbar$-adic nonlocal vertex subalgebra of
$\E_{\hbar}(W)$ satisfying the condition that $U\subset V$ and
$[V]=V$. It is straightforward to see that $\<U\>\subset V$.
\end{proof}

For a topologically free $\C[[\hbar]]$-module $W$, we say a subset
$T$ {\em spans $W$ $\hbar$-adically}  if $W=(\C T)[[\hbar]]'$. The
following is an $\h$-adic version of (\cite{li-qva1}, Theorem 6.3)
which is an analogue of a theorem of Frenkel-Kac-Radual-Wang
\cite{fkrw} and of Muerman-Primc \cite{mp} (cf. \cite{ll}):

\bt{tfkrw-mp} Let $V$ be a topologically free $\C[[\hbar]]$-module,
$U$ a subset of $V$, ${\bf 1}$ a vector in $V$, $\D$ a
$\C[[\hbar]]$-module endomorphism of $V$, and $Y^{0}$ a map
$$Y^{0}:U \rightarrow \E_{\hbar}(V);\ \
u\mapsto Y^{0}(u,x)=u(x)=\sum_{m\in \Z}u_{m}x^{-m-1}.$$ Assume that
all the following conditions hold:
\begin{eqnarray}
& &\D {\bf 1}=0,\\
& &Y^{0}(u,x){\bf 1}\in V[[x]]\;\;\mbox{ and }\;\;\;
\lim_{x\rightarrow 0}Y^{0}(u,x){\bf 1}=u,\\
& &[\D,Y^{0}(u,x)]={d\over dx}Y^{0}(u,x)\;\;\;\mbox{ for }u\in U,
\end{eqnarray}
$U(x)=\{ u(x)\;|\; u\in U\}$ is $\hbar$-adically compatible, and
 $V$ is $\hbar$-adically spanned by vectors
\begin{eqnarray}
u^{(1)}_{m_{1}}\cdots u^{(r)}_{m_{r}}{\bf 1}
\end{eqnarray}
for $r\ge 0,\; u^{(i)}\in U,\; m_{i}\in \Z$. In addition we assume
that there exists a $\C[[\h]]$-module morphism  $\psi$ {}from $V$ to
$\<U(x)\>\subset \E_{\hbar}(V)$ such that $\psi({\bf 1})=1_{V}$ and
\begin{eqnarray}
\psi(u_{m}v)=u(x)_{m}\psi(v)\;\;\;\mbox{ for }u\in U,\; v\in V,\;
m\in \Z.
\end{eqnarray}
Then $Y^{0}$ extends uniquely  to a $\C[[\hbar]]$-module map $Y$
from $V$ to $\E_{\hbar}(V)$ such that $(V,Y,{\bf 1})$ carries the
structure of an $\h$-adic nonlocal vertex algebra. \et

\begin{proof} The uniqueness is obvious, so it remains to establish the existence.
 As $U(x)$ is an $\hbar$-adically compatible
 subset of $\E_{\hbar}(V)$, by Theorem \ref{tmain-canonical-U}
we have an $\h$-adic nonlocal vertex algebra $\<U(x)\>$ with $V$ as
a faithful $\<U(x)\>$-module where
$Y_{V}(\alpha(x),x_{0})=\alpha(x_{0})$ for $\alpha(x)\in \<U(x)\>$.
For $u\in U$, we have
\begin{eqnarray*}
& &Y_{V}(u(x),x_{0}){\bf 1}=u(x_{0}){\bf 1}\in V[[x_{0}]],\\
& &[\D,Y_{V}(u(x),x_{0})]=[\D,Y^{0}(u,x_{0})]
=\frac{d}{dx_{0}}Y^{0}(u,x_{0}) =\frac{d}{dx_{0}}Y_{V}(u(x),x_{0}).
\end{eqnarray*}
By Lemma \ref{lvacuum-like}, the map $\phi$ from $\<U(x)\>$ to $V$,
defined by $\phi(\alpha(x))=\Res_{x}x^{-1}\alpha(x){\bf 1}$, is a
$\<U(x)\>$-module homomorphism. We see that $\phi(1_{V})={\bf 1}$
and that for $u\in U,\; \alpha(x)\in \<U(x)\>$,
$$\phi(Y_{\E}(u(x),x_{0})\alpha(x))=Y_{V}(u(x),x_{0})\phi(\alpha(x))
=u(x_{0})\phi(\alpha(x)),$$ which amounts to
$$\phi(u(x)_{m}\alpha(x))=u_{m}\phi(\alpha(x))\ \ \ \mbox{ for }m\in \Z.$$
It follows that $\phi\circ \psi=1_{V}$. Thus $\psi$ is a
$\C[[\hbar]]$-module isomorphism from $V$ onto $\psi(V)\subset
\<U(x)\>$. For $u\in U$, we have
$$\psi (u)=\psi(u_{-1}{\bf 1})=u(x)_{-1}1_{V}=u(x)
\ \ \mbox{ and }\ \ \ \phi(u(x))=\phi (\psi(u))=u.$$ Inside
$\<U(x)\>$, $\psi(V)$ is $\hbar$-adically spanned by vectors
\begin{eqnarray*}
u^{(1)}(x)_{m_{1}}\cdots u^{(r)}(x)_{m_{r}}1_{V}
\end{eqnarray*}
for $r\ge 0,\; u^{(i)}\in U,\; m_{i}\in \Z$.

For $a\in V$, we define $Y(a,x)\in (\End V)[[x,x^{-1}]]]$ by
$$Y(a,x_{0})b=\phi \left(Y_{\E}(\psi(a)(x),x_{0})\psi(b)(x)\right)
\  \mbox{ for }b\in V.$$ As $\phi$ is a $\<U(x)\>$-module
homomorphism with $\phi \circ \psi=1$, we have
$$Y(a,x_{0})b=Y_{V}(\psi(a)(x),x_{0})b=\psi(a)(x_{0})b,$$
so that $Y(a,x)=\psi(a)(x)\in \E_{\hbar}(V)$. In particular, for
$u\in U$,
$$Y(u,x_{0})=Y_{V}(\psi(u),x_{0})=Y_{V}(u(x),x_{0})=u(x_{0})=Y^{0}(u,x_{0}),$$
so the map $Y$ extends $Y_{0}$.  For $v\in V$, we have
$$Y({\bf 1},x_{0})v=Y_{V}(1_{V},x_{0})v=1_{V}(v)=v,$$
$$Y(v,x_{0}){\bf 1}=Y_{V}(\psi(v)(x),x_{0}){\bf 1}\in V[[x_{0}]],$$ and
$$\lim_{x_{0}\rightarrow 0}Y(v,x_{0}){\bf 1}
=\lim_{x_{0}\rightarrow 0}Y_{V}(\psi(v)(x),x_{0}){\bf 1}=\phi (\psi
(v))=v.$$ To prove that $(V,Y,{\bf 1})$ is an $\hbar$-adic nonlocal
vertex algebra, we show that for every positive integer $n$,
$V/\hbar^{n} V$ with the reduced structures is a nonlocal vertex
algebra over $\C$. For any $n\ge 1$,  $\<U(x)\>/\hbar^{n}\<U(x)\>$
is a nonlocal vertex algebra over $\C$ and $\phi$ reduces to a
homomorphism $\bar{\phi}$ {}from $\<U(x)\>/\hbar^{n}\<U(x)\>$ to
$V/\hbar^{n}V$. On the other hand, the map $\psi$ reduces to a map
$\bar{\psi}$ {}from $V/\hbar^{n}V$ to $\<U(x)\>/\hbar^{n}\<U(x)\>$
such that $\bar{\phi}\circ \bar{\psi}=1$. We see that the image
$\overline{\psi(V)}$ of $\psi(V)$ in $\<U(x)\>/\hbar^{n}\<U(x)\>$ is
a nonlocal vertex subalgebra and $\overline{\psi(V)}\simeq
V/\hbar^{n}V$ through the maps $\bar{\phi}$ and $\bar{\psi}$. It
follows that $V/\hbar^{n}V$ is a nonlocal vertex algebra over $\C$.
Therefore, $(V,Y,{\bf 1})$ is an $h$-adic nonlocal vertex algebra.
This establishes the existence, concluding the proof.
\end{proof}

\bd{dZx12} {\em Let $W$ be a topologically free $\C[[\hbar]]$-module
as before. We define a $\C[[\hbar]]$-module map
\begin{eqnarray*}
Z(x_{1},x_{2}): \E_{\hbar}(W)\hat{\otimes}
\E_{\hbar}(W)\hat{\otimes} \C((x))[[\hbar]]\rightarrow (\End
W)[[x_{1}^{\pm 1},x_{2}^{\pm 1}]]
\end{eqnarray*}
by
\begin{eqnarray}
Z(x_{1},x_{2})(a(x)\otimes b(x)\otimes f(x))
=f(x_{1}-x_{2})a(x_{1})b(x_{2}).
\end{eqnarray}} \ed

\bl{literate-formula} Let $a(x),b(x)\in \E_{\hbar}(W),\; B(x)\in
\E_{\hbar}(W)\hat{\otimes}\E_{\hbar}(W)\hat{\otimes}\C((x))[[\hbar]]$
such that
$$a(x_{1})b(x_{2})\sim Z(x_{2},x_{1})(B(x)).$$
Then $(a(x),b(x))$ is $\hbar$-adically compatible and
\begin{eqnarray}\label{eabiterate-formula}
&&Y_{\E}(a(x),x_{0})b(x)\nonumber\\
&=&\Res_{x_{1}}x_{0}^{-1}\delta\left(\frac{x_{1}-x}{x_{0}}\right)a(x_{1})b(x)
-x_{0}^{-1}\delta\left(\frac{x-x_{1}}{-x_{0}}\right)Z(x,x_{1})(B(x)).
\end{eqnarray}
\el

\begin{proof} By definition, for any positive integer $n$, there exists
$k\in \N$ such that
$$(x_{1}-x_{2})^{k}\pi_{n}(a(x_{1}))\pi_{n}(b(x_{2}))
=(x_{1}-x_{2})^{k}\pi_{n}^{(2)}(Z(x_{2},x_{1})(B(x))).$$ {}From
\cite{li-qva1} we have
\begin{eqnarray*}
&&Y_{\E}(\pi_{n}(a(x)),x_{0})\pi_{n}(b(x))\nonumber\\
&=&\Res_{x_{1}}x_{0}^{-1}\delta\left(\frac{x_{1}-x}{x_{0}}\right)\pi_{n}(a(x_{1}))\pi_{n}(b(x))
-x_{0}^{-1}\delta\left(\frac{x-x_{1}}{-x_{0}}\right)\pi_{n}^{(2)}Z(x,x_{1})(B(x)).
\end{eqnarray*}
Then it follows.
\end{proof}

\bd{dYEx12} {\em Let $V$ be an $\hbar$-adic nonlocal vertex
subalgebra of $\E_{\hbar}(W)$. Define  a $\C[[\hbar]]$-module map
$Y_{\E}(x_{2},x_{1})$ from $V\hat{\otimes} V\hat{\otimes}
\C((x))[[\hbar]]$ to $(\End V)[[x_{1}^{\pm 1},x_{2}^{\pm 1}]]$ by
\begin{eqnarray}
Y_{\E}(x_{2},x_{1})(a(x)\otimes b(x)\otimes f(x))
=f(x_{2}-x_{1})Y_{\E}(a(x),x_{2})Y_{\E}(b(x),x_{1})
\end{eqnarray}
for $a(x),b(x)\in V,\; f(x)\in \C((x))[[\hbar]]$.} \ed

\bl{ph-transp-slocal} Let $V$ be an $\h$-adic nonlocal vertex
subalgebra of $\E_{\hbar}(W)$, and let
$$u(x),v(x)\in V, \ \ A(x)\in V\hat{\otimes} V\hat{\otimes} \C((x))[[\hbar]].$$
Suppose that
\begin{eqnarray}
u(x_{1})v(x_{2}) \sim Z(x_{2},x_{1})(A(x))
\end{eqnarray}
in $(\End W)[[x_{1}^{\pm 1},x_{2}^{\pm 1}]]$. Then
\begin{eqnarray}
Y_{\E}(u(x),x_{1})Y_{\E}(v(x),x_{2})\sim Y_{\E}(x_{2},x_{1})(A(x))
\end{eqnarray}
in $(\End V)[[x_{1}^{\pm 1},x_{2}^{\pm 1}]]$ and
\begin{eqnarray}
&&x_{0}^{-1}\delta\left(\frac{x_{1}-x_{2}}{x_{0}}\right)Y_{\E}(u(x),x_{1})Y_{\E}(v(x),x_{2})
-x_{0}^{-1}\delta\left(\frac{x_{2}-x_{1}}{-x_{0}}\right)Y_{\E}(x_{2},x_{1})(A(x))\nonumber\\
&&\hspace{2.5cm}=x_{2}^{-1}\delta\left(\frac{x_{1}-x_{0}}{x_{2}}\right)
Y_{\E}(Y_{\E}(u(x),x_{0})v(x),x_{2}).
\end{eqnarray}
 \el

\begin{proof} Let $n$ be a positive integer.
We have a nonlocal vertex algebra $\pi_{n}(V)\subset
\E(W/\hbar^{n}W)$ over $\C$ with $\ker \pi_{n}=V\cap
\hbar^{n}\E_{\hbar}(W)\;(=\hbar^{n}V)$. {}From assumption, there
exists a nonnegative integer $k$ such that
$$(x_{1}-x_{2})^{k}\pi_{n}(u(x_{1}))\pi_{n}(v(x_{2}))=
(x_{1}-x_{2})^{k}Z(x_{2},x_{1})(\pi_{n}(A(z))).$$ {}From
\cite{li-qva1} (Proposition 3.13), we have
$$(x_{1}-x_{2})^{k}Y_{\E}(\pi_{n}(u(x)),x_{1})Y_{\E}(\pi_{n}(v(x)),x_{2})
=(x_{1}-x_{2})^{k}Y_{\E}(x_{2},x_{1})\pi_{n}(A(z))$$ as desired.
\end{proof}

\bd{dh-Slocal} {\em A subset $U$ of $\E_{\h}(W)$ is said to be {\em
$\hbar$-adically $\S$-local} if for any $a(x),b(x)\in U$, there
exists $A(x)\in (\C U\otimes \C U\otimes \C((x)))[[\hbar]]$
 such that
\begin{eqnarray}
a(x_{1})b(x_{2}) \sim Z(x_{2},x_{1})(A(x)).
\end{eqnarray}}
\ed

\bl{lclosedness} Every $\hbar$-adically $\S$-local subset of
$\E_{\hbar}(W)$ is $\hbar$-adically compatible. \el

\begin{proof} Let $U$ be an $\hbar$-adically $\S$-local subset of $\E_{\hbar}(W)$.
For every positive integer $n$, we see that $\pi_{n}(\C[[\hbar]]U)$
is an $\S$-local subset of $\E(W/\hbar^{n}W)$, so that
$\pi_{n}(\C[[\hbar]]U)$ is compatible. By definition, $\C[[\hbar]]U$
is an $\hbar$-adically compatible subset of $\E_{\hbar}(W)$. Thus
$U$ is an $\hbar$-adically compatible subset.
\end{proof}

The following is a refinement of Theorem \ref{tfkrw-mp} (cf.
\cite{li-qva2}, Theorem 2.9):

\bt{tqva-construction} Let $V$ be a topologically free
$\C[[\hbar]]$-module, $U$ a subset of $V$, ${\bf 1}$ a vector in
$V$, and $Y^{0}$ a map
$$Y^{0}:U \rightarrow \E_{\hbar}(V);\ \
u\mapsto Y^{0}(u,x)=u(x)=\sum_{m\in \Z}u_{m}x^{-m-1}.$$ Assume that
all the following conditions hold:
\begin{eqnarray}
Y^{0}(u,x){\bf 1}\in V[[x]]\;\;\mbox{ and }\;\;\; \lim_{x\rightarrow
0}Y^{0}(u,x){\bf 1}=u \;\;\;\mbox{ for }u\in U,
\end{eqnarray}
$U(x)=\{ u(x)\;|\; u\in U\}$ is $\hbar$-adically $\S$-local, and
 $V$ is $\hbar$-adically spanned by vectors
\begin{eqnarray}
u^{(1)}_{m_{1}}\cdots u^{(r)}_{m_{r}}{\bf 1}
\end{eqnarray}
for $r\ge 0,\; u^{(i)}\in U,\; m_{i}\in \Z$. In addition we assume
that there exists a $\C[[\h]]$-module morphism  $\psi$ {}from $V$ to
$\<U(x)\>\subset \E_{\hbar}(V)$ such that $\psi({\bf 1})=1_{V}$ and
\begin{eqnarray}
\psi(u_{m}v)=u(x)_{m}\psi(v)\;\;\;\mbox{ for }u\in U,\; v\in V,\;
m\in \Z.
\end{eqnarray}
Then the map $Y^{0}$ extends uniquely to a $\C[[\hbar]]$-module map
$Y$ {}from $V$ to $\E_{\hbar}(V)$ such that $(V,Y,{\bf 1})$ carries
the structure of an $\h$-adic weak quantum vertex algebra.
 \et

\begin{proof} We shall slightly modify the proof of Theorem \ref{tfkrw-mp}.
By Lemma \ref{lclosedness}, $U(x)$ is $\hbar$-adically compatible,
so it generates an $\hbar$-adic nonlocal vertex algebra $\<U(x)\>$.
Set
$$K=(\cup_{n\ge 1}U(x)^{(n)})[[\hbar]]'\subset \<U(x)\>\subset \E_{\hbar}(V).$$
By Proposition \ref{pvacuum-generating},  we have a
$\C[[\hbar]]$-map $\phi: K\rightarrow V$ such that
$$\phi(1_{V})={\bf 1},\ \ \ \phi (u(x)_{n}a(x))=u_{n}\phi(a(x))
\ \ \ \mbox{ for }u\in U,\; n\in \Z,\; a(x)\in K.$$ Then continue
with the proof of Theorem \ref{tfkrw-mp} to see that  $Y^{0}$
extends uniquely to a $\C[[\hbar]]$-module map $Y$ from $V$ to
$\E_{\hbar}(V)$ such that $(V,Y,{\bf 1})$ carries the structure of
an $\h$-adic nonlocal vertex algebra. From Lemma
\ref{ph-transp-slocal}, $U$ is an $\hbar$-adically $\S$-local subset
of $V$ and then by Proposition \ref{pgenerator-Slocal}, $V$ is an
$\hbar$-adic weak quantum vertex algebra.
\end{proof}

Notice that compared with the corresponding theorem in \cite{fkrw}
and \cite{mp}, Theorems \ref{tfkrw-mp} and \ref{tqva-construction}
have an extra assumption on the existence of the map $\psi$. By
using classical linear algebra, it is not hard to see that in the
general noncommutative situation, an assumption like this is indeed
necessary. This assumption means that $V$ is a universal vacuum
module for a certain algebra.

The following results are companions in practical applications:

\bl{lpre-hadic-slocal} Let $U$ be a subset of $\E_{\hbar}(W)$
satisfying the condition that for $a(x),b(x)\in U$, there exist
$B(z)\in (\C U\otimes \C U\otimes \C((x)))[[\hbar]]$ and
$p(x,\hbar)\in \C[x,\hbar]$ with $p(x,0)\ne 0$ such that
\begin{eqnarray}\label{ecommon-quantity}
p(x_{1}-x_{2},\hbar)a(x_{1})b(x_{2})
=p(x_{1}-x_{2},\hbar)Z(x_{2},x_{1})(B(x)).
\end{eqnarray}
Then $U$ is $\hbar$-adically $\S$-local. \el

\begin{proof} Let $a(x),b(x)\in U$. By assumption there exist $B$ and
$p(x,\hbar)$ with all the assumed properties. With $p(x,0)\ne 0$, we
have $p(x,\hbar)=f(x)-\hbar g(x,\hbar)$, where $0\ne f(x)\in
\C[x],\; g(x,\hbar)\in \C[x,\hbar]$. Expand $p(x,\hbar)^{-1}$ in the
nonnegative powers of $\hbar$ as follows
$$p(x,\hbar)^{-1}=\sum_{i\ge 0}f(x)^{-1-i}g(x,\hbar)^{i}\hbar^{i}\in \C((x))[[\hbar]],$$
where $f(x)^{-i-1}$ is understood as an element of $\C((x))$. Let
$n$ be a positive integer. Then
$$p(x,\hbar)^{-1}\equiv
\sum_{i=0}^{n-1}f(x)^{-1-i}g(x,\hbar)^{i}\hbar^{i} \;\;\;(\mod\;
\h^{n}\C((x))[[\hbar]]).$$  Let $k$ be a nonnegative integer such
that $x^{k}f(x)^{-n}\in \C[[x]]$, so that
$$(x_{1}-x_{2})^{k}f(x_{1}-x_{2})^{-1-i}
=(-x_{2}+x_{1})^{k}f(-x_{2}+x_{1})^{-1-i}\ \ \mbox{ for all }
i=0,\dots,n-1.$$ Set $A=p(x_{1}-x_{2},\hbar)a(x_{1})b(x_{2})$, the
common quantity of both sides of (\ref{ecommon-quantity}). We have
\begin{eqnarray*}
& &(x_{1}-x_{2})^{k}a(x_{1})b(x_{2})\nonumber\\
&=&(x_{1}-x_{2})^{k}p(x_{1}-x_{2}, \h)^{-1}A\nonumber\\
&\equiv&\left((x_{1}-x_{2})^{k}\sum_{i=0}^{n-1}f(x_{1}-x_{2})^{-1-i}
g(x_{1}-x_{2},\hbar)^{i}\hbar^{i}\right)
A\;\;\;(\mod\; \h^{n})\nonumber\\
&=&(-x_{2}+x_{1})^{k}\sum_{i=0}^{n-1}
f(-x_{2}+x_{1})^{-1-i}g(-x_{2}+x_{1},\hbar)^{i}\hbar^{i}
A\nonumber\\
&\equiv & (-x_{2}+x_{1})^{k}\sum_{i\ge 0}
f(-x_{2}+x_{1})^{-1-i}g(-x_{2}+x_{1},\hbar)^{i}\hbar^{i}
A\;\;\;(\mod\; \h^{n})\nonumber\\
&=&(-x_{2}+x_{1})^{k}Z(x_{2},x_{1})B(x).
\end{eqnarray*}
This proves
\begin{eqnarray*}
a(x_{1})b(x_{2})\sim Z(x_{2},x_{1})B(x).
\end{eqnarray*}
Thus $U$ is $\hbar$-adically $\S$-local.
\end{proof}

\bp{phadic-slocal} Let $V$ be an $\hbar$-adic nonlocal vertex
subalgebra of $\E_{\hbar}(W)$. Suppose that $a(x),b(x)\in V$,
$B(x)\in V\hat{\otimes} V\hat{\otimes} \C((x))[[\hbar]]$, and
$p(x,\hbar)\in \C[x,\hbar]$ with $p(x,0)\ne 0$ satisfy
\begin{eqnarray}
p(x_{1}-x_{2},\hbar)a(x_{1})b(x_{2})
=p(x_{1}-x_{2},\hbar)Z(x_{2},x_{1})(B(x)).
\end{eqnarray}
Then
\begin{eqnarray}
p(x_{1}-x_{2},\hbar)Y_{\E}(a(x),x_{1})Y_{\E}(b(x),x_{2})
=p(x_{1}-x_{2},\hbar)Y_{\E}(x_{2},x_{1})(B(x)).
\end{eqnarray}
 \ep

\begin{proof} In view of Lemma \ref{lpre-hadic-slocal}, we have
$$a(x_{1})b(x_{2})\sim Z(x_{2},x_{1})(B(x)).$$
Furthermore, by Lemma \ref{ph-transp-slocal} we have
\begin{eqnarray*}
Y_{\E}(a(x),x_{1})Y_{\E}(b(x),x_{2})\sim Y_{\E}(x_{2},x_{1})(B(x)).
\end{eqnarray*}
Then the following Jacobi identity holds:
\begin{eqnarray*}
& &x_{0}^{-1}\delta\left(\frac{x_{1}-x_{2}}{x_{0}}\right)
Y_{\E}(a(x),x_{1})Y_{\E}(b(x),x_{2})
-x_{0}^{-1}\delta\left(\frac{x_{2}-x_{1}}{-x_{0}}\right)Y_{\E}(x_{2},x_{1})(B(x))\\
&&\ \ \ \ \ \
=x_{1}^{-1}\delta\left(\frac{x_{2}+x_{0}}{x_{1}}\right)
Y_{\E}(Y_{\E}(a(x),x_{0})b(x),x_{2}).
\end{eqnarray*}
By Proposition \ref{panother-version} we have
$$p(x_{0},\hbar)Y_{\E}(a(x),x_{0})b(x)=\left(p(x_{1}-x,\hbar)a(x_{1})b(x)\right)|_{x_{1}=x+x_{0}},$$
which involves only nonnegative integer powers of $x_{0}$.
Multiplying the both sides of the Jacobi identity by
$p(x_{0},\hbar)$, and then applying $\Res_{x_{0}}$ we obtain the
desired identity.
\end{proof}

We also have:

\bp{ph-transp} Let $W$ be given as before and let $V$ be an
$\h$-adic nonlocal vertex subalgebra of $\E_{\hbar}(W)$ such that
$V$ is $\hbar$-adically compatible, and let
$$m\in \Z,\; u(x),v(x), c^{0}(x),c ^{1}(x),\dots \in V,
\ \ A(x)\in V\hat{\otimes} V\hat{\otimes} \C((x))[[\hbar]]$$
satisfying the condition that for every positive integer $n$, there
exists a nonnegative integer $r$ such that $c^{j}(x)\in \hbar^{n}V$
for $j\ge r$. Suppose that
\begin{eqnarray}
& &(x_{1}-x_{2})^{m}u(x_{1})v(x_{2})
-(-x_{2}+x_{1})^{m}Z(x_{2},x_{1})(A(x))
\nonumber\\
& &\ \ \ \ \ \ = \sum_{j\ge 0}c^{j}(x_{2})\frac{1}{j!}
\left(\frac{\partial}{\partial x_{2}}\right)^{j}
x_{2}^{-1}\delta\left(\frac{x_{1}}{x_{2}}\right).
\end{eqnarray}
Then
\begin{eqnarray}\label{esecond}
& &(x_{1}-x_{2})^{m}Y_{\E}(u(x),x_{1})Y_{\E}(v(x),x_{2})
-(-x_{2}+x_{1})^{m}Y_{\E}(x_{2},x_{1})(A(x))
\nonumber\\
& &\ \ \ \ \ \ \ =\sum_{j\ge 0}Y_{\E}(c^{j}(x),x_{2})\frac{1}{j!}
\left(\frac{\partial}{\partial x_{2}}\right)^{j}
x_{2}^{-1}\delta\left(\frac{x_{1}}{x_{2}}\right).
\end{eqnarray}
\ep

\begin{proof} Since $V$ is $\hbar$-adically compatible,
$W$ is a faithful $V$-module with
$Y_{W}(\alpha(x),x_{0})=\alpha(x_{0})$ for $\alpha(x)\in V$. Then it
follows immediately {}from Proposition
\ref{phmodule-algebra-relations}.
\end{proof}

\section{$\hbar$-deformations of quantum vertex algebras $V_{\bf Q}$}
In this section we construct a family of $\hbar$-adic quantum vertex
algebras as $\hbar$-deformations of certain quantum vertex algebras
which were studied in \cite{kl}. One special case gives rise to an
$\hbar$-deformed $\beta\gamma$-system, while another special case
gives rise to an $\hbar$-deformation of the vertex operator
superalgebra $V_{L}$ associated to the rank-one lattice $L=\Z
\alpha$ with $\<\alpha,\alpha\>=1$. We essentially deal with the
same algebras as in \cite{kl} with a formal parameter $\hbar$,
instead of a nonzero complex number $q$.

We first recall the quantum vertex algebras of Zamolodchikov-Faddeev
type from \cite{kl}. Let $l$ be a positive integer and let ${\bf
Q}=(q_{ij})_{i,j=1}^{l}$ be a complex matrix such that
\begin{eqnarray}\label{eq-skewsymmetry}
q_{ij}q_{ji}=1\ \ \ \mbox{ for }1\le i,j\le l.
\end{eqnarray}
Define $\A_{\bf Q}$ to be the associative algebra with identity
(over $\C$) with generators
$$X_{i,n},\; Y_{i,n}\;\;\ \ (i=1,\dots,l,\; n\in \Z),$$
subject to relations
\begin{eqnarray}
& &X_{i,m}X_{j,n}=q_{ij}X_{j,n}X_{i,m},
\ \ \ \ Y_{i,m}Y_{j,n}=q_{ij}Y_{j,n}Y_{i,m},\nonumber\\
&&X_{i,m}Y_{j,n}-q_{ji}Y_{j,n}X_{i,m}=\delta_{i,j}\delta_{m+n+1,0}
\end{eqnarray}
for $1\le i,j\le l,\; m,n\in \Z$. A vector $w$ in an $\A_{\bf
Q}$-module is called a {\em vacuum vector} if
$$X_{i,n}w=Y_{i,n}w=0\ \ \ \mbox{ for }1\le i\le l,\; n\ge 0,$$ and an
$\A_{\bf Q}$-module $W$ equipped with a vacuum vector that generates
$W$ is called a {\em vacuum $\A_{\bf Q}$-module}.

Let $J_{\bf Q}$ be the left ideal of $\A_{\bf Q}$, generated by
$X_{i,n},\; Y_{i,n}$ $(1\le i\le l,\; n\ge 0)$, that is,
\begin{eqnarray*}
J_{\bf Q}=\sum_{i=1}^{l}\sum_{n\ge 0}(\A_{\bf Q}X_{i,n}+\A_{\bf
Q}Y_{i,n}).
\end{eqnarray*}
Set
\begin{eqnarray}
V_{\bf Q}=\A_{\bf Q}/J_{\bf Q},
\end{eqnarray}
a left $\A_{\bf Q}$-module, and set
$${\bf 1}=1+J_{\bf Q}\in V_{\bf Q}.$$
Then ${\bf 1}$ is a vacuum vector and $V_{\bf Q}$ equipped with
${\bf 1}$ is a vacuum $\A_{\bf Q}$-module. For $1\le i\le l$, set
\begin{eqnarray}
u^{(i)}=X_{i,-1}{\bf 1},\ \ v^{(i)}=Y_{i,-1}{\bf 1}\in V_{\bf Q}
\end{eqnarray}
and set
\begin{eqnarray}
X_{i}(x)=\sum_{n\in \Z}X_{i,n}x^{-n-1},\ \ Y_{i}(x)=\sum_{n\in
\Z}Y_{i,n}x^{-n-1}\in \A_{\bf Q}[[x,x^{-1}]].
\end{eqnarray}
It was proved therein that there exists a unique quantum vertex
algebra structure on $V_{\bf Q}$ with ${\bf 1}$ as the vacuum vector
and with
$$Y(u^{(i)},x)=X_{i}(x),\ \
Y(v^{(i)},x)=Y_{i}(x)\ \ \ \mbox{ for }1\le i\le l.$$ It was also
proved that $V_{\bf Q}$ is nondegenerate.

Next, we are going to construct a family of $\hbar$-adic quantum
vertex algebras by deforming $V_{\bf Q}$. For this purpose we shall
need the following notion (cf. \cite{li-pseudo}):

\bd{dhadicpseudo-end} {\em Let $V$ be a general nonlocal vertex
algebra. An {\em $\hbar$-adic pseudo-endomorphism} of $V$ is a
linear map
$$\Phi(x): V\rightarrow (V\otimes \C((x)))[[\hbar]]$$
satisfying the condition that $\Phi(x){\bf 1}={\bf 1}\otimes 1$,
\begin{eqnarray}\label{edef-pseudo}
\Phi(x_{1})Y(v,x_{2})=Y(\Phi(x_{1}-x_{2})v,x_{2})\Phi(x_{1}) \ \ \ \
\mbox{ for }v\in V,
\end{eqnarray}
where the map $Y$ is canonically extended. An $\hbar$-adic
pseudo-endomorphism $\Phi(x)$ is called a {\em pseudo-automorphism}
if there exists an $\hbar$-adic pseudo-endomorphism $\Psi(x)$ such
that $\Phi(x)\Psi(x)v=v=\Psi(x)\Phi(x)v$ for $v\in V$. We say that
$\hbar$-adic pseudo-endomorphisms $\Phi(x)$ and $\Psi(x)$ {\em
commute} if $\Phi(x_{1})\Psi(x_{2})=\Psi(x_{2})\Phi(x_{1})$. } \ed

The following is an $\hbar$-adic version of Lemma 3.15 of \cite{kl}:

\bl{lpseudo-auto-pseudo} Let ${\bf Q}=(q_{ij})$ be an $l\times l$
matrix as before and let $p_{1}(x,\hbar),\dots,p_{l}(x,\hbar)$ be
any sequence in $\C((x))[[\hbar]]$ with $p_{i}(x,0)\ne 0$ for $1\le
i\le l$ (so that $p_{i}(x,\hbar)$ are invertible). Then there exists
an $\hbar$-adic pseudo-automorphism $\Phi(x)$ of $V_{{\bf Q}}$ such
that
\begin{eqnarray*}
\Phi(x)(u^{(i)})=u^{(i)}\otimes p_{i}(x,\hbar),\ \ \ \
\Phi(x)(v^{(i)})=v^{(i)}\otimes p_{i}(x,\hbar)^{-1}\ \ \ \mbox{ for
}1\le i\le l.
\end{eqnarray*}
Furthermore, all such pseudo-automorphisms mutually commute.
 \el

\begin{proof} Note that $V_{\bf Q}\otimes \C((x))[[\hbar]]\subset (V_{\bf
Q}\otimes \C((x)))[[\hbar]]$. As $\C((x))[[\hbar]]$ is a commutative
associative algebra over $\C$ with $-\frac{d}{dx}$ as a derivation,
$\C((x))[[\hbar]]$ becomes a vertex algebra with $1$ as the vacuum
vector and with
$$Y(f,z)g=(e^{-z\frac{d}{dx}}f)g\ \ \ \ \mbox{ for }f,g\in
\C((x))[[\hbar]].$$ We then equip $V_{\bf Q}\otimes
\C((x))[[\hbar]]$ with the tensor product vertex algebra structure
with $Y_{ten}$ denoting the vertex operator map. Then
$$Y_{ten}(v\otimes f(x),z)=Y(v,z)\otimes f(x-z)
\ \ \ \mbox{ for }v\in V,\; f(x)\in \C((x))[[\hbar]].$$ For $A(x)\in
V\otimes \C((x))[[\hbar]]$, we have $Y_{ten}(A(x),z)=Y(A(x-z),z)$,
noting that as in (\ref{edef-pseudo}), $Y$ is
$\C((x))[[\hbar]]$-bilinear. An $\hbar$-adic pseudo-endomorphism
{}from $V_{\bf Q}$ to $V_{\bf Q}\otimes \C((x))[[\hbar]]$ exactly
amounts to a vertex algebra homomorphism.

It is straightforward to check that with $X_{i}(z)$ and $Y_{i}(z)$
($1\le i\le l$) acting on $V_{\bf Q}\otimes \C((x))[[\hbar]]$ as
$Y(u^{(i)}\otimes p_{i}(x,\hbar),z)$ and $Y(v^{(i)}\otimes
p_{i}(x,\hbar)^{-1},z)$, respectively, $V_{\bf Q}\otimes
\C((x))[[\hbar]]$ becomes an $\A_{\bf Q}$-module with ${\bf
1}\otimes 1$ as a vacuum vector. Since $V_{\bf Q}$ is a universal
vacuum $\A_{\bf Q}$-module, it follows that there exists an $\A_{\bf
Q}$-module homomorphism $\theta$ from $V_{\bf Q}$ to $V_{\bf
Q}\otimes \C((x))[[\hbar]]$ such that $\theta({\bf 1})={\bf
1}\otimes 1$. Because $V_{\bf Q}$ as a quantum vertex algebra is
generated by $u^{(i)},v^{(i)}$ ($1\le i\le l$), it follows that
$\theta$ is a vertex algebra homomorphism. We have
$$\theta (u^{(i)})=u^{(i)}\otimes p_{i}(x,\hbar),\ \ \
\theta_{i} (v^{(i)})=v^{(i)}\otimes p_{i}(x,\hbar)^{-1}\ \ \mbox{
for }1\le i\le l.$$ Denoting $\theta$ alternatively by $\Phi(x)$, we
see that $\Phi(x)$ satisfies all the properties.
\end{proof}

Set
\begin{eqnarray}
G(x,\hbar)=\left\{ \frac{p(x,\hbar)}{q(x,\hbar)}\;|\;
p(x,\hbar),q(x,\hbar)\in \C[x,\hbar] \ \ \mbox{ with
}p(x,0),q(x,0)\ne 0\right\},
\end{eqnarray}
an abelian group. We also consider $G(x,\hbar)$ as a
(multiplicative) subgroup of $\C((x))[[\hbar]]$.

\bt{tpre-diag} Let ${\bf Q}=(q_{ij})_{i,j=1}^{l}$ be given as before
and let
\begin{eqnarray}
p_{ij}(x,\hbar)\in G(x,\hbar)\subset \C((x))[[\hbar]]
\end{eqnarray}
such that $p_{ij}(x,0)=1$ for $1\le i,j\le l$. For $1\le i\le l$,
let $\Phi_{i}(x)$ be the pseudo-automorphism of $V_{\bf Q}$ such
that
\begin{eqnarray*}
& &\Phi_{i}(x)(u^{(j)})=u^{(j)}\otimes p_{ij}(x,\hbar),\ \ \ \
\Phi_{i}(x)(v^{(j)})=v^{(j)}\otimes p_{ij}(x,\hbar)^{-1} \ \ \
\mbox{ for }1\le j\le l,\ \ \ \ \
\end{eqnarray*}
obtained in Lemma \ref{lpseudo-auto-pseudo}. Then there exists a
unique $\hbar$-adic quantum vertex algebra structure on $V_{\bf
Q}[[\hbar]]$ with ${\bf 1}$ as the vacuum vector and with
$$Y_{\hbar}(u^{(i)},x)=Y(u^{(i)},x)\Phi_{i}(x),\ \ \ \
Y_{\hbar}(v^{(i)},x)=Y(v^{(i)},x)\Phi_{i}(x)^{-1} \ \ \ \mbox{ for
}1\le j\le l.$$ Furthermore, $V_{\bf Q}[[\hbar]]$ is non-degenerate
and generated by $u^{(i)},v^{(i)}$ ($1\le i\le l$), and the
following relations hold for $1\le i,j\le l$:
\begin{eqnarray*}
&&p_{ij}(x_{1}-x_{2},\hbar)^{-1}Y_{\hbar}(u^{(i)},x_{1})Y_{\hbar}(u^{(j)},x_{2})
=q_{ji}p_{ji}(x_{2}-x_{1},\hbar)^{-1}
Y_{\hbar}(u^{(j)},x_{2})Y_{\hbar}(u^{(i)},x_{1}),\\
&&p_{ij}(x_{1}-x_{2},\hbar)^{-1}Y_{\hbar}(v^{(i)},x_{1})Y_{\hbar}(v^{(j)},x_{2})
=q_{ij}p_{ji}(x_{2}-x_{1},\hbar)^{-1}
Y_{\hbar}(v^{(j)},x_{2})Y_{\hbar}(v^{(i)},x_{1}),\\
&&p_{ij}(x_{1}-x_{2},\hbar)Y_{\hbar}(u^{(i)},x_{1})Y_{\hbar}(v^{(j)},x_{2})
-q_{ji}p_{ji}(x_{2}-x_{1},\hbar)Y_{\hbar}(v^{(j)},x_{2})Y_{\hbar}(u^{(i)},x_{1})\\
&&\hspace{4cm}=\delta_{ij}x_{2}^{-1}\delta\left(\frac{x_{1}}{x_{2}}\right).
\end{eqnarray*}
 \et

\begin{proof} Note that from Lemma \ref{lpseudo-auto-pseudo}, pseudo-automorphisms
$\Phi_{i}(x)$ ($1\le i\le l$) are mutually commuting. For $1\le i\le
l$, set
\begin{eqnarray*}
a^{(i)}(x)=Y(u^{(i)},x)\Phi_{i}(x),\ \ \ \
b^{(i)}(x)=Y(v^{(i)},x)\Phi_{i}^{-1}(x).
\end{eqnarray*}
We have
\begin{eqnarray*}
&&p_{ij}(x_{1}-x_{2},\hbar)^{-1}a^{(i)}(x_{1})a^{(j)}(x_{2})\\
&=&p_{ij}(x_{1}-x_{2},\hbar)^{-1}Y(u^{(i)},x_{1})\Phi_{i}(x_{1})Y(u^{(j)},x_{2})\Phi_{j}(x_{2})\\
&=&Y(u^{(i)},x_{1})Y(u^{(j)},x_{2})\Phi_{i}(x_{1})\Phi_{j}(x_{2})\\
&=&q_{ij}Y(u^{(j)},x_{2})Y(u^{(i)},x_{1})\Phi_{j}(x_{2})\Phi_{i}(x_{1})\\
&=&q_{ij}p_{ji}(x_{2}-x_{1},\hbar)^{-1}a^{(j)}(x_{2})a^{(i)}(x_{1}),
\end{eqnarray*}
\begin{eqnarray*}
&&p_{ij}(x_{1}-x_{2},\hbar)^{-1}b^{(i)}(x_{1})b^{(j)}(x_{2})\\
&=&p_{ij}(x_{1}-x_{2},\hbar)^{-1}Y(v^{(i)},x_{1})\Phi_{i}^{-1}(x_{1})Y(v^{(j)},x_{2})\Phi_{j}^{-1}(x_{2})\\
&=&Y(v^{(i)},x_{1})Y(v^{(j)},x_{2})
\Phi_{i}^{-1}(x_{1})\Phi_{j}^{-1}(x_{2})\\
&=&q_{ij}Y(v^{(j)},x_{2})Y(v^{(i)},x_{1})
\Phi_{j}^{-1}(x_{2})\Phi_{i}^{-1}(x_{1})\\
&=&q_{ij}p_{ji}(x_{2}-x_{1},\hbar)^{-1}b^{(j)}(x_{2})b^{(i)}(x_{1}),
\end{eqnarray*}
\begin{eqnarray*}
& &p_{ij}(x_{1}-x_{2},\hbar)a^{(i)}(x_{1})b^{(j)}(x_{2})
-q_{ji}p_{ji}(x_{2}-x_{1},\hbar)b^{(j)}(x_{2})a^{(i)}(x_{1})\\
&=&Y(u^{(i)},x_{1})Y(v^{(j)},x_{2})\Phi_{i}(x_{1})\Phi_{j}^{-1}(x_{2})
-q_{ji}Y(v^{(j)},x_{2})Y(u^{(i)},x_{1})
\Phi_{j}^{-1}(x_{2})\Phi_{i}(x_{1})\\
&=&\delta_{ij}x_{2}^{-1}\delta\left(\frac{x_{1}}{x_{2}}\right)
\Phi_{j}^{-1}(x_{2})\Phi_{i}(x_{1})\\
&=&\delta_{ij}x_{2}^{-1}\delta\left(\frac{x_{1}}{x_{2}}\right)
\Phi_{j}^{-1}(x_{2})\Phi_{i}(x_{2})\\
&=&\delta_{ij}x_{2}^{-1}\delta\left(\frac{x_{1}}{x_{2}}\right).
\end{eqnarray*}
Set $T=\{a^{(i)}(x),\;b^{(i)}(x)\;|\; 1\le i\le l\}\subset
\E_{\hbar}(V_{\bf Q}[[\hbar]])$. {}From Lemma
\ref{lpre-hadic-slocal}, $T$ is $\hbar$-adically $\S$-local and
hence $\hbar$-adically compatible by Lemma \ref{lclosedness}. By
Theorem \ref{tmain-canonical-U}, $T$ generates an $\hbar$-adic
nonlocal vertex algebra $\<T\>$ inside $\E_{\hbar}(V_{\bf
Q}[[\hbar]])$ with $V_{\bf Q}[[\hbar]]$ as a module.

We are going to apply Theorem \ref{tqva-construction} with $V=V_{\bf
Q}[[\hbar]]$, $U=\{ u^{(i)},v^{(i)}\;|\; 1\le i\le l\}$, and
$$Y_{0}(u^{(i)},x)=a^{(i)}(x),\ \ \ Y_{0}(v^{(i)},x)=b^{(i)}(x).$$
 We claim that $V_{{\bf Q}}[[\hbar]]$ is generated from
${\bf 1}$ by field operators $a^{(i)}(x),b^{(i)}(x)$ ($1\le i\le
l$). Let $W$ be the $\C[[\hbar]]$-submodule of $V_{{\bf
Q}}[[\hbar]]$ generated from ${\bf 1}$ by field operators
$a^{(i)}(x),b^{(i)}(x)$ ($1\le i\le l$). We have $\Phi_{i}(x){\bf
1}={\bf 1}\otimes 1$,
\begin{eqnarray*}
&&\Phi_{i}(x)a^{(j)}(x_{1})=p_{ij}(x-x_{1},\hbar)a^{(j)}(x_{1})\Phi_{i}(x),
\\
&&\Phi_{i}(x)b^{(j)}(x_{1})=p_{ij}(x-x_{1},\hbar)^{-1}b^{(j)}(x_{1})\Phi_{i}(x)
\end{eqnarray*}
for $1\le i,j\le l$. It follows from induction that
$\Phi_{i}(x)W\subset W[[x,x^{-1}]]$ for $1\le i\le l$. Similarly, we
have $\Phi_{i}(x)^{-1}W\subset W[[x,x^{-1}]]$. As
$$Y(u^{(i)},x)=a^{(i)}(x)\Phi_{i}(x)^{-1},\ \ \
Y(v^{(i)},x)=b^{(i)}(x)\Phi_{i}(x),$$ it follows that $W$ is closed
under vertex operators $Y(u^{(i)},x)$ and $Y(v^{(i)},x)$ for $1\le
i\le l$. Consequently, $W=V_{{\bf Q}}[[\hbar]]$, as claimed.

By Proposition \ref{pvacuum-generating}, there exists a $K$-module
homomorphism $\pi$ from $K$ to $V_{\bf Q}[[\hbar]]$, sending
$1_{V_{\bf Q}[[\hbar]]}$ to ${\bf 1}$. We are going to prove that
$\pi$ is in fact an isomorphism. First, we see that $\pi$ is
surjective and it gives rise to a surjective linear map $\pi_{0}:
K/\hbar K\rightarrow V_{\bf Q}\; (=V_{\bf Q}[[\hbar]]/\hbar V_{\bf
Q}[[\hbar]])$. By Lemma \ref{lzf-need}, which follows next, we have
\begin{eqnarray*}
&&p_{ij}(x_{1}-x_{2},\hbar)Y_{\E}(a^{(i)}(x),x_{1})Y_{\E}(a^{(j)}(x),x_{2})\\
&&\hspace{1cm}=q_{ij}p_{ji}(x_{2}-x_{1},\hbar)Y_{\E}(a^{(j)}(x),x_{2})Y_{\E}(a^{(i)}(x),x_{1}),\\
&&p_{ij}(x_{1}-x_{2},\hbar)Y_{\E}(b^{(i)}(x),x_{1})Y_{\E}(b^{(j)}(x),x_{2})\\
&&\hspace{1cm}=q_{ij}p_{ji}(x_{2}-x_{1},\hbar)Y_{\E}(b^{(j)}(x),x_{2})Y_{\E}(b^{(i)}(x),x_{1}),\\
&&p_{ij}(x_{1}-x_{2},\hbar)Y_{\E}(a^{(i)}(x),x_{1})Y_{\E}(b^{(j)}(x),x_{2})\\
&&\hspace{2cm}-q_{ji}p_{ji}(x_{2}-x_{1},\hbar)Y_{\E}(b^{(j)}(x),x_{2})Y_{\E}(a^{(i)}(x),x_{1})\\
&&\hspace{1cm}=\delta_{ij}x_{2}^{-1}\delta\left(\frac{x_{1}}{x_{2}}\right)
\end{eqnarray*}
for $1\le i,j\le l$. {}From this, it follows that $K/\hbar K$ is a
vacuum $\A_{\bf Q}$-module with $X_{i}(z)$ and $Y_{i}(z)$ ($1\le
i\le l$) acting as $Y_{\E}(a^{(i)}(x),z)$ and
$Y_{\E}(b^{(i)}(x),z)$, respectively. Furthermore, we see that
$\pi_{0}$ is a surjective $\A_{\bf Q}$-module homomorphism from
$K/\hbar K$ to $V_{\bf Q}$. As every nonzero vacuum $\A_{\bf
Q}$-module is irreducible from \cite{kl}, $\pi_{0}$ must be an
isomorphism.  With $K$ separated and with $V_{\bf Q}[[\hbar]]$
torsion-free, it follows from a result of Enriquez (\cite{en}, Lemma
3.14) that $\pi$ is an isomorphism. Now, by Theorem
\ref{tqva-construction}, there exists an $\hbar$-adic weak quantum
vertex algebra structure on $V_{\bf Q}[[\hbar]]$ with ${\bf 1}$ as
the vacuum vector and with
\begin{eqnarray*}
Y_{\hbar}(u^{(i)},x)=a^{(i)}(x),\ \ \ \
Y_{\hbar}(v^{(i)},x)=b^{(i)}(x)\ \ \ \mbox{ for }1\le i\le l.
\end{eqnarray*}
Then the last assertion follows immediately. As $V_{\bf Q}$ is
nondegenerate, $V_{\bf Q}[[\hbar]]$ is an $\hbar$-adic quantum
vertex algebra. Now, the proof is complete.
\end{proof}

The following is the result we need in the proof of Theorem
\ref{tpre-diag}:

\bl{lzf-need} Let $W$ be a topologically free $\C[[\hbar]]$-module
and let $V$ be an $\hbar$-adic nonlocal vertex subalgebra of
$\E_{\hbar}(W)$. Assume that
$$a(x),b(x)\in V,\; p(x,\hbar),q(x,\hbar)\in \C[x,\hbar],\; f(x,\hbar)\in \C((x))[[\hbar]]$$
with $p(x,0),q(x,0)\ne 0$ such that
\begin{eqnarray}\label{epqabrelation}
\left(\frac{p(x_{1}-x_{2},\hbar)}{q(x_{1}-x_{2},\hbar)}\right)
a(x_{1})b(x_{2})-f(x_{2}-x_{1},\hbar)b(x_{2})a(x_{1})
=\lambda x_{2}^{-1}\delta\left(\frac{x_{1}}{x_{2}}\right),
\end{eqnarray}
where $\lambda$ is a complex number. Then
\begin{eqnarray}\label{epqrelation}
&&\left(\frac{p(x_{1}-x_{2},\hbar)}{q(x_{1}-x_{2},\hbar)}\right)Y_{\E}(a(x),x_{1})Y_{\E}(b(x),x_{2})
-f(x_{2}-x_{1},\hbar)Y_{\E}(b(x),x_{2})Y_{\E}(a(x),x_{1})\nonumber\\
&&\hspace{5cm}=\lambda
x_{2}^{-1}\delta\left(\frac{x_{1}}{x_{2}}\right).
\end{eqnarray}
\el

\begin{proof} From (\ref{epqabrelation}), we have
\begin{eqnarray*}
&&(x_{1}-x_{2})p(x_{1}-x_{2},\hbar)a(x_{1})b(x_{2})
=(x_{1}-x_{2})q(x_{1}-x_{2},\hbar)f(x_{2}-x_{1},\hbar)b(x_{2})a(x_{1})\\
&&\hspace{2.5cm}=(x_{1}-x_{2})p(x_{1}-x_{2},\hbar)
\left(\frac{q(x_{1}-x_{2},\hbar)f(x_{2}-x_{1},\hbar)}{p(-x_{2}+x_{1},\hbar)}\right)b(x_{2})a(x_{1}).
\end{eqnarray*}
 In view of Lemma \ref{lpre-hadic-slocal}, we have
$$a(x_{1})b(x_{2}) \sim
\left(\frac{q(x_{1}-x_{2},\hbar)f(x_{2}-x_{1},\hbar)}{p(-x_{2}+x_{1},\hbar)}\right)b(x_{2})a(x_{1}).$$
By Lemma \ref{literate-formula}, we have
\begin{eqnarray*}
&&Y_{\E}(a(x),x_{0})b(x)
=\Res_{x_{1}}x_{0}^{-1}\delta\left(\frac{x_{1}-x}{x_{0}}\right)
a(x_{1})b(x)\\
&&\hspace{2cm}
-\Res_{x_{1}}x_{0}^{-1}\delta\left(\frac{x-x_{1}}{-x_{0}}\right)
\left(\frac{q(x_{1}-x,\hbar)f(x-x_{1},\hbar)}{p(-x+x_{1},\hbar)}\right)b(x)a(x_{1}).
\end{eqnarray*}
Multiplying both sides by $\frac{p(x_{0},\hbar)}{q(x_{0},\hbar)}$
(viewed as an element of $\C((x_{0}))[[\hbar]]$) we obtain
\begin{eqnarray*}
&&\left(\frac{p(x_{0},\hbar)}{q(x_{0},\hbar)}\right)Y_{\E}(a(x),x_{0})b(x)=\Res_{x_{1}}
x_{0}^{-1}\delta\left(\frac{x_{1}-x}{x_{0}}\right)
\left(\frac{p(x_{1}-x,\hbar)}{q(x_{1}-x,\hbar)}\right)a(x_{1})b(x)\\
&&\hspace{3cm}
-\Res_{x_{1}}x_{0}^{-1}\delta\left(\frac{x-x_{1}}{-x_{0}}\right)
f(x-x_{1},\hbar)b(x)a(x_{1}).
\end{eqnarray*}
Furthermore, for $n\ge 0$, applying $\Res_{x_{0}}x_{0}^{n}$ to both
sides and then using (\ref{epqabrelation}) we get
\begin{eqnarray*}
\Res_{x_{0}}x_{0}^{n}\left(\frac{p(x_{0},\hbar)}{q(x_{0},\hbar)}\right)Y_{\E}(a(x),x_{0})b(x)
=\lambda
\Res_{x_{1}}(x_{1}-x)^{n}x^{-1}\delta\left(\frac{x_{1}}{x}\right)
=\delta_{n,0}\lambda.
\end{eqnarray*}
On the other hand, by Lemma \ref{ph-transp-slocal} we have the
following Jacobi identity in $V$:
\begin{eqnarray*}
&&x_{0}^{-1}\delta\left(\frac{x_{1}-x_{2}}{x_{0}}\right)
Y_{\E}(a(x),x_{1})Y_{\E}(b(x),x_{2})\\
&&\hspace{1cm}-x_{0}^{-1}\delta\left(\frac{x_{2}-x_{1}}{-x_{0}}\right)
\left(\frac{q(x_{1}-x_{2},\hbar)f(x_{2}-x_{1},\hbar)}{p(-x_{2}+x_{1},\hbar)}\right)
Y_{\E}(b(x),x_{2})Y_{\E}(a(x),x_{1})\\
&=&x_{1}^{-1}\delta\left(\frac{x_{2}+x_{0}}{x_{1}}\right)Y_{\E}(Y_{\E}(a(x),x_{0})b(x),x_{2}).
\end{eqnarray*}
Multiplying both sides by $\frac{p(x_{0},\hbar)}{q(x_{0},\hbar)}$
and then taking $\Res_{x_{0}}$ we obtain
\begin{eqnarray*}
&&\left(\frac{p(x_{1}-x_{2},\hbar)}{q(x_{1}-x_{2},\hbar)}\right)Y_{\E}(a(x),x_{1})Y_{\E}(b(x),x_{2})
-f(x_{2}-x_{1},\hbar)Y_{\E}(b(x),x_{2})Y_{\E}(a(x),x_{1})\\
&=&\Res_{x_{0}}x_{1}^{-1}\delta\left(\frac{x_{2}+x_{0}}{x_{1}}\right)
\left(\frac{p(x_{0},\hbar)}{q(x_{0},\hbar)}\right)Y_{\E}(Y_{\E}(a(x),x_{0})b(x),x_{2})\\
&=&\lambda x_{1}^{-1}\delta\left(\frac{x_{2}}{x_{1}}\right),
\end{eqnarray*}
proving (\ref{epqrelation}).
\end{proof}

\bex{quantumbg} {\em Consider the special case with $l=1$,
$q_{11}=1$, and $p_{11}(x,\hbar)=\frac{x+\h}{x}$. In this case (with
${\bf Q}=q_{11}=1$), $\A_{\bf Q}$ is a Weyl algebra and $V_{\bf Q}$
is a vertex algebra. By Theorem \ref{tpre-diag}, there exists an
$\hbar$-adic quantum vertex algebra $V_{\bf Q}[[\hbar]]$ with
generators $u,v$ such that
\begin{eqnarray*}
&&\left(\frac{x_{1}-x_{2}}{x_{1}-x_{2}+\hbar}\right)Y(u,x_{1})Y(u,x_{2})
=\left(\frac{x_{2}-x_{1}}{x_{2}-x_{1}+\hbar}\right)Y(u,x_{2})Y(u,x_{1}),\\
&&\left(\frac{x_{1}-x_{2}}{x_{1}-x_{2}+\hbar}\right)Y(v,x_{1})Y(v,x_{2})
=\left(\frac{x_{2}-x_{1}}{x_{2}-x_{1}+\hbar}\right)Y(v,x_{2})Y(v,x_{1}),\\
&&\left(\frac{x_{1}-x_{2}+\hbar}{x_{1}-x_{2}}\right)Y(u,x_{1})Y(v,x_{2})-
\left(\frac{x_{2}-x_{1}+\hbar}{x_{2}-x_{1}}\right)Y(v,x_{2})Y(u,x_{1})
=x_{2}^{-1}\delta\left(\frac{x_{1}}{x_{2}}\right).
\end{eqnarray*}
This gives an $\hbar$-deformed $\beta\gamma$-system (cf.
\cite{efk}).}\eex

\bex{lattice} {\em Consider the case with $l=1$, $q_{11}=-1$, and
$p_{11}(x,\hbar)=\frac{x+\h}{x}$. In this case, $\A_{\bf Q}$ is a
Clifford algebra and $V_{\bf Q}$ is a vertex superalgebra (cf.
\cite{ffr}). Theorem \ref{tpre-diag} asserts that there exists an
$\hbar$-adic quantum vertex algebra $V_{\bf Q}[[\hbar]]$ with
generators $u,v$ such that
\begin{eqnarray*}
&&\left(\frac{x_{1}-x_{2}}{x_{1}-x_{2}+\hbar}\right)Y(u,x_{1})Y(u,x_{2})
=-\left(\frac{x_{2}-x_{1}}{x_{2}-x_{1}+\hbar}\right)Y(u,x_{2})Y(u,x_{1}),\\
&&\left(\frac{x_{1}-x_{2}}{x_{1}-x_{2}+\hbar}\right)Y(v,x_{1})Y(v,x_{2})
=-\left(\frac{x_{2}-x_{1}}{x_{2}-x_{1}+\hbar}\right)Y(v,x_{2})Y(v,x_{1}),\\
&&\left(\frac{x_{1}-x_{2}+\hbar}{x_{1}-x_{2}}\right)Y(u,x_{1})Y(v,x_{2})+
\left(\frac{x_{2}-x_{1}+\hbar}{x_{2}-x_{1}}\right)Y(v,x_{2})Y(u,x_{1})
=x_{2}^{-1}\delta\left(\frac{x_{1}}{x_{2}}\right).
\end{eqnarray*}
Let $L=\Z\alpha$ be a rank-one lattice with $\<\alpha,\alpha\>=1$.
Associated to $L$, one has a vertex superalgebra $V_{L}$ (cf.
\cite{dl}). It is known that vertex superalgebras $V_{\bf Q}$ and
$V_{L}$ are isomorphic with $u=e^{\alpha}$ and $v=e^{-\alpha}$. In
view of this, $V_{\bf Q}[[\hbar]]$ is an $\hbar$-deformation of
$V_{L}$.}\eex

\section{$\hbar$-adic quantum vertex algebras associated with double Yangians}

In this section, we shall associate the centrally extended double
Yangian of $\sl_{2}$ (see \cite{kh}) with $\hbar$-adic quantum
vertex algebras. This can be viewed as an $\hbar$-adic version of
the corresponding result of \cite{li-infinity} for the centerless
double Yangian of $\sl_{2}$ with $\hbar$ evaluated as a nonzero
complex number.

The following is a variant of the centrally extended double Yangian
$\widehat{DY_{\hbar}(\sl_{2})}$ in \cite{kh}:

\bd{ddy-withc} {\em Define $\widehat{DY_{\hbar}(\sl_{2})}$ to be the
topological associative algebra over $\C[[\h]]$ with generators
$$e_{m},\ \ f_{m},\ \ h_{m},\ \ c,\ \ d\ \ (m\in \Z),$$
which are grouped together in terms of generating functions
\begin{eqnarray*}
& &e(x)=\sum_{m\in \Z}e_{m}x^{-m-1},\ \ \ \ \
f(x)=\sum_{m\in \Z}f_{m}x^{-m-1},\\
& &h^{+}(x)=1+\h\sum_{m\ge 0}h_{m}x^{-m-1}, \ \ \ \ \
h^{-}(x)=1-\h\sum_{m< 0}h_{m}x^{-m-1},
\end{eqnarray*}
with $c$ a central element, subject to relations
\begin{eqnarray}\label{enew-yd}
& &[d,e(x)]=\frac{d}{dx}e(x), \ \ \ \ [d,f(x)]=\frac{d}{dx}f(x),
\ \ \ \ [d,h^{\pm}(x)]=\frac{d}{dx}h^{\pm}(x), \nonumber\\
& &e(x)e(y)=\frac{y-x-\h}{y-x+\h} e(y)e(x),\nonumber\\
& &f(x)f(y)=\frac{y-x+\h}{y-x-\h}f(y)f(x),\nonumber\\
& &h^{+}(x)e(y)=\frac{x-y+\h}{x-y-\h}e(y)h^{+}(x),\nonumber\\
& &h^{+}(x)f(y)=\frac{x-y-(1+c)\h}{x-y+(1-c)\h}f(y)h^{+}(x),\nonumber\\
& &h^{-}(x)e(y)=\frac{y-x-\h}{y-x+\h}e(y)h^{-}(x),\nonumber\\
& &h^{-}(x)f(y)=\frac{y-x+\h}{y-x-\h}f(y)h^{-}(x),\nonumber\\
& &h^{\pm}(x)h^{\pm}(y)=h^{\pm}(y)h^{\pm}(x),\nonumber\\
& &h^{+}(x)h^{-}(y)=\frac{x-y+\h}{x-y-\h}\cdot
\frac{x-y-(1+c)\h}{x-y+(1-c)\h}
\cdot h^{-}(y)h^{+}(x),\nonumber\\
& &[e(x),f(y)] =\frac{1}{\h}\left( x^{-1}\delta\left(\frac{y+\h
c}{x}\right)h^{+}(x)
-x^{-1}\delta\left(\frac{y}{x}\right)h^{-}(y)\right),
\end{eqnarray}
where by convention for $a\in \C[c]$,
\begin{eqnarray*}
\frac{1}{x-y-a\h}
&=&\sum_{n\ge 0}a^{n}(x-y)^{-n-1}\hbar^{n}\in \C[c,(x-y)^{-1}][[\hbar]],\\
 \frac{1}{y-x-a \h}&=&\sum_{n\ge 0}a^{n}(y-x)^{-n-1}\hbar^{n}
 \in \C[c,(y-x)^{-1}][[\hbar]].
\end{eqnarray*}}
\ed

\br{rdetails} {\em  Here we give some details for the definition.
Let $T$ be the free associative algebra over $\C$ with generators
$e_{m},f_{m},h_{m}$ $(m\in \Z)$, $c$, $d$. Set
$$\deg c=\deg d=0,\ \deg e_{m}=\deg f_{m}=\deg h_{m}=m
\ \ \mbox{for }m\in \Z,$$ to make $T$ a $\Z$-graded algebra
$T=\coprod_{n\in \Z}T_{n}$. For $n\in \Z$, set
$$F_{n}(T)=\sum_{m\ge n} T_{m}\subset T.$$
This defines a decreasing filtration $\{F_{n}(T)\}_{n\in \Z}$ of $T$
with $\cap_{n\in \Z}F_{n}(T)=0$. Denote by $\overline{T}$ the formal
completion of $T$ with respect to this filtration. Then
$\overline{T}[[\hbar]]$ is an algebra over $\C[[\hbar]]$. The
algebra $\widehat{DY_{\hbar}(\sl_{2})}$ is simply the quotient
algebra of $\overline{T}[[\hbar]]$ modulo all the relations
corresponding to (\ref{enew-yd}). The standard double Yangian
$\widehat{DY_{\hbar}(\sl_{2})}$ (see \cite{kh}) was related to the
increasing filtration $\{C_{k}\}_{k\in \Z}$, where $C_{k}=\sum_{m\le
k}T_{m}$ for $k\in \Z$. } \er

\br{rsimple-fact10} {\em Let $W$ be a $\C[[\hbar]]$-module and let
$a(x),b(x)\in \E_{\hbar}(W)$. Assume that $a(x)\in (\End
W)[[x^{-1}]]$. Note that the equalities
\begin{eqnarray*}
a(x_{1})b(x_{2})&=&\frac{x_{1}-x_{2}-\hbar}{x_{1}-x_{2}+\hbar}b(x_{2})a(x_{1}),\\
a(x_{1})b(x_{2})&=&\frac{x_{2}-x_{1}+\hbar}{x_{2}-x_{1}-\hbar}b(x_{2})a(x_{1})
\;\left(=\frac{-x_{2}+x_{1}-\hbar}{-x_{2}+x_{1}+\hbar}b(x_{2})a(x_{1})\right)
\end{eqnarray*}
both make sense, but they are not equivalent in general. On the
other hand, the following equalities both make sense and are
equivalent
\begin{eqnarray*}
a(x_{1})b(x_{2})&=&\frac{x_{1}-x_{2}-\hbar}{x_{1}-x_{2}+\hbar}b(x_{2})a(x_{1}),\\
b(x_{2})a(x_{1})&=&\frac{x_{1}-x_{2}+\hbar}{x_{1}-x_{2}-\hbar}a(x_{1})b(x_{2}).
\end{eqnarray*}
} \er

By {\em a $\widehat{DY_{\hbar}(\sl_{2})}$-module} we mean a
$\C[[\hbar]]$-module $W$ on which $e(x),f(x),h^{\pm}(x)$ and $c,d$
act such that
\begin{eqnarray*}
e(x), f(x), h^{\pm }(x)\in \E_{\hbar}(W)
\end{eqnarray*}
and such that all the defining relations in (\ref{enew-yd}) hold.
 A $\widehat{DY_{\hbar}(\sl_{2})}$-module $W$ is said to be of {\em
level} $\ell\in \C$ if $c$ acts on $W$ as scalar $\ell$. Let $W$ be
a $\widehat{DY_{\hbar}(\sl_{2})}$-module of level $\ell$.
 Set
\begin{eqnarray}
U_{W}=\{1_{W},e(x),f(x), h^{+}(x),h^{-}(x)\} \subset \E_{\hbar}(W).
\end{eqnarray}
Note that from the last defining relation we have
$$(x_{1}-x_{2})(x_{1}-x_{2}-\ell\hbar)[e(x_{1}),f(x_{2})]=0.$$
In view of Lemma \ref{lpre-hadic-slocal}, $U_{W}$ is
$\hbar$-adically $\S$-local. Then by Theorem \ref{tmain-canonical},
$U_{W}$ generates an $\hbar$-adic nonlocal vertex algebra $V_{W}$
inside $\E_{\hbar}(W)$. We are going to show that the space $V_{W}$
is naturally a module for a variant of
$\widehat{DY_{\hbar}(\sl_{2})}$.

\bd{dcover-algebra} {\em Let $\widetilde{DY_{\hbar}(\sl_{2})}$ be
the topological associative algebra over $\C[[\hbar]]$ with
generators
\begin{eqnarray*}
\tilde{e}_{m},\ \ \tilde{f}_{m}, \ \ \tilde{h}^{\pm}_{m}\ \ (m\in
\Z),\ \ c, \ \ d,
\end{eqnarray*}
and with generating functions
\begin{eqnarray}
& &\tilde{e}(x)=\sum_{m\in \Z}\tilde{e}_{m}x^{-m-1}, \ \ \ \ \
\tilde{f}(x)=\sum_{m\in \Z}\tilde{f}_{m}x^{-m-1},\nonumber\\
& &\tilde{h}^{+}(x)=1+ \hbar \sum_{m\in
\Z}\tilde{h}^{+}_{m}x^{-m-1}, \ \ \ \ \tilde{h}^{-}(x)=1-\hbar
\sum_{m\in \Z}\tilde{h}^{-}_{m}x^{-m-1},
\end{eqnarray}
subject to relations
\begin{eqnarray}\label{enew-algebra}
& &[d,\tilde{e}(x)]=\frac{d}{dx}\tilde{e}(x), \ \ \ \
[d,\tilde{f}(x)]=\frac{d}{dx}\tilde{f}(x), \ \ \ \
[d,\tilde{h}^{\pm}(x)]
=\frac{d}{dx}\tilde{h}^{\pm}(x), \nonumber\\
& &\tilde{e}(x)\tilde{e}(y)
=\frac{y-x-\h}{y-x+\h} \tilde{e}(y)\tilde{e}(x),\nonumber\\
& &\tilde{f}(x)\tilde{f}(y)
=\frac{y-x+\h}{y-x-\h}\tilde{f}(y)\tilde{f}(x),\nonumber\\
& &\tilde{e}(y)\tilde{h}^{+}(x)
=\frac{x-y-\h}{x-y+\h}\tilde{h}^{+}(x)\tilde{e}(y),\nonumber\\
& &\tilde{f}(y)\tilde{h}^{+}(x)
=\frac{x-y+(1-c)\h}{x-y-(1+c)\h}\tilde{h}^{+}(x)\tilde{f}(y),\nonumber\\
& &\tilde{h}^{-}(x)\tilde{e}(y)
=\frac{y-x-\h}{y-x+\h}\tilde{e}(y)\tilde{h}^{-}(x),\nonumber\\
& &\tilde{h}^{-}(x)\tilde{f}(y)
=\frac{y-x+\h}{y-x-\h}\tilde{f}(y)\tilde{h}^{-}(x),\nonumber\\
& &[\tilde{h}^{\pm}(x),\tilde{h}^{\pm}(y)]=0,\nonumber\\
& &\tilde{h}^{-}(y)\tilde{h}^{+}(x)=\frac{x-y-\h}{x-y+\h}\cdot
\frac{x-y+(1-c)\h}{x-y-(1+c)\h}
\cdot \tilde{h}^{+}(x)\tilde{h}^{-}(y),\nonumber\\
& &[\tilde{e}(x),\tilde{f}(y)] =\frac{1}{\h}\left(
x^{-1}\delta\left(\frac{y+\h c}{x}\right)\tilde{h}^{+}(x)
-x^{-1}\delta\left(\frac{y}{x}\right)\tilde{h}^{-}(y)\right).
\end{eqnarray}} \ed

Similarly, we define a $\widetilde{DY_{\hbar}(\sl_{2})}$-module to
be a $\C[[\hbar]]$-module on which $\tilde{e}(x),\tilde{f}(x)$,
$\tilde{h}^{\pm}(x)$ and $c,d$ act with $\tilde{e}(x),\tilde{f}(x)$,
$\tilde{h}^{\pm}(x)\in \E_{\hbar}(W)$, satisfying all the defining
relations. A vector $w$ in a
$\widetilde{DY_{\hbar}(\sl_{2})}$-module is called a {\em vacuum
vector} if
$$dw=0,\ \ \tilde{e}_{m}w=\tilde{f}_{m}w=\tilde{h}^{\pm}_{m}w=0
\ \ \mbox{ for }m\ge 0.$$ A {\em vacuum
$\widetilde{DY_{\hbar}(\sl_{2})}$-module} is a module equipped with
a vacuum vector which generates the whole module.

\br{rquotient-algebra} {\em In view of Remark \ref{rsimple-fact10},
the following relations hold in $\widehat{DY_{\hbar}(\sl_{2})}$:
\begin{eqnarray*}
& &e(y)h^{+}(x)=\frac{x-y-\h}{x-y+\h}h^{+}(x)e(y),\\
& &f(y)h^{+}(x)
=\frac{x-y+(1-c)\h}{x-y-(1+c)\h}h^{+}(x)f(y),\\
& &h^{-}(y)h^{+}(x) =\frac{x-y-\h}{x-y+\h}\cdot
\frac{x-y+(1-c)\h}{x-y-(1+c)\h} \cdot h^{+}(x)h^{-}(y).
\end{eqnarray*}
On the other hand, these three relations also imply the
corresponding original relations. Thus
$\widehat{DY_{\hbar}(\sl_{2})}$ is isomorphic to the quotient
algebra of $\widetilde{DY_{\hbar}(\sl_{2})}$ modulo the relations
\begin{eqnarray}
\tilde{h}^{+}_{m}=0\;\;\;\mbox{ for }m<0 \ \mbox{ and }\
\tilde{h}^{-}_{n}=0\;\;\;\mbox{ for }n\ge 0.
\end{eqnarray}
Consequently, every $\widehat{DY_{\hbar}(\sl_{2})}$-module of level
$\ell$ is naturally a $\widetilde{DY_{\hbar}(\sl_{2})}$-module of
level $\ell$.} \er

Set
\begin{eqnarray}
\tilde{h}^{+}(x)'=\sum_{k\ge
0}\tilde{h}_{k}x^{-k-1}=\frac{1}{\hbar}(\tilde{h}^{+}(x)-1),\ \ \ \
\tilde{h}^{-}(x)'=\sum_{k<0}\tilde{h}_{k}x^{-k-1}=\frac{1}{\hbar}(1-\tilde{h}^{-}(x)).
\end{eqnarray}
The following is straightforward:

\bl{lnew-hoperators} In terms of $\tilde{h}^{\pm}(x)'$, those
relations involving $\tilde{h}^{\pm}(x)$ in (\ref{enew-algebra})
become
\begin{eqnarray*}
& &\tilde{h}^{-}(x)'\tilde{e}(y)
=\frac{y-x-\hbar}{y-x+\hbar}\tilde{e}(y)\tilde{h}^{-}(x)'
+\frac{2}{y-x+\hbar}\tilde{e}(y),\\
& &\tilde{h}^{-}(x)'\tilde{f}(y)
=\frac{y-x+\hbar}{y-x-\hbar}\tilde{f}(y)\tilde{h}^{-}(x)'
+\frac{-2}{y-x-\hbar}\tilde{f}(y), \\
&&\tilde{e}(y)\tilde{h}^{+}(x)' =\frac{x-y-\hbar}{x-y+\hbar}\cdot
\tilde{h}^{+}(x)'\tilde{e}(y)+
\frac{-2}{x-y+\hbar}\tilde{e}(y),\\
& &\tilde{f}(y)\tilde{h}^{+}(x)'
=\frac{x-y+(1-c)\hbar}{x-y-(1+c)\hbar}\cdot
\tilde{h}^{+}(x)'\tilde{f}(y)+
\frac{2}{x-y-(1+c)\hbar}\tilde{f}(y),\\
& &[\tilde{h}^{\pm}(x)',\tilde{h}^{\pm}(y)']=0,\\
&&\tilde{h}^{-}(y)'\tilde{h}^{+}(x)'
=F(\hbar)\tilde{h}^{+}(x)'\tilde{h}^{-}(y)'\\
&&\hspace{2cm}+
\frac{-2c}{(x-y+\h)(x-y-(1+c)\h)}\left(1+\hbar\tilde{h}^{+}(x)'-\hbar\tilde{h}^{-}(y)'\right),
\end{eqnarray*}
where
$$F(\hbar)=\frac{x-y-\h}{x-y+\h}\cdot \frac{x-y+(1-c)\h}{x-y-(1+c)\h},$$
 and
 \begin{eqnarray*}
&&[\tilde{e}(x),\tilde{f}(y)]\\
&=& x^{-1}\delta\left(\frac{y}{x}\right)
\left(\tilde{h}^{+}(x)'+\tilde{h}^{-}(x)'\right) +\sum_{k\ge
1}\frac{c^{k} \hbar^{k-1}}{k!}(1+\hbar \tilde{h}^{+}(x)')
\left(\frac{\partial}{\partial
y}\right)^{k}x^{-1}\delta\left(\frac{y}{x}\right).
\end{eqnarray*}
\el

\bp{ptransfer} Let $W$ be a topologically free
$\widetilde{DY_{\hbar}(\sl_{2})}$-module of level $\ell$. Set
$$\widetilde{U}_{W}
=\{1_{W}, \tilde{e}(x), \tilde{f}(x),\tilde{h}^{+}(x)',
\tilde{h}^{-}(x)'\}\subset \E_{\hbar}(W).$$ Then $\widetilde{U}_{W}$
is $\hbar$-adically $\S$-local and the $\hbar$-adic nonlocal vertex
algebra  $\widetilde{V}_{W}$ generated by $\widetilde{U}_{W}$ is a
$\widetilde{DY_{\hbar}(\sl_{2})}$-module of level $\ell$ with
$\tilde{e}(x_{0}),\; \tilde{f}(x_{0})$ and $\tilde{h}^{\pm}(x_{0})$
acting as $Y_{\E}(\tilde{e}(x),x_{0}),\; Y_{\E}(\tilde{f}(x),x_{0})$
and $Y_{\E}(\tilde{h}^{\pm}(x),x_{0})$, respectively, and with $d$
acting as $\D=d/dx$ and $c$ acting as scalar $\ell$. Furthermore,
$1_{W}$ generates a vacuum $\widetilde{DY_{\hbar}(\sl_{2})}$-module
of level $\ell$. \ep

\begin{proof}  With
the commutation relations (\ref{enew-algebra}) and with Lemma
\ref{lnew-hoperators}, by Lemma \ref{lpre-hadic-slocal},
$\widetilde{U}_{W}$ is $\hbar$-adically $\S$-local. Note that the
$\D$-operator of the $\hbar$-adic nonlocal vertex algebra
$\widetilde{V}_{W}$ is exactly the formal differential operator
$d/dx$ and we have
\begin{eqnarray*}
& &[\D,Y_{\E}(\tilde{e}(x),x_{0})]
=\frac{d}{dx_{0}}Y_{\E}(\tilde{e}(x),x_{0}), \ \ \ \
[\D,Y_{\E}(\tilde{f}(x),x_{0})]
=\frac{d}{dx_{0}}Y_{\E}(\tilde{f}(x),x_{0}), \\
& &\ \ \ \ [\D,Y_{\E}(\tilde{h}^{\pm}(x),x_{0})]
=\frac{d}{dx_{0}}Y_{\E}(\tilde{h}^{\pm}(x),x_{0}).
 \end{eqnarray*}
 With
the commutation relations (\ref{enew-algebra}), by Proposition
\ref{phadic-slocal} we have
\begin{eqnarray*}
& &Y_{\E}(\tilde{e}(x),x_{1})Y_{\E}(\tilde{e}(x),x_{2})
=\frac{x_{2}-x_{1}-\h}{x_{2}-x_{1}+\h}\cdot
Y_{\E}(\tilde{e}(x),x_{2})Y_{\E}(\tilde{e}(x),x_{1}),\nonumber\\
& &Y_{\E}(\tilde{f}(x),x_{1})Y_{\E}(\tilde{f}(x),x_{2})
=\frac{x_{2}-x_{1}+\h}{x_{2}-x_{1}-\h}\cdot
Y_{\E}(\tilde{f}(x),x_{2})Y_{\E}(\tilde{f}(x),x_{1}),\nonumber\\
& &Y_{\E}(\tilde{e}(x),x_{2})Y_{\E}(\tilde{h}^{+}(x),x_{1})
=\frac{x_{1}-x_{2}-\h}{x_{1}-x_{2}+\h}\cdot
Y_{\E}(\tilde{h}^{+}(x),x_{1})Y_{\E}(\tilde{e}(x),x_{2}),\nonumber\\
& &Y_{\E}(\tilde{f}(x),x_{2})Y_{\E}(\tilde{h}^{+}(x),x_{1})
=\frac{x_{1}-x_{2}+(1-\ell)\h}{x_{1}-x_{2}-(1+\ell)\h} \cdot
Y_{\E}(\tilde{h}^{+}(x),x_{1})Y_{\E}(\tilde{f}(x),x_{2}),\nonumber\\
& &Y_{\E}(\tilde{h}^{-}(x),x_{1})Y_{\E}(\tilde{e}(x),x_{2})
=\frac{x_{2}-x_{1}-\h}{x_{2}-x_{1}+\h}\cdot
Y_{\E}(\tilde{e}(x),x_{2})Y_{\E}(\tilde{h}^{-}(x),x_{1}),\nonumber\\
& &Y_{\E}(\tilde{h}^{-}(x),x_{1})Y_{\E}(\tilde{f}(x),x_{2})
=\frac{x_{2}-x_{1}+\h}{x_{2}-x_{1}-\h}\cdot
Y_{\E}(\tilde{f}(x),x_{2})Y_{\E}(\tilde{h}^{-}(x),x_{1}),\nonumber\\
&
&[Y_{\E}(\tilde{h}^{\pm}(x)',x_{1}),Y_{\E}(\tilde{h}^{\pm}(x)',x_{2})]=0,\\
& &Y_{\E}(\tilde{h}^{-}(x),x_{2})Y_{\E}(\tilde{h}^{+}(x),x_{1})\nonumber\\
& &\ \ \ \ =\frac{x_{1}-x_{2}-\h}{x_{1}-x_{2}+\h}\cdot
\frac{x_{1}-x_{2}+(1-\ell)\h}{x_{1}-x_{2}-(1+\ell)\h} \cdot
Y_{\E}(\tilde{h}^{+}(x),x_{1})Y_{\E}(\tilde{h}^{-}(x),x_{2}).
\end{eqnarray*}
Recall that $W$ is a faithful $\widetilde{V}_{W}$-module. Using
Proposition \ref{phmodule-algebra-relations-v2} we get
\begin{eqnarray*}
& &[Y_{\E}(\tilde{e}(x),x_{1}),Y_{\E}(\tilde{f}(x),x_{2})] =
x_{1}^{-1}\delta\left(\frac{x_{2}}{x_{1}}\right)
Y_{\E}\left(\tilde{h}^{+}(x)'+\tilde{h}^{-}(x)',x_{1}\right)\\
& &\hspace{3cm} +\sum_{k\ge
1}\frac{\ell^{k}\hbar^{k-1}}{k!}Y_{\E}(1+\hbar
\tilde{h}^{+}(x)',x_{1})\left(\frac{\partial}{\partial
x_{2}}\right)^{k}
x_{1}^{-1}\delta\left(\frac{x_{2}}{x_{1}}\right)\\
& &\ \ \ \  =\frac{1}{\h}\left(
x_{1}^{-1}\delta\left(\frac{x_{2}+\ell \h
}{x_{1}}\right)Y_{\E}(\tilde{h}^{+}(x),x_{1})
-x_{1}^{-1}\delta\left(\frac{x_{2}}{x_{1}}\right)Y_{\E}(\tilde{h}^{-}(x),x_{2})\right).
\end{eqnarray*}
This shows that $V_{W}$ is a
$\widetilde{DY_{\hbar}(\sl_{2})}$-module of level $\ell$. Clearly,
${\bf 1}$ is a vacuum  vector of $\widetilde{V}_{W}$ viewed as a
$\widetilde{DY_{\hbar}(\sl_{2})}$-module, so it generates a vacuum
$\widetilde{DY_{\hbar}(\sl_{2})}$-module.
\end{proof}

\bd{dlimit-lie-algebra} {\em Define $K$ to be a Lie algebra over
$\C$ with a basis $\{c, E_{m},F_{m},I_{m},J_{m}\;|\; m\in \Z\}$ and
with the Lie bracket relations
\begin{eqnarray*}
& &[c, K]=0,\\
& &[E(x),E(y)]=0, \ \ \ \ [F(x),F(y)]=0,\\
& &[I(x),I(y)]=0,\ \ \ \ [J(x),J(y)]=0,\ \ \ \
[I(x),J(y)]=\frac{2c}{(x-y)^{2}},\\
 & &[I(x),E(y)]=\frac{2}{x-y}E(y),\ \ \ \
[I(x),F(y)]=\frac{-2}{x-y}F(y),\\
& &[J(x),E(y)]=\frac{2}{y-x}E(y),\ \ \ \
[J(x),F(y)]=\frac{-2}{y-x}F(y),\\
& &[E(x),F(y)]=x^{-1}\delta\left(\frac{y}{x}\right)(I(y)+J(y))
+c\frac{\partial}{\partial y}x^{-1}\delta\left(\frac{y}{x}\right),
\end{eqnarray*}
where $a(x)=\sum_{m\in \Z}a_{m}x^{-m-1}$ for $a=E,F,I,J$.} \ed

\br{rlimit-lie-algebra} {\em  Consider the product Lie algebra
$\sl_{2}\oplus \C z$. Extend the normalized Killing form on
$\sl_{2}$ to a symmetric (invariant) bilinear form on $\sl_{2}\oplus
\C z$ by
$$\< \sl_{2},z\>=0, \ \ \ \< z,z\>=0.$$
Then we have an affine Lie algebra $\widehat{\sl_{2}\oplus \C z}$.
It is readily to see that Lie algebra $K$ is isomorphic to
$\widehat{\sl_{2}\oplus \C z}$ with
\begin{eqnarray*}
E(x)=e(x),\ \ \ F(x)=f(x),\ \ \ I(x)=h(x)^{+}+z(x),\ \ \
J(x)=h(x)^{-}-z(x),
\end{eqnarray*}
and $c={\bf k}$ (the central element of $\widehat{\sl_{2}\oplus \C
z}$), where for $a\in \sl_{2}\oplus\C z$,
$$a(x)=\sum_{m\in \Z}a(m)x^{-m-1}$$
and
$$h(x)^{+}=\sum_{m\ge 0}h(m)x^{-m-1},\ \ \ \ \
h(x)^{-}=\sum_{m< 0}h(m)x^{-m-1}.$$ } \er

A vector $w$ in a $K$-module is called a {\em vacuum vector} if
$$E_{m}w=F_{m}w=I_{m}w=J_{m}w=0\ \ \mbox{ for }m\ge 0.$$
A {\em vacuum $K$-module} is a module equipped with a vacuum vector
which generates the whole module. Denote by $K_{\ge 0}$ the
 linear span of $c, E_{m},F_{m},I_{m},J_{m}$
for $m\ge 0$. It is clear that $K_{\ge 0}$ is a Lie subalgebra. Let
$\ell$ be a complex number. Denote by $\C_{\ell}$ the
$1$-dimensional $K_{\ge 0}$-module with $c$ acting as $\ell$ and
with all the other generators of $K_{\ge 0}$ acting trivially. Form
the induced $K$-module
\begin{eqnarray}
V_{K}(\ell,0)=U(K)\otimes_{U(K_{\ge 0})} \C_{\ell}.
\end{eqnarray}
Set ${\bf 1}=1\otimes 1\in V_{K}(\ell,0)$. Then $V_{K}(\ell,0)$ is a
vacuum $K$-module of level $\ell$, which is universal in the obvious
sense. In view of the connection of $K$ with $\widehat{\sl_{2}\oplus
\C z}$, $V_{K}(\ell,0)$ is also a universal vacuum
$\widehat{\sl_{2}\oplus \C z}$-module of level $\ell$. If $\ell$ is
generic, it is well known that the universal vacuum
$\widehat{\sl_{2}\oplus \C z}$-module of level $\ell$ is
irreducible. It follows that  $V_{K}(\ell,0)$ is an irreducible
$K$-module if $\ell$ is generic.

Set
$$E=E(-1){\bf 1}, \ \ F=F(-1){\bf 1},\ \ I=I(-1){\bf 1},\ \ J=J(-1){\bf
1}\in V_{K}(\ell,0).$$ Clearly, $\{1, E(x),F(x),I(x),J(x)\}$ is
$\S$-local. It follows from \cite{li-qva2} that there exists a
(unique) non-degenerate quantum vertex algebra structure on
$V_{K}(\ell,0)$ over $\C$ with ${\bf 1}$ as the vacuum vector and
with
$$Y(E,x)=E(x),\ \ Y(F,x)=F(x),\ \ Y(I,x)=I(x),\ \ Y(J,x)=J(x).$$

\bp{plimit} Let $W$ be any $\widetilde{DY_{\hbar}(\sl_{2})}$-module
of level $\ell$. Then $W/\hbar W$ is a $K$-module of level $\ell$
with $E_{m},F_{m},I_{m},J_{m}$ for $m\in \Z$ acting as
$\tilde{e}_{m},\tilde{f}_{m},\tilde{h}^{+}_{m},\tilde{h}^{-}_{m}$,
respectively. If $W$ is a vacuum
$\widetilde{DY_{\hbar}(\sl_{2})}$-module of level $\ell$, then
$W/\hbar W$ is a vacuum $K$-module of level $\ell$. \ep

\begin{proof} The first assertion follows immediately from the defining
relations in (\ref{enew-algebra}) and the relations in Lemma
\ref{lnew-hoperators}. If $W$ is a vacuum
$\widetilde{DY_{\hbar}(\sl_{2})}$-module of level $\ell$, we see
that $W/\hbar W$ is a vacuum $K$-module of level $\ell$.
\end{proof}

\bt{tdouble-yangian-qva-module} Let $\ell$ be a complex number.
Assume that there exists a vacuum
$\widetilde{DY_{\hbar}(\sl_{2})}$-module
$V(\widetilde{DY_{\hbar}(\sl_{2})},\ell)$ of level $\ell$ which is
universal in the obvious sense and which is topologically free. Then
there exists a unique $\hbar$-adic weak quantum vertex algebra
structure on $V(\widetilde{DY_{\hbar}(\sl_{2})},\ell)$ with ${\bf
1}$ as the vacuum vector and with
\begin{eqnarray}
Y(\tilde{e}_{-1}{\bf 1},x)=\tilde{e}(x),\ \ Y(\tilde{f}_{-1}{\bf
1},x)=\tilde{f}(x),\ \ Y(\tilde{h}^{\pm}_{-1}{\bf
1},x)=\tilde{h}^{\pm}(x)'.
\end{eqnarray}
If $\ell$ is generic, then $V(\widetilde{DY_{\hbar}(\sl_{2})},\ell)$
is a non-degenerate $\hbar$-adic quantum vertex algebra.
Furthermore, on any $\widetilde{DY_{\hbar}(\sl_{2})}$-module $W$ of
level $\ell$, there exists a unique
$V(\widetilde{DY_{\hbar}(\sl_{2})},\ell)$-module structure such that
$$Y_{W}(\tilde{e}_{-1}{\bf 1},x)=\tilde{e}(x),\ \ \
Y_{W}(\tilde{f}_{-1}{\bf 1},x)=\tilde{f}(x),\ \ \
Y_{W}(\tilde{h}^{\pm}_{-1}{\bf 1},x)=\tilde{h}^{\pm}(x)'.$$  \et

\begin{proof} Let $W$ be a $\widetilde{DY_{\hbar}(\sl_{2})}$-module of
level $\ell$. By Proposition \ref{ptransfer}, the $\hbar$-adic
nonlocal vertex algebra $\widetilde{V}_{W}$ generated by
$\tilde{U}_{W}$ is a $\widetilde{DY_{\hbar}(\sl_{2})}$-module of
level $\ell$ and the submodule generated by $1_{W}$ is a vacuum
$\widetilde{DY_{\hbar}(\sl_{2})}$-module of level $\ell$ with vacuum
vector $1_{W}$. It follows that there exists a
$\widetilde{DY_{\hbar}(\sl_{2})}$-module homomorphism $\psi_{W}$
from $V(\widetilde{DY_{\hbar}(\sl_{2})},\ell)$ to
$\widetilde{V}_{W}$, sending ${\bf 1}$ to $1_{W}$. Specializing
$W=V(\widetilde{DY_{\hbar}(\sl_{2})},\ell)$ and then applying
Theorem \ref{tqva-construction} with
$V=V(\widetilde{DY_{\hbar}(\sl_{2})},\ell)$, $U=\{ {\bf
1},\tilde{e}_{-1}{\bf 1},\tilde{f}_{-1}{\bf
1},\tilde{h}^{\pm}_{-1}{\bf 1}\}$, we obtain the first assertion.
For a general $W$, from Theorem \ref{tmain-canonical-U}, $W$ is a
canonical module for $\widetilde{V}_{W}$. The $\C[[\hbar]]$-module
map $\psi_{W}$ satisfies
$$\psi_{W}({\bf 1})=1_{W},\ \ \psi_{W}(u_{m}v)=u(x)_{m}\psi_{W}(v)
\ \ \ \mbox{ for } u\in U,\; v\in V,\; m\in \Z.$$ As ${\bf 1}$
generates $V(\widetilde{DY_{\hbar}(\sl_{2})},\ell)$ as a
$\widetilde{DY_{\hbar}(\sl_{2})}$-module, it follows that $\psi_{W}$
is a homomorphism of $\hbar$-adic nonlocal vertex algebras. Then the
last assertion follows.

As for the second assertion, since $\ell$ is generic,
$V_{K}(\ell,0)$ is an irreducible $K$-module.  Because
$V(\widetilde{DY_{\hbar}(\sl_{2})},\ell)/\hbar
V(\widetilde{DY_{\hbar}(\sl_{2})},\ell)$ is a vacuum $K$-module of
level $\ell$ by Proposition \ref{plimit}, it follows that
$$V(\widetilde{DY_{\hbar}(\sl_{2})},\ell)/\hbar
V(\widetilde{DY_{\hbar}(\sl_{2})},\ell)\simeq V_{K}(\ell,0)$$ as a
$K$-module. It follows that this $K$-module isomorphism is also an
isomorphism of nonlocal vertex algebras. As $V_{K}(\ell,0)$ is
(irreducible) non-degenerate,
$V(\widetilde{DY_{\hbar}(\sl_{2})},\ell)$ is a non-degenerate
$\hbar$-adic quantum vertex algebra.
\end{proof}

\end{document}